# Embedded contact homology and Seiberg-Witten Floer cohomology I

Clifford Henry Taubes[†]

Department of Mathematics
Harvard University
Cambridge, MA  O2138

chtaubes@math.harvard.edu

This is the first of five papers that construct an isomorphism between the embedded contact homology and Seiberg-Witten Floer cohomology of a compact 3-manifold with a given contact 1-form.  This paper describes what is involved in the construction.

[†]Supported in part by the National Science Foundation

# 1. Introduction

The purpose of this article is to describe an isomorphism between the Seiberg-Witten Floer cohomology of a compact, oriented 3-manifold and the embedded contact homology as defined by a given contact 1-form on the 3-manifold. What follows momentarily is a very brief description of how these coholmology/homology groups are defined. A more detailed description is provided for both later in the article.

Consider first the Seiberg-Witten side of the story. Let M denote the 3-manifold in question. The Seiberg-Witten Floer homology/cohomology is defined with the choice of a $Spin_{\mathbb{C}}$ structure on M. With a Riemannian metric chosen, the latter is an equivalence class of lifts of the oriented, orthonormal frame bundle to a principle $Spin_{\mathbb{C}}$ bundle. Each $Spin_{\mathbb{C}}$ structure has an associated cohomology class in $H^2(M; \mathbb{Z})$, this its first Chern class. Let $p \in \{2, 4, \ldots\}$ denote the divisibility of this class. (It is always divisible by 2). Each $Spin_{\mathbb{C}}$ structure with non-torsion first Chern class defines $\mathbb{Z}/p\mathbb{Z}$ graded homology and cohomology groups. These are the associated Seiberg-Witten Floer homology and cohomology. These groups are always finitely generated. When the first Chern class is torsion, the associated Seiberg-Witten Floer homology and cohomology is $\mathbb{Z}$ graded. The Seiberg-Witten Floer homology in this case is finitely generated in each degree. As explained in [KM], the degrees in which the Seiberg-Witten Floer homology is non-zero are bounded from above, but never from below. This group is designated by $\widehat{HM}_*$ [KM]. There is a corresponding Seiberg-Witten Floer cohomology group as well; this designated by $\widehat{HM}^*$. In all cases, the generators of these homology and cohomology groups are the solutions to certain versions of the Seiberg-Witten equations on M; and the differentials are defined via a weighted count of certain sorts of solutions to the Seiberg-Witten equations on $\mathbb{R} \times M$. (The book [KM] is taken here to be the reference bible for the Seiberg-Witten side of the story.)

Consider next the contact homology story. Embedded contact homology was invented by Michael Hutchings (see [HS], [HT1]). The definition requires first the choice of a contact 1-form on M that is compatible with the orientation. Thus, the form, a, is such that $a \wedge da$ is nowhere zero and orients the 3-manifold. A (suitably generic) contact form of this sort defines a version of embedded contact homology and its associated cohomology for each $Spin_{\mathbb{C}}$ structure on M. These groups have $\mathbb{Z}/p\mathbb{Z}$ grading when the first Chern class is not torsion, and they are $\mathbb{Z}$ graded otherwise. The generators for these groups consist of finite sets where any given element is a pair that consists of a closed integral curve of the vector field that generates the kernel of da and a positive integer. Note that closed orbits with hyperbolic return map are paired only with the integer 1. The differential in each case is defined via a weighted count of certain embedded, pseudoholomorphic curves in the symplectization, $\mathbb{R} \times M$, of M.

The theorem that follows states formally what is said in the opening paragraph.

**Theorem 1**: *Let* M *denote a compact, oriented three dimensional manifold and let* a *denote a suitably generic contact 1-form on* M *that gives the chosen orientation. Fix a Spin$_\mathbb{C}$ structure on* M. *Then there is an isomorphism between the associated embedded contact homology and the Seiberg-Witten Floer cohomology that reverses the sign of the relative $\mathbb{Z}/p\mathbb{Z}$ grading.*

Note that this theorem implies two conjectures by Michael Hutchings: The embedded contact homology does not depend on the contact 1-form; and the embedded contact homology is finitely generated in each degree.

There are circumstances where both the Seiberg-Witten Floer cohomology and embedded contact homology, have additional structure. This additional structure is not discussed further until the fifth paper in this series, [T5], except for the remark that the isomorphism that is described here is compatible with these additional structures.

Theorem 1 can be viewed as the 3-manifold analog of the equivalence proved in [T1], [T2] between the Seiberg-Witten invariants of a compact symplectic manifold and certain of the Gromov invariants that are computed by counting its pseudoholomorphic curves. Theorem 1 can also be viewed as a generalization of [T3] and [T4] which use the existence of certain non-trivial Seiberg-Witten Floer homology classes to find closed orbits of the vector field that generates the kernel of da.

The proof of Theorem 1 occupies Section 4 of this article plus its three immediate sequels [T6], [T7] and [T8]. The proof uses many of the constructions and observations that are used in [T2]. These parts of the argument are summarized by Theorems 4.2 and 4.3 to come. The subsequent papers in this series contain the proofs of Theorems 4.2 and 4.3. Ideas from [T3] and [T4] also play a central role to the proof of Theorem 1. Theorems 4.4 and 4.5 contain most of the input to Theorem 1 from [T3] and [T4]. These last two theorems are proved below in Section 4h.

What follows is a table of contents for the remaining parts of this article.

<u>Section</u> 2: This section gives the definition of embedded contact homology. Also stated here is Proposition 2.5; this very useful proposition asserts that any given contact 1-form has a suitable deformation to one with properties that very much simplify the subsequent analysis.

<u>Section</u> 3: This section introduces the Seiberg-Witten equations and then the Seiberg-Witten Floer cohomology. It also describes the special versions of the Seiberg-Witten equations that can be defined with the help of a contact 1-form.

Section 4: This section proves Theorem 1 modulo two technical results, Theorems 4.2 and 4.3. The various parts of Theorems 4.2 and 4.3 are proved in the sequels, [T6], [T7] and [T8].

Section 5: Section 5 is meant to give a very rough picture of two key maps that are supplied by Theorems 4.2 and 4.3. Section 5 also indicates how the proof of Theorem 1 would proceed without the approximation result from Proposition 2.5.

Appendix: This contains the proof of Proposition 2.5.

Before continuing, the author hereby acknowledges the immense debt owed to Michael Hutchings, Peter Kronheimer and Tom Mrowka for sharing their thoughts and knowledge about the subject matter in this paper.

**2. Embedded contact homology**

The purpose of this section is to give the definition of embedded contact homology. As noted in the introduction, this homology theory was introduced by Michael Hutchings. Most of what follows here paraphrases parts of the accounts in [HS] and [HT1].

**a) Reeb orbits**

Use v in what follows to denote the vector field on M that generates the kernel of da and pairs with a so as to equal 1. It is a traditional to call v the *Reeb* vector field. A *Reeb orbit* denotes here an embedded circle with tangent v, thus a closed integral curve of v. A Reeb orbit is implicitly oriented by v.

Let $\gamma$ denote a Reeb orbit. The integral of the contact 1-form along $\gamma$ is denoted by $\ell_\gamma$. This integral, a positive number, is called the *symplectic action* of $\gamma$. The set of Reeb orbits whose symplectic action is bounded by any given positive number is a compact subset in $C^\infty(S^1; M)$.

Fix an almost complex structure, J, on the kernel of a so that $da(\cdot, J(\cdot))$ defines a Riemannian metric on the kernel of a. A Reeb orbit $\gamma$ has a neighborhood that is parametrized by the product of $S^1$ and a disk $D \subset \mathbb{C}$ about the origin by an embedding $\varphi$: $S^1 \times D \to M$ which makes a, da, and the Reeb vector field v appear as

- $\frac{2\pi}{\ell_\gamma} \varphi^* a = (1 - 2\nu|z|^2 - \mu \bar{z}^2 - \bar{\mu} z^2)\, dt + \frac{i}{2}(z\, d\bar{z} - \bar{z}\, dz) + \cdots$,
- $\frac{2\pi}{\ell_\gamma} da = i\, dz \wedge d\bar{z} - 2(\nu z + \mu \bar{z})\, d\bar{z} \wedge dt - 2(\nu \bar{z} + \bar{\mu} z)\, dz \wedge dt + \cdots$,
- $\frac{\ell_\gamma}{2\pi} v = \frac{\partial}{\partial t} + 2i(\nu z + \mu \bar{z})\frac{\partial}{\partial z} - 2i(\nu \bar{z} + \bar{\mu} z)\frac{\partial}{\partial \bar{z}} + \cdots$.

(2.1)

Here, ν and μ are respectively real and complex valued functions on $S^1$. The unwritten terms in the top equation are $\mathcal{O}(|z|^3)$ and those in the lower two equations are $\mathcal{O}(|z|^2)$. In (2.1) and in what follows, the circle $S^1$ is implicitly identified with $\mathbb{R}/(2\pi\mathbb{Z})$ and $t \in \mathbb{R}/(2\pi\mathbb{Z})$ is used to denote its affine coordinate. These coordinates are such that the vector field $\frac{\partial}{\partial z}$ at $z = 0$ pushes forward via φ so as to generate the +i eigenspace of J on kernel(a).

It follows as a consequence of (2.1) that the integral curves of v appear in this coordinate chart as the graphs of maps from an interval in $\mathbb{R}$ to D that obey an equation of the form

$$\tfrac{i}{2} \tfrac{d}{dt} z + \nu z + \mu \bar{z} = \mathfrak{r}$$

(2.2)

where $\mathfrak{r}$ is a smooth function of t and z with $|\mathfrak{r}| \leq c_0 |z|^2$ and $|d\mathfrak{r}| \leq c_0 |z|$.

The left hand side of (2.2) defines a first order, $\mathbb{R}$-linear symmetric operator on $C^\infty(\mathbb{R}; \mathbb{C})$, this the operator that takes a function $t \to z(t)$ to

$$Lz = \tfrac{i}{2} \tfrac{d}{dt} z + \nu z + \mu \bar{z}$$

(2.3)

Such an operator is defined given any pair $(\nu, \mu) \in C^\infty(S^1; \mathbb{R} \oplus \mathbb{C})$. When z is written in terms of real functions x and y as $z = x + iy$ and then any function in the kernel of (2.3) can be written as

$$\begin{pmatrix} x(t) \\ y(t) \end{pmatrix} = U \begin{pmatrix} x(0) \\ y(0) \end{pmatrix} \quad \textit{where } U|_t \in SL(2; \mathbb{R}) \textit{ for each } t \in \mathbb{R} \;.$$

(2.4)

As t varies in $[0, 2\pi]$, the map $t \to U|_t$ defines a path in SL(2: $\mathbb{R}$) from the identity. (The matrix $U|_{2\pi}$ is the linearization of the Reeb flow on kernel(a) along the Reeb orbit.)

A pair of functions (ν, μ) is said to be *non-degenerate* when the corresponding matrix U has trace($U|_{2\pi}$) ≠ 2. The pair is deemed to be *elliptic* when |trace($U|_{2\pi}$)| < 2 and *hyperbolic* when |trace($U|_{2\pi}$)| > 2. Note that when (ν, μ) is hyperbolic, then the k'th power of $U|_{2\pi}$ does not have eigenvalue 1 for any k. Such is the case because $U|_{2\pi}$ in this case has two real eigenvalues, one with absolute value greater than 1 and the other of the same sign with absolute value less than 1. When elliptic, the pair (ν, μ) is said to be *n-elliptic* when the k'th power of $U|_{2\pi}$ does not have eigenvalue 1 for all k ≤ n. Note that a matrix in SL(2; $\mathbb{R}$) whose trace has absolute value less than 2 has two complex eigenvalues, these are on the unit circle and one is the conjugate of the other. A Reeb orbit γ is said to be respectively *non-degenerate*, *hyperbolic*, or *n-elliptic* when such is the case for the functions (ν, μ) that come from (2.2). Note that the labeling of γ as either hyperbolic or n-elliptic is an intrinsic property of γ; it does not depend on the choice of φ or the almost complex structure on the kernel of a.

The notion of hyperbolic or n-elliptic can be viewed as a condition on the operator L in (2.3). In particular, L has trivial kernel on the space of maps from $S^1$ to $\mathbb{C}$ if and only if (ν, μ) is non-degenerate. If (ν, μ) is hyperbolic, then L has trivial kernel on the space of 2πk-periodic maps from $\mathbb{R}$ to $\mathbb{C}$ for any positive integer k. Meanwhile, if (ν, μ) is n-elliptic, then L has trivial kernel on the space of 2πk periodic maps from $\mathbb{R}$ to $\mathbb{C}$ for all k ∈ {1, 2, …, n}.

A complex valued function on $\mathbb{R}$ is said to be an eigenvector of L if L sends the function to a constant, real multiple of itself. The constant in question is called the eigenvalue. An eigenvector is said to be 2πk-periodic for a given integer k > 1 if it is 2πk-periodic but not 2πk´ periodic for any positive integer k´ < k. Two non-trivial eigenvectors can have the same eigenvalue only if they have the same periodicity and the same degree as a map from $S^1$ to $\mathbb{C}$−{0}. (A non-trivial eigenvector is nowhere zero.)

The definition of embedded contact homology requires the Reeb orbits to be non-degenerate, and the elliptic ones to be n-elliptic for all n. The following lemma asserts a well known fact that there exist such contact forms.

**Lemma 2.1**: *There exists a residual set of contact forms in* $C^\infty(M; T^*M)$ *whose associated Reeb orbits are either hyperbolic or else are n-elliptic for all positive integers n. In fact, given a positive integer n, and* L ≥ 1, *there is a dense open subset of contact forms in* $C^\infty(M; T^*M)$ *with the following property: If γ is an associated Reeb orbit with* $\ell_\gamma$ ≤ L, *then γ is hyperbolic or n-elliptic.*

A non-degenerate Reeb orbit is isolated in the following sense: There is an open, concentric disk D´ ⊂ D such that there are no Reeb orbits in φ($S^1$ × D´) except γ that generate the homology of φ($S^1$ × D´). Hyperbolic Reeb orbits have a stronger isolation property: There are no Reeb orbits except γ in φ($S^1$ × D´). In the case when γ is n-elliptic, then the following is true: But for multiple covers of γ, there are no Reeb orbits that generate the class of kγ in $H_1$(φ($S^1$ × D´); $\mathbb{Z}$) for any k ∈ {1, …, n}.

**b) Pseudoholomorphic subvarieties**

The manifold $\mathbb{R} \times M$ has a family of almost complex structures that it inherits from the contact geometry of M. These almost complex structures are characterized by the following properties: They are invariant with respect to translations along the $\mathbb{R}$ factor of $\mathbb{R} \times M$; they map the generator, $\frac{\partial}{\partial s}$, of these translations to ν; and they preserve the kernel of a. Such an almost complex structure endows M with a metric that sets the Hodge star of da equal to 2a and gives v norm 1. Almost complex structures of this sort are said to be *compatible with* a. These are the only ones used in this article.

Let J now denote an a-compatible, almost complex structure. An irreducible, pseudoholomorphic subvariety in $\mathbb{R} \times M$ is defined to be a closed subset with the following two properties: First, the complement of a finite set of points is a connected, 2-dimensional submanifold whose tangent space is J-invariant. Second, the integral of da over this submanifold is finite. A pseudoholomorphic subvariety is defined to be a finite union of irreducible, pseudoholomorphic subvarieties. What follows describes some of the salient properties of irreducible, pseudoholomorphic subvarieties $\mathbb{R} \times M$. The basic story on such curves is presented in a series of seminal papers of Hofer (see [Ho1]-[Ho4] and Hofer, Wysocki and Zehnder [HWZ1]-[HWZ3]).

Agree to use s in what follows to denote the Euclidean coordinate on the $\mathbb{R}$ factor of $\mathbb{R} \times M$. The first point to note is that these subvarieties are well behaved where |s| is large. To say more, suppose for the moment that the contact form comes from the residual subset that is described in Lemma 2.1. With an almost complex structure fixed, let $\Sigma$ denote a given pseudoholomorphic subvariety. Then there exists $s_0 > 1$ such that the $|s| \geq s_0$ portion of $\Sigma$ is a disjoint union of properly embedded cylinders to which the function s restricts without critical points. Each such cylinder is called an *end* of $\Sigma$. The ends on which $s \geq s_0$ are called the positive side ends, and those where $s \leq -s_0$ are called the negative side ends.

The subvariety $\Sigma$ may contain irreducible components of the form $\mathbb{R} \times \gamma$ with $\gamma$ a Reeb orbit. Such cylinders are the only $\mathbb{R}$-invariant pseudoholomorphic subvarieties where $\mathbb{R}$ is understood to act on $\mathbb{R} \times M$ as the constant translations along the $\mathbb{R}$ factor.

Let $\mathcal{E}$ denote an end of $\Sigma$ that is not in an $\mathbb{R}$-invariant cylinder. Then there is a Reeb orbit, $\gamma = \gamma_\mathcal{E}$, and a positive integer $q_\mathcal{E}$ such that the following is true: Each constant s slice of $\mathcal{E}$, thus $\mathcal{E}|_s \subset M$, is a braid in the $S^1 \times D$ tubular neighborhood of $\gamma$ that projects as a $q_\mathcal{E}$ to 1 covering map to the central circle. Moreover, as $|s| \to \infty$, these braids converge pointwise to $\gamma$.

To say more about this, note that a tubular neighborhood map $\varphi: S^1 \times D \to M$ for $\gamma$ can be chosen so as to have the following additional properties: First, the vector field $\frac{\partial}{\partial z}$ along D pushes forward to define a type (1,0) tangent vector along $\mathbb{R} \times \gamma$ in $\mathbb{R} \times M$ with length $2^{-1/2}$. Second, the C-valued 1-form

$$dz - 2i(\nu z + \mu \bar{z})dt$$

(2.5)

differs from φ-pull back of a 1-form in $T^{1,0}(\mathbb{R} \times M)$ by $\mathcal{O}(|z|^2)\,dt$, $\mathcal{O}(|z|^2)\,dz$ and $\mathcal{O}(|z|)\,d\bar{z}$. Note that it follows from the second line in (2.1) that this form has length $(\frac{4\pi}{\ell_\gamma})^{1/2}$ to $\mathcal{O}(|z|^2)$. Henceforth, all tubular neighborhood maps are assumed to be of this sort.

Now, suppose that $\mathcal{E}$ denotes a $s \ll -1$ end of $\Sigma$ that is not part of an $\mathbb{R}$-invariant cylinder. Let $\gamma = \gamma_\mathcal{E}$ and $q_\mathcal{E}$ be as described above. The end $\mathcal{E}$ can be viewed using the tubular neighborhood map as a subvariety in $\mathbb{R} \times (S^1 \times D)$, this the image of a map from $(-\infty, -s_0] \times \mathbb{R}/(2\pi q_\mathcal{E}\mathbb{Z})$ into $\mathbb{R} \times (S^1 \times D)$ that sends any given $(s, t)$ to the point $(s, t, z(s, t))$ where $z$ is a certain $\mathbb{C}$-valued function. To say more about $z(\cdot)$, introduce the set $\text{div}_\mathcal{E} \subset \{1, 2, \ldots, q\}$ with the following two properties: First, any given $q' \in \text{div}_\mathcal{E}$ evenly divides $q$. Second, there is a $2\pi q'$ periodic eigenvector, $\varsigma_{q'}$, of (2.3) with eigenvalue $\lambda_{q'}$ such that $0 > \lambda_{q'} \geq \lambda_{q_\mathcal{E}}$. The function $z(\cdot)$ is given in terms of this data as

$$z(s, t) = \sum_{q' \in \text{div}_\mathcal{E}} (\varsigma_{q'}(t) + \mathfrak{r}_{q'})\, e^{-2\lambda_{q'} s} ,$$

(2.6)

where $\mathfrak{r}_{q'}$ is $2\pi q'$ periodic and its norms and that of its derivative are bounded by $e^{-\varepsilon|s|}$ with $\varepsilon$ a positive constant. Note in this regard that $\varsigma_{q_\mathcal{E}} \neq 0$. A positive end of $\Sigma$ is described by (2.6) but with each $\lambda_{q_\mathcal{E}} \geq \lambda_{q'} > 0$. Equation (2.6) can be derived from what is done, for example, in [S].

The set, $\mathcal{M}$, of irreducible, pseudoholomorphic subvarieties has nice properties also. To say more, endow $\mathcal{M}$ with the topology whereby the neighborhoods of a given element $\Sigma \in \mathcal{M}$ are generated by sets of the following sort: A given basis set is labeled by $\varepsilon > 0$ and it consists of the subvarieties $\Sigma' \subset \mathcal{M}$ with the following two properties:

- $\sup_{z \in \Sigma} \text{dist}(z, \Sigma') + \sup_{z' \in \Sigma'} \text{dist}(\Sigma, z') < \varepsilon$ .
- *If $\varpi$ is a compactly supported 2-form on $\mathbb{R} \times M$, then* $|\int_\Sigma \varpi - \int_{\Sigma'} \varpi| \leq \varepsilon \sup_{\mathbb{R} \in M} |\varpi|$.

(2.7)

What follows is the basic structure theorem for $\mathcal{M}$.

**Lemma 2.2**: *Fix $\Sigma \subset \mathcal{M}$ with the following property: Let $\mathcal{E}$ be any given end of $\Sigma$ and let $\gamma$ denote the Reeb orbit that is approached by the $|s| \to \infty$ limit of the constant s slices of $\mathcal{E}$. Then $\gamma$ is either hyperbolic or $q_\mathcal{E}$-elliptic. Assuming that $\Sigma$ has these properties, there exists a Fredholm operator, $\mathcal{D}_\Sigma$, a ball $B \subset \text{kernel}(\mathcal{D}_\Sigma)$ a smooth map $f: B \to \text{cokernel}(\mathcal{D}_\Sigma)$ and a homeomorphism from $f^{-1}(0)$ to a neighborhood of $\Sigma$ in $\mathcal{M}$. Here, $f(0) = 0$ and the homeomorphism sends 0 to $\Sigma$. Furthermore:*
- *Let $\mathcal{M}_{\text{reg}} \subset \mathcal{M}$ denote the set that consists of those $\Sigma$ with $\text{cokernel}(\mathcal{D}_\Sigma) = 0$. This set $\mathcal{M}_{\text{reg}}$ is open and it has the structure of a smooth manifold*

- *Let $\Sigma \in \mathcal{M}_{reg}$. Then the just described homeomorphism from $B \subset \text{kernel}(\mathcal{D}_\Sigma)$ into $\mathcal{M}$ gives a smooth coordinate chart for a neighborhood of $\Sigma$.*
- *If the contact form comes from Lemma 2.1's residual set, then there is a residual set of compatible almost complex structures for which $\mathcal{M}_{reg} = \mathcal{M}$.*

To say more about $\mathcal{D}_\Sigma$, remember that an irreducible, pseudoholomorphic subvariety, C, has a *model curve*; this a complex curve, $C_0$, together with an almost everywhere 1-1 pseudoholomorphic map $\phi: C_0 \to \mathbb{R} \times M$ whose image is C. Assuming that C has only immersion singularities, there is a well defined pull-back normal bundle over $C_0$. This is the bundle, $N \to C_0$, whose fiber at any given point is the normal 2-plane in $T(\mathbb{R} \times M)$ at the point's $\phi$-image to an embedded disk in C. The composition of an exponential map with a section of a suitable disk sub-bundle of N defines a deformation of C in $\mathbb{R} \times M$. Here, an exponential map is a smooth map from a uniform radius disk subbundle in N to $\mathbb{R} \times M$ that restricts as $\phi$ to the zero section and has surjective differential along the zero section.

The almost complex structure gives N a complex bundle structure, and the induced metric from $\mathbb{R} \times M$ gives N a hermitian structure and thus the structure of a holomorphic line bundle. As such, there is an associated d-bar operator that maps sections of N to those of $N \otimes T^{1,0}C_0$. This d-bar operator enters the story in the following manner: A deformation of C that preserves to first order the pseudoholomorphic condition is the image via an exponential map of a section of N that is annihilated by an operator that can be viewed as the d-bar operator with an extra, $\mathbb{R}$-linear, zero'th order term. It can be identified as operator that sends a section $\zeta$ to

$$\mathcal{D}_C \zeta = \bar{\partial}\zeta + \nu_C \zeta + \mu_C \bar{\zeta} ,$$

(2.8)

where $\nu_C$ is a section of $T^{0,1}C_0$ and where $\mu_C$ is one of $N^2 \otimes T^{0,1}C_0$. Note that the parametrization given in (2.6) for any given end of C induces a trivialization of N and $TC_0$ on such an end with the following property: When written using this trivialization, the pair $(\nu_C, \mu_C)$ converges as $|s| \to \infty$ on the end to the pair $(\nu, \mu)$ that appears in the associated version of (2.1).

The operator $\mathcal{D}_C$ defines a bounded, $\mathbb{R}$-linear Fredolm map from the Sobolev space $L^2_1(C_0; N)$ to $L^2(C_0; N \otimes T^{0,1}C_0)$ if the following is true: The constant s slices of each end of C limit as $|s| \to \infty$ as some integer q-fold cover of a Reeb orbit that is either hyperbolic or q-elliptic.

A more complicated version of this $\mathbb{R}$-linear operator defines $\mathcal{D}_C$ in the cases when C has non-immersion singularities. The operator $\mathcal{D}_\Sigma$ is $\mathcal{D}_C$ if $\Sigma = C$. If each irreducible component of $\Sigma$ is an immersed curve of the sort just described, then $\mathcal{D}_\Sigma$ is the direct sum of the corresponding operators with domain and range the direct sum of the corresponding Sobolev spaces.

**c) Embedded contact homology**

The following definition of Hutching's embedded contact homology is taken from Section 11 in [HS]. To set the stage for the definition, note that each class in $H_1(M; \mathbb{Z})$ labels a version of this homology. This understood, fix a class $\Gamma$. Assume that the contact form a is from Lemma 2.1's residual set.

The Chain complex: The chain complex for $\Gamma$'s version of embedded contact homology is the free $\mathbb{Z}$ module that is generated by equivalance classes of pairs $(\Theta, \mathfrak{o})$ where $\Theta$ and $\mathfrak{o}$ are as follows. First $\Theta$ is a finite set of pairs of the form $(\gamma, m)$ where $\gamma$ is a Reeb orbit and m is a positive integer subject to three constraints. First, no two pairs have the same Reeb orbit. Second, $m = 1$ when $\gamma$ is hyperbolic. Third, the formal sum $[\Theta] = \sum_{(\gamma, m) \in \Theta} m \gamma$ should define a closed cycle that generates the class $\Gamma$ in $H_1(M; \mathbb{Z})$. Note that the empty set $\Theta = \emptyset$ defines a generator in the case when $\Gamma$ is the trivial class.

Meanwhile, $\mathfrak{o}$ is an ordering of the pairs in $\Theta$ whose Reeb orbit component is hyperbolic and whose version of (2.4) has matrix U with $\text{trace}(U|_{2\pi}) > 2$. The equivalence relation identifies $(\Theta, \mathfrak{o})$ with $\pm(\Theta, \mathfrak{o}')$ where the $\pm$ factor is the image in $\{\pm 1\}$ of the permutation that takes $\mathfrak{o}$ to $\mathfrak{o}'$.

The free $\mathbb{Z}$ module so generated is denoted by $\mathcal{C}_{\text{ech}}$ in what follows. Keep in mind that it depends on $\Gamma$. In most of what follows, the pair $(\Theta, \mathfrak{o})$ will be denoted as $\Theta$ with the presence of the ordering $\mathfrak{o}$ implicit.

The grading: Let $K^{`} \subset TM$ denote the 2-plane bundle given by the kernel of a. Orient this bundle using da, and let $-\hat{e}_K \in H^2(M; \mathbb{Z})$ denote its Euler class. Let $P(\Gamma) \in H^2(M; \mathbb{Z})$ denote the Poincaré dual of $\Gamma$ and let p denote the divisibility of $-\hat{e}_K + 2P(\Gamma)$. The $\mathbb{Z}$ module $\mathcal{C}_{\text{ech}}$ has a relative $\mathbb{Z}/p\mathbb{Z}$ grading whose definition is given in the five steps that follow.

*Step 1*: The path U: $[0, 2\pi] \to SL(2; \mathbb{R})$ can be used to assign a *rotation number* to any hyperbolic or n-elliptic path $(\nu, \mu)$. This rotation number is defined as follows: When the pair $(\nu, \mu)$ is hyperbolic, there is a homotopy of U through a 1-parameter family of paths such that the $t = 2\pi$ element of each path on this family has $|\text{trace}(\cdot)| \geq 2$, and so that the end member path is a path of pure rotations. As such, the end member path rotates $\mathbb{C}$ by a total of $\pi k$ radians for some $k \in \mathbb{Z}$. This integer k is the rotation number for $(\nu, \mu)$. In the case when $(\nu, \mu)$ is n-elliptic, there is a similar homotopy of U, now through a family of paths such that each path in this family has its $t = 2\pi$ point conjugate to U(2π). The end-member of this homotopy is again a path of rotations, this time rotating $\mathbb{C}$ by an angle $2\pi R$. This angle R is the rotation number. Note in this regard that the n-elliptic condition means that kR is not in $\mathbb{Z}$ when k ∈ {1, …, n}. The numbers k or R depend on φ, but not k mod(2) or R mod($\mathbb{Z}$).

Suppose that the pair (ν, μ) is defined from a given Reeb orbit γ and coordinate map φ. Use $z_{\gamma,1}$ to denote the rotation number k ∈ $\mathbb{Z}$ when γ is hyperbolic. When γ is m-elliptic and q ∈ {1, …, m}, use $z_{\gamma,q}$ to denote 1 plus twice the greatest integer less qR.

*Step 2*: Let $\Theta_-$ and $\Theta_+$ define generators in $\mathcal{C}_{\text{ech}}$. Fix a tubular neighborhood embedding as described above for each Reeb orbit from the pairs that comprise $\Theta_-$ and $\Theta_+$. As $\Theta_-$ and $\Theta_+$ define the class Γ in the manner just described, there is a smooth, oriented, properly immersed surface Z ⊂ $\mathbb{R} \times M$ with transversal self-intersections that has the following properties: The |s| >> 1 portion of this surface is a disjoint union of embedded cylinders on which s restricts as a function with no critical points. The cylinders that sit where s >> 1 are distinguished in part by the elements in $\Theta_+$. In particular, a given pair (γ, m) labels m such cylinders. Given γ's tubular neighborhood map φ, then each of the m cylinders sits in $\mathbb{R} \times \varphi(S^1 \times D)$ as the image of the graph over $\mathbb{R} \times S^1$ of the function that sends (s, t) to $e^{-2\lambda s + ix} \in \mathbb{C}$ where λ > 0 and $x \in \mathbb{R}/(2\pi\mathbb{Z})$. However, if C and C´ are two such cylinders, then the corresponding points $x$ and $x´$ must define distinct points in the circle. There is an analogous correspondence between the cylinders that sit where s << -1 with the elements in $\Theta_-$. The only difference is that λ is now required to be negative.

*Step 3*: Let Z denote a surface as described in Step 2. Then Z has a relative self-intersection pairing which defined to be its intersection number with a deformation, Z´, whose restriction to any given |s| >> 1 cylinder deforms the latter so as to change the parameter $x$ to $x + \varepsilon$ with ε > 0 but very small. With Z denoting the surface in question, $Q_Z$ is used to denote this self-intersection number.

*Step 4*: The surface Z also has a well defined pairing with the Euler class of the bundle $K^{-1}$. This pairing is defined by the usual count of the zeros of a section of this bundle over the surface with the proviso that the section should restrict to each |s| >> 1 cylinder so as to be non-zero, and to be constant with respect to the trivialization of K on $\varphi(S^1 \times D)$ that is given by the coordinate vector field $\frac{\partial}{\partial z}$ on D. Note in this regard that K and the tangent space to D agree along $S^1 \times \{0\}$. This pairing is denoted here by $-\langle c_1, Z \rangle$.

*Step 5*: With Σ chosen as in Step 2, introduce the integer

$$I(\Theta_-, \Theta_+; Z) = -\langle c_1, Z \rangle + Q_Z + \sum_{(\gamma,m)\in\Theta_+} \sum_{1\leq q\leq m} z_{\gamma,q} - \sum_{(\gamma,m)\in\Theta_-} \sum_{1\leq q\leq m} z_{\gamma,q}.$$
(2.9)

Although the various tubular neighborhood embeddings of $S^1 \times D$ are needed to make sense of the terms in (2.9), the value of $I(\Theta_-, \Theta_+; Z)$ does not, in fact, depend on them. Moreover, the image of I in $\mathbb{Z}/p\mathbb{Z}$ depends only on the ordered pair $(\Theta_-, \Theta_+)$. This is proved by Hutchings in [Hu]; see also [HS]. Hutchings also proves that this image in $\mathbb{Z}/(p\mathbb{Z})$ obeys the sum rule $I(\Theta_1, \Theta_2) + I(\Theta_2, \Theta_3) = I(\Theta_1, \Theta_3)$. The relative $\mathbb{Z}/p\mathbb{Z}$ degree assignments to the generators of $\mathcal{C}_{ech}$ are made so that $I(\Theta_-, \Theta_+) = \text{degree}(\Theta_-) - \text{degree}(\Theta_+)$.

The differential: The differential that is used by Hutchings decreases the $\mathbb{Z}/p\mathbb{Z}$ grading by 1. Its definition is given in the three steps that follow.

*Step 1*: The almost complex structure should be chosen so as to be generic in the sense given in Lemma 2.2. Suppose that $\Theta_-$ and $\Theta_+$ are generators of $\mathcal{C}_{ech}$. Introduce $\mathcal{M}_1(\Theta_-, \Theta_+)$ to denote the set whose elements are finite sets of pairs of the form (C, m) where C is a pseudoholomorphic subvariety and m is a positive integer. The elements in this set are constrained as follows: First, m = 1 unless $C = \mathbb{R} \times \gamma$ with $\gamma$ a closed Reeb orbit. Second, if (C, m) and (C´, m´) are distinct pairs, then C is not a translate of C´ along the $\mathbb{R}$ factor of $\mathbb{R} \times M$. To state the third constraint, let $\pi_M$ denote the projection from $\mathbb{R} \times M$ to M. Here is the third constaint: A given element $\Sigma \in \mathcal{M}_1(\Theta_-, \Theta_+)$ defines the formal sum $\sum_{(C,m)\in\Sigma} m \pi_M(C)$, here viewed as a 2-cycle. The boundary of this 2-cycle must be $\sum_{(\gamma,m)\in\Theta_+} m\gamma - \sum_{(\gamma,m)\in\Theta_-} m\gamma$. To state the fourth constraint, let $H_2(M, \Theta_-, \Theta_+)$ denote the set of that consists of the relative homology classes of 2-chains $z \subset M$ with $\partial z = \sum_{(\gamma,m)\in\Theta_+} m\gamma - \sum_{(\gamma,m)\in\Theta_-} m\gamma$. To be more explicit, chains z and z´ define the same class in $H_2(M, \Theta_-, \Theta_+)$ when the closed cycle z - z´ is the boundary of a 3-cycle in M. Thus, $H_2(M, \Theta_-, \Theta_+)$ is an affine space modeled on $H_2(M; \mathbb{Z})$. Let $Z \subset \mathbb{R} \times M$ denote a surface as described above that can be used to define the invariant $I(\Theta_-, \Theta_+, \cdot)$, but one such that $\pi_M(Z)$ and $\sum_{(C,m)\in\Sigma} m\pi_M(C)$ define the same element in $H_2(M, \Theta_-, \Theta_+)$. Here is the fourth constraint: $I(\Theta_-, \Theta_+, Z) = 1$. Note in this regard that $I(\Theta_-, \Theta_+, \cdot)$ takes identical values on surfaces Z and Z´ that give the same element in $H_2(M, \Theta_-, \Theta_+)$. Take $\mathcal{M}_1(\Theta_-, \Theta_+) = \emptyset$ when $I(\Theta_-, \Theta_+) \neq 1$.

The set $\mathcal{M}_1(\Theta_-, \Theta_+)$ inherits a topology and a local structure of the sort described by Lemma 2.2 using its tautological embedding into a disjoint union of products of $\mathcal{M}$.

*Step 2*: As noted by Hutchings, there exists a residual set of almost complex structures for which the resulting version of $\mathcal{M}_1(\Theta_-, \Theta_+)$ has the following properties:

1. *If $\Sigma \in \mathcal{M}_1(\Theta_-, \Theta_+)$, then $\cup_{(C,m)\in\Sigma} C$ is an embedded, pseudoholomorphic subvariety.*

2. *The space $\mathcal{M}_1(\Theta_-, \Theta_+)$ has a finite set of components and each component is a smooth, 1-dimensional manifold.*

3. *The $\mathbb{R}$ action on $\mathbb{R} \times M$ induces a free $\mathbb{R}$ action on each component. As a consequence, each element in $\mathcal{M}_1(\Theta_-, \Theta_+)$ consists of a disjoint union of $\mathbb{R}$-invariant cylinders with integer weights and one pseudoholomorphic submanifold that is not $\mathbb{R}$-invariant.*

4. *Let $\Sigma \subset \mathcal{M}_1(\Theta_-, \Theta_+)$ and let $\mathcal{E} \subset \Sigma$ denote an end. Let $q_\mathcal{E}$ and $\mathrm{div}_\mathcal{E}$ denote the data that appear in $\mathcal{E}$'s version of (2.6). Then $\mathrm{div}_\mathcal{E} = \{q_\mathcal{E}\}$.*

5. *Let $\Sigma \subset \mathcal{M}_1(\Theta_-, \Theta_+)$ and let $\mathcal{E}$ and $\mathcal{E}'$ denote distinct pairs of either positive or negative ends of $\Sigma$ with $\gamma_\mathcal{E} = \gamma_{\mathcal{E}'}$ and $q_\mathcal{E} = q_{\mathcal{E}'}$. Let $\varsigma_{q_\mathcal{E}}$ and $\varsigma_{q_{\mathcal{E}'}}$ denote the $2\pi q$-periodic eigenvector that appears in the respective $\mathcal{E}$ and $\mathcal{E}'$ versions of (2.6). Then $\varsigma_{q_\mathcal{E}}|_t \neq \varsigma_{q_{\mathcal{E}'}}|_{t+2\pi k}$ for any $t \in S^1$ and $k \in \mathbb{Z}$.*

(2.10)

Properties 1-3 are proved in [Hu]; see also [HS]. Properties 4 and 5 follow from what is done in Section 3 of [HT2]. Property 4 also needs some facts that can be derived from what is said in Section 11 of [Hu] and its Remark 1.

Let $\mathcal{J}_a'$ denote the residual set of a-compatible almost complex structures that have the properties listed in (2.10) and lie in Lemma 2.2's residual set. Fix an almost complex structure from $\mathcal{J}_a'$. Let $\Theta_-$ and $\Theta_+$ denote generators of the embedded contact homology chain complex. Hutching uses constructions of Bourgeois and Mohnke [BM] to associate a sign, $\pm 1$, to each component of $\mathcal{M}_1(\Theta_-, \Theta_+)$. The full details of this are given in Section 9 and especially Section 9.5 of [HT2]. Let $\sigma(\Theta_-, \Theta_+)$ denote the sum of these signs when $\mathcal{M}_1(\Theta_-, \Theta_+)$ is non-empty, and 0 otherwise.

To say a wee bit more about these signs, note that the sign that is associated to any given component of $\mathcal{M}_1(\Theta_-, \Theta_+)$ is obtained by comparing two natural orientations. The first is that induced by the generator of the $\mathbb{R}$ action. The second is defined using ideas of Quillen [Q] about determinant line bundles of parametrized families of Fredholm operators. As noted by Hutchings and explained in Section 9.5 of [HT2], these ideas of Quillen can be used along lines explained in [BM] so as to define a second orienation to each component of $\mathcal{M}_1(\Theta_-, \Theta_+)$. (The respective parts of $\Theta_-$ and $\Theta_+$ that involve the ordering of their even rotation number hyperbolic Reeb orbits is needed solely to define this second orientation of $\mathcal{M}_1(\Theta_-, \Theta_+)$.)

*Step 3*: The differential, $\delta$, on $\mathcal{C}_{\mathrm{ech}}$ is defined on any given generator $\Theta$ by the formula $\delta\Theta = \sum_{\Theta' \in \mathcal{C}_{\mathrm{ech}}} \sigma(\Theta', \Theta) \Theta'$. The proof that $\delta^2 = 0$ appears in [HT1] and [HT2]. The proof required that the almost complex structure come from a certain residual subset of $\mathcal{J}_a'$. This last residual set is denoted by $\mathcal{J}_a$ in what follows. Almost complex structures will always be chosen from $\mathcal{J}_a$ unless explicitly noted otherwise.

The embedded contact homology for the class $\Gamma$ is defined to be the homology of $\delta$ on the $\mathcal{C}_{ech}$. As defined, $\delta$ decreases the $\mathbb{Z}/p\mathbb{Z}$ degree by 1 so this homology is $\mathbb{Z}/p\mathbb{Z}$ graded. The homology defined here is denoted in what follows by $\mathcal{ECH}$ where it is understood that the class $\Gamma$ is fixed in advance and not subsequently changed.

The Filtration: As noted in Section 2a, the set of Reeb orbits with an apriori symplectic action bound is compact. If all such orbits are non-degenerate, then there are but a finite set with symplectic action less than any given amount. Granted this last point, fix $L > 0$ and let $\mathcal{C}_{ech}^L \subset \mathcal{C}_{ech}$ denote the submodule that is generated by elements $\Theta$ that obey $\sum_{(\gamma,m) \in \Theta} m \ell_\gamma < L$. This is a finitely generated chain complex.

If $\Sigma \subset \mathcal{M}_1(\Theta, \Theta')$, then $\sum_{(\gamma,m) \in \Theta} m \ell_\gamma \leq \sum_{(\gamma,m) \in \Theta'} m \ell_\gamma$. As a consequence, the differential on $\mathcal{C}_{ech}$ maps $\mathcal{C}_{ech}^L$ to itself. This understood, let $\mathcal{ECH}^L$ denote the homology that is defined by $\delta$ on $\mathcal{C}_{ech}^L$. Then $\mathcal{ECH} = \mathrm{dir\,lim}_{L \to \infty} \mathcal{ECH}^L$ where the homomorphisms for this direct limit are induced by the L and $L' > L$ versions of the submodule inclusion homomorphism from $\mathcal{C}_{ech}^L$ into $\mathcal{C}_{ech}^{L'}$.

**d) Changing the contact structure**

As it turns out, the proof of Theorem 1 is considerably shorter when the contact structure is approximated by one which has a canonical form near some of its Reeb orbits. The following lemma describes these canonical forms.

**Lemma 2.3**: *Suppose that $(\nu, \mu) \in C^\infty(S^1; \mathbb{R} \oplus \mathbb{C})$.*
The elliptic case: *Suppose that $(\nu, \mu)$ is elliptic with rotation angle $R \in \mathbb{R}$. There is a homotopy of $(\nu, \mu)$ through elliptic pairs with rotation angle R to the pair $(\frac{1}{2} R, 0)$.*
The hyperbolic case: *Suppose that $(\nu, \mu)$ is hyperbolic with rotation number k. If $\varepsilon > 0$ is small, there is a homotopy of $(\nu, \mu)$ through hyperbolic pairs to the pair $(\frac{1}{4} k, i\varepsilon\, e^{ikt})$.*

*Proof of Lemma 2.3*: The statement in the elliptic case is straightforward; it follows readily from the geometry of $SL(2; \mathbb{R})$ that any two elliptic pairs with the same rotation number are homotopic through a family of constant rotation number elliptic pairs. In the hyperbolic case, remark first that any two hyperbolic pairs with the same rotation number are homotopic through hyperbolic pairs. As a consequence, it is enough to verify that the pair $(\frac{1}{2} k, i\varepsilon\, e^{-ikt})$ is hyperbolic with rotation number k when $\varepsilon$ is sufficiently small. The calculation is straight forward and left to the reader.

Assume that the contact structure a is from Lemma 2.1's residual set and the almost complex structure, J, is from $\mathcal{J}_a$. Fix $L \geq 1$ and $\delta > 0$. A pair $(\hat{a}, \hat{J})$ of contact structure on M and compatible almost complex structure on $\mathbb{R} \times M$ is said to be a $(\delta, L)$-*approximation* for $(a, J)$ when the following is true: There is a smooth, 1-parameter family $\{(a_\tau, J_\tau)\}_{\tau \in [0,1]}$ of pairs of contact structure and compatible almost complex structure with $(a_0, J_0) = (a, J)$ and $(a_1, J_1) = (\hat{a}, \hat{J})$; and such that

1. *For each* $\tau \in [0, 1]$, *the respective sets of* a *and* $a_\tau$ *Reeb orbits with symplectic action less than* L *are identical.*
2. *Let* γ *denote a Reeb orbit for* a *with* $\ell_\gamma < L$. *If* γ *is elliptic or hyperbolic as defined using* a, *then it is respectively elliptic or hyperbolic as defined using any* $\tau \in [0, 1]$ *version of* $a_\tau$.
3. *Let* $\Theta_-$ *and* $\Theta_+$ *denote generators of* $\mathcal{C}_{ech}{}^L$. *For each* $\tau \in [0, 1]$, *there is a 1-1 correspondence between the components of the respective J and* $J_\tau$ *versions of the space* $\mathcal{M}_1(\Theta_-, \Theta_+)$ *such that partnered components contribute the same sign to the respective J and* $J_\tau$ *versions of* $\sigma(\Theta_-, \Theta_+)$.
4. *Let* γ *denote a Reeb orbit with* $\ell_\gamma < L$. *There is a coordinate embedding* $\varphi: S^1 \times D \to$ M *of the sort described in the preceding with the following property*: *If* γ *is hyperbolic with rotation number* k, *then the* $\hat{a}$*-version of the pair* $(\nu, \mu)$ *is equal to* $(\frac{1}{2} k, i\varepsilon e^{-ikt})$ *for some* $\varepsilon \in (0, \delta)$. *If* γ *is elliptic with rotation number* R, *then*

    i) $\frac{2\pi}{\ell_\gamma} \varphi^*\hat{a} = (1 - R|z|^2)dt + \frac{i}{2}(z d\bar{z} - \bar{z} dz)$.

    ii) *The* $\varphi^*$-*pull back of the* $\hat{J}$-*version of* $T^{1,0}(\mathbb{R} \times M)$ *is spanned by the forms* $ds + i\hat{a}$ *and* $\frac{\ell_\gamma}{2\pi}(dz - i R z dt)$.

5. *The contact structure* $\hat{a}$ *comes from Lemma 2.1's residual set and the almost complex structure* $\hat{J}$ *comes from the set* $\mathcal{J}_{\hat{a}}$.

(2.11)

The next subsection gives a first indication as to why pairs $(\hat{a}, \hat{J})$ as just described are easy to work with.

The proposition that follows asserts that the homology of $\mathcal{C}_*^L$ as defined by a pair $(a, J)$ is isomorphic to that defined by a $(\delta, L)$ approximation.

**Proposition 2.4**: *Let* a *denote a contact 1-form from the residual set given in Lemma 2.1 and let J denote a complex structure from* $\mathcal{J}_a$. *Fix* $\delta > 0$ *and* $L \geq 1$ *such that there is no generator* $\Theta$ *of* $\mathcal{C}_{ech}$ *with* $\sum_{(\gamma,m) \in \Theta} m \ell_\gamma = L$. *Let* $(\hat{a}, \hat{J})$ *denote a* $(\delta, L)$ *approximation to the given pair* $(a, J)$. *Then the identification provided by the first item in (2.11) between the Reeb orbits with symplectic action less than* L *induces a degree preserving isomorphism between the* a *and* $\hat{a}$ *versions of* $\mathcal{C}_{ech}{}^L$ *that intertwines the respective* *differentials. Thus, it induces an isomorphism between the respective* (a, J) *and* (â, Ĵ) *versions of* $\mathcal{ECH}^L$.

*Proof of Proposition 2.4*: The fact that the isomorphism preserves degree follows from the third item in (2.11). The fact that it intertwines the differential follows from the fourth item in (2.11).

The final proposition asserts that there are in all cases (δ, L) approximations.

**Proposition 2.5**: *Let* a *denote a contact form from Lemma 2.1's residual set and let* J ∈ $\mathcal{J}_a$. *Fix* δ > 0 *and* L ≥ 1 *such that there is no generator* Θ ∈ $\mathcal{C}_{ech}$ *with* $\sum_{(\gamma,m)\in\Theta} m\,\ell_\gamma = L$. *Then there exist* (δ, L) *approximations to* (a, J).

This proposition is proved in the Appendix to this article.

### e) Pseudoholomorphic subvarieties for (δ, L) approximating pairs

This last subsection is an aside of sorts whose purpose is to say something about the pseudoholomorphic curves for a pair (a, J) of contact form on $S^1 \times \mathbb{C}$ and compatible almost complex structure on $\mathbb{R} \times S^1 \times \mathbb{C}$ where

$$a = \tfrac{\ell}{2\pi}((1 - R|z|^2)\,dt + \tfrac{i}{2}(z\,d\bar{z} - \bar{z}\,dz)),$$

(2.12)

and where $T^{1,0}(\mathbb{R} \times S^1 \times \mathbb{C})$ is spanned by ds + ia and dz – iR zdt. Here, $\ell > 0$ and R are constant. A straightforward calculation verifies the following: Let $w = \tfrac{2\pi}{\ell}s - \tfrac{1}{2}|z|^2$. Then

- dw + i dt ∈ $T^{1,0}(\mathbb{R} \times M)$; *thus the constant* (w, t) *planes are pseudoholomorphic*.
- *Subvarieties* z = f(w, t) *are pseudoholomophic if and only if* $(\tfrac{\partial}{\partial w} + i\tfrac{\partial}{\partial t} + R)f = 0$.

(2.13)

Note that the complex structure in this case is integrable. Local holomorphic corrdinates are u = w + it and x = $e^{Rw}$ z.

## 3. Seiberg-Witten Floer (co)homology

The purpose of this section is to say more about the relevant versions of Seiberg-Witten Floer homology and cohomology. As a complete treatment of the subject is given by Kronheimer and Mrowka in [KM], what follows focuses for the most part on those aspects of the story that are relevant to the case when M comes with a contact 1-form. In any event, much of what is said below paraphrases the definitions and discussion in Kronheimer and Mrowka's book [KM].

**a) The Seiberg-Witten equations on M and $\mathbb{R} \times M$.**

Fix a Riemannian metric on M so as to define the bundle of oriented, orthonormal frames for TM. Let Fr $\to$ M denote this principle SO(3) bundle. A $\text{Spin}_\mathbb{C}$ lift of this bundle denotes here a principle U(2) bundle, F $\to$ M such that F/U(1) = Fr. Such a lift is called a $\text{Spin}_\mathbb{C}$ structure. Two lifts, F and F´, are deemed equivalent if there is a bundle isomorphism from F to F´ that covers the projections ot Fr. The set of equivalence classes of lifts can be put in 1-1 correspondence with elements in $H^2(M; \mathbb{Z})$.

Let F $\to$ M denote now a $\text{Spin}_\mathbb{C}$ structure. Use $\mathbb{S}$ to denote the associated $\mathbb{C}^2$ bundle $F \times_{U(2)} \mathbb{C}^2$. Use $\det(\mathbb{S})$ in what follows to denote the complex hermitian line bundle $F \times_{U(1)} \mathbb{C}$. Having fixed a $\text{Spin}_\mathbb{C}$ structure, the associated Seiberg-Witten equations constitute a system of equations for a pair $(A, \psi)$ where A here denotes a connection on $\det(\mathbb{S})$ and $\psi$ denotes a section of $\mathbb{S}$.

To say more about these equations, introduce the Clifford multiplication homomorphism cl: $T^*M \to \text{End}(\mathbb{S})$. This homomorphism is such that $\text{cl}(b)^\dagger = -\text{cl}(b)$ and $\text{cl}(b)\text{cl}(b´) = -\text{cl}(*(b \wedge b´)) - \langle b, b´\rangle$. Here, $\langle\, ,\,\rangle$ denotes the metric inner product and $*$ denotes the associated Hodge star. The Seiberg-Witten equations involve two related homomorphisms. The first, $\hat{c}: \mathbb{S} \otimes T^*M \to \mathbb{S}$, is defined so as to send any given decomposable element $\eta \otimes b$ to $\text{cl}(b)\eta$. The second is a quadratic, bundle preserving map from $\mathbb{S}$ to $iT^*M$. The image of any given $\eta \in \mathbb{S}$ under the latter map is written in what follows as $\eta^\dagger \tau \eta$. It is defined by the rule $\langle b, \eta^\dagger \tau \eta\rangle = \eta^\dagger \text{cl}(b)\eta$.

Let A now denote a connection on $\det(\mathbb{S})$. In what follows, the Hodge star of its curvature 2-form is denoted by $B_A$, this a section of $iT^*M$. The connection A and the Levi-Civita connection on TM define a Hermitian connection on $\mathbb{S}$. The associated covariant derivative is denoted in what follows by $\nabla_A$. This covariant derivative is used to define the Dirac operator, $D^A = \hat{c}(\nabla_A): C^\infty(M; \mathbb{S}) \to C^\infty(M; \mathbb{S})$.

A pair $(A, \psi)$ of connection on $\det(\mathbb{S})$ and section $\psi$ of $\mathbb{S}$ obeys the simplest version of the Seiberg-Witten equations when

$$B_A - \psi^\dagger \tau \psi = 0 \quad \textit{and} \quad D^A \psi = 0 \,.$$

(3.1)

A rigorous definition of the Seiberg-Witten Floer homology involves solutions to perturbed versions of the equations in (3.1). The description of these perturbed equations requires a brief digression to set the stage.

To start the digression, remark that the equations in (3.1) are gauge invariant in the following sense: If u is a smooth map from M to U(1), then the pair $(A - 2u^{-1}du, u\psi)$ solves (3.1) if and only if $(A, \psi)$ does. A function, $\mathfrak{g}$, of pairs consisting of a connection on det($\mathbb{S}$) and a section of $\mathbb{S}$ is deemed gauge invariant when $\mathfrak{g}(A - 2u^{-1}du, u\psi) = \mathfrak{g}(A, \psi)$ for all $u \in C^\infty(M; U(1))$. The allowed sorts of functions from what Kronheimer and Mrowka call a large, separable Banach space of tame perturbations. Such a Banach space is described in Chapter 11 of [KM]. Somewhat more is said below about this. This Banach space of tame perturbations that is used here is denoted by $\mathcal{P}$. If $\mathfrak{g} \in \mathcal{P}$, then the differential of $\mathfrak{g}$ at any given $(A, \psi)$ defines section $(\mathfrak{T}|_{(A, \psi)}, \mathfrak{S}|_{(A, \psi)})$ of $iT^*M \oplus \mathbb{S}$ by writing $\frac{d}{dt} \mathfrak{g}(A+tb, \psi + t\eta)|_{t=0}$ as $\int_M (b \wedge *\mathfrak{T} - \frac{1}{2}(\eta^\dagger \mathfrak{S} + \mathfrak{S}^\dagger \eta))$. Each $\mathfrak{g} \in \mathcal{P}$ gives the equation

$$B_A - \psi^\dagger \tau \psi - \mathfrak{T}_{(A,\psi)} = 0 \quad and \quad D^A \psi - \mathfrak{S}|_{(A, \psi)} = 0.$$

(3.2)

Note that if u is a smooth map from M to $U(1)$, then $(A - 2u^{-1}du, u\psi)$ solves (3.2) if and only if $(A, \psi)$ does. Pairs of connection and section that are related in this way are said to be *gauge equivalent*.

There are corresponding Seiberg-Witten equations on $\mathbb{R} \times M$ that constitute a system of equations for a pair $\mathfrak{d} = (A, \psi)$ where A now denotes a map from $\mathbb{R}$ into the space of Hermitian connections on det($\mathbb{S}$) and $\psi$ denotes a map from $\mathbb{R}$ into the space of sections of $\mathbb{S} \to M$. With $s \in \mathbb{R}$ denoting the Euclidean coordinate, these equations read

- $\frac{\partial}{\partial s} A + B_A - \psi^\dagger \tau^k \psi - \mathfrak{T}(A, \psi) = 0$,
- $\frac{\partial}{\partial s} \psi + D^A \psi - \mathfrak{S}(A, \psi) = 0$.

(3.3)

Of particular interest here are *instanton* solutions. An instanton is a solution to (3.3) with $s \to +\infty$ limit and $s \to -\infty$ limit, each a solution to (3.1). If u is a smooth map from M to $U(1)$ and $(A, \psi)$ is a solution to (3.3), then so is $(A - 2u^{-1}du, u\psi)$.

**b) An overview of Seiberg-Witten Floer homology/cohomology**

This subsection very briefly summarizes the story from [KM]. To start, Kronheimer and Mrowka prove that (3.2) has but a finite set of solutions up to gauge equivalence if $\mathfrak{g}$ is chosen from a certain residual subset in the Banach space $\mathcal{P}$. With one caveat, these equivalence classes form a basis for the chain complex that defines the Seiberg-Witten Floer homology. The caveat concerns the case when the first Chern class of det($\mathbb{S}$) is a torsion class. The situation here is more complicated by virtue of the fact that (3.2) admits solutions with $\psi$ identically zero when $c_1(det(\mathbb{S}))$ is torsion. These $\psi = 0$ solutions are deemed to be *reducible*, and those with $\psi$ somewhere non-zero are deemed to be *irreducible*. Here is the salient distinction: The group $C^\infty(M; U(1))$ acts with trivial stabilizer on any pair $(A, \psi)$ with $\psi$ somewhere non-zero, but it acts with stabilizer $U(1)$ on any (A, 0). Here, $U(1) \subset C^\infty(M; U(1))$ is identified with the constant maps. This distinction makes for a chain complex when $c_1(\det(\mathbb{S}))$ is torsion that has one generator for each gauge equivalence class of irreducible solution to (3.2), and a countable set of generators for each gauge equivalence class of reducible solution to (3.2). The chain complex for the Seiberg-Witten Floer homology is denoted in what follows by $\mathcal{C}_{SW}$. This $\mathbb{Z}$ module is finite when $c_1(\det(\mathbb{S}))$ is not torsion, but not finitely generated otherwise.

The complex $\mathcal{C}_{SW}$ has a natural, relative $\mathbb{Z}/p\mathbb{Z}$ grading, where p here denotes the divisibility of the class $c_1(\det(\mathbb{S}))$ in $H^2(M; \mathbb{Z})$/torsion. The complex is $\mathbb{Z}$ graded when $c_1(\det(\mathbb{S}))$ is torsion. This grading is described in some detail momentarily. Suffice it to say for now that the relative grading between two irreducible generators is defined to be minus the spectral flow for a certain 1-parameter family of unbounded, self-adjoint, operators (with compact resolvent) on $L^2(M; iT^*M \oplus \mathbb{S} \oplus i\mathbb{R})$. This family is constructed from a path, parametrized by [0, 1], of pairs $(A, \psi)$ with A a connection on $\det(\mathbb{S})$ and $\psi$ a section of $\mathbb{S}$. This path starts at the first irreducible solution, and ends at the second. Meanwhile, the operator that is parametrized by any such pair $(A, \psi)$ is a self map of $C^\infty(M; iT^*M \oplus \mathbb{S} \oplus i\mathbb{R})$; it is defined from the linearization of (3.2) at $(A, \psi)$.

In the case where $c_1(\det(\mathbb{S}))$ is torsion, the countable set of cycles that correspond to any given irreducible generator can be labled by a set of the form {k, k - 2, …} where $k \in \mathbb{Z}$. The relative grading between cycles k - 2j and k - 2j´ is 2(j - j´). The integer k can be fixed once a fiducial, irreducible configuration is chosen to define the zero point for the grading. Given such a choice, k is then minus the spectral flow for a family of self-adjoint differntial operators that starts at a certain operator that is parametrized by the fiducial configuration and ends at one parametrized by a suitable irreducible configuration lying very near the given reducible solution. In the case when $c_1(\det(\mathbb{S}))$ is torsion, the $\mathbb{Z}$-module $\mathcal{C}_{SW}$ is finitely generated in each degree.

The differential that defines the Seiberg-Witten Floer homology is defined using an algebraic count of instanton solutions to (3.3) as defined by any $\mathfrak{g}$ from a certain residual subset in $\mathcal{P}$. The differential decrease that $\mathbb{Z}/p\mathbb{Z}$ degree by 1. To say more about the differential, note that Kronheimer and Mrowka prove the following: Let $\mathfrak{c}_-$ and $\mathfrak{c}_+$ denote pairs of irreducible solutions to (3.2) and introduce $\mathcal{M}(\mathfrak{c}_-, \mathfrak{c}_+)$ to denote the set of instanton solutions to (3.3) with $s \to -\infty$ limit equal to $\mathfrak{c}_-$ and with $s \to \infty$ limit equal to $u \cdot \mathfrak{c}_+$ with $u \in C^\infty(M; S^1)$. This set depends only on the gauge equivalence classes of $\mathfrak{c}_-$ and $\mathfrak{c}_+$ in the following sense: Suppose that $u \in C^\infty(M; U(1))$. Let $\mathfrak{d} = (A, \psi) \in \mathcal{M}(\mathfrak{c}_-, \mathfrak{c}_+)$, then $u \cdot \mathfrak{d} \in \mathcal{M}(u \cdot \mathfrak{c}_+, \mathfrak{c}_-)$ where $u \cdot (A, \psi)$ is shorthand for $(A - 2u^{-1}du, u\psi)$.

Given that $\mathfrak{g}$ comes from a certain residual subset of $\mathcal{P}$, this $\mathcal{M}(\mathfrak{c}_-, \mathfrak{c}_+)$ has the structure of a smooth, finite dimensional manifold. There is one zero dimensional component if $\mathfrak{c}_-$ is gauge equivalent to $\mathfrak{c}_+$; and in this case, $\mathcal{M}(\mathfrak{c}_-, \mathfrak{c}_+)$ consists of the constant map $s \to \mathfrak{c}_-$. There are no zero dimensional components otherwise. Meanwhile, there is a finite set of 1-dimensional components of $\mathcal{M}(\mathfrak{c}_-, \mathfrak{c}_+)$; and each component is an orbit of the $\mathbb{R}$ action that is induced by translation along the $\mathbb{R}$ factor of $\mathbb{R} \times M$. Such 1-dimensional components exist only in the case where the degree of $\mathfrak{c}_+$ is one less than that of $\mathfrak{c}_-$. Use $\mathcal{M}_1(\mathfrak{c}_-, \mathfrak{c}_+)$ in what follows to denote the space of 1-dimensional components of $\mathcal{M}(\mathfrak{c}_-, \mathfrak{c}_+)$

Each component of $\mathcal{M}_1(\mathfrak{c}_-, \mathfrak{c}_+)$ has a corresponding sign. This sign is obtained by comparing the orientation given by the generator of the $\mathbb{R}$ action with an orientation that is defined using Quillen's ideas about determinant line bundles. Some what more is said about this below, but in any event, the full story is given in [KM]. Let $\sigma(\mathfrak{c}_-, \mathfrak{c}_+)$ denote the sum of these signs when $\mathcal{M}_1(\mathfrak{c}_-, \mathfrak{c}_+)$ is non-empty, or zero when it is. In the case when $c_1(\det(\mathbb{S}))$ is not torsion, the differential that defines the Seiberg-Witten Floer homology acts on any given generator $\mathfrak{c}$ as $\delta \mathfrak{c} = \sum_{\mathfrak{c}' \in \mathcal{C}_{SW}} \sigma(\mathfrak{c}, \mathfrak{c}') \mathfrak{c}'$. In the case when $c_1(\det(\mathbb{S}))$ is torsion, what is written here defines the part of the differential that involves the irreducible generators. Only this part is needed for the proofs of the theorems in the introduction. This being the case, the reader can consult [KM] to see how the rest of the differential is defined.

The homology of this differential on $\mathcal{C}_{SW}$ is denoted by $\widehat{\mathcal{H}}_*$ in [KM], and so denoted by $\widehat{\mathcal{H}}_*$ here. This homology is finitely generated in the case when $c_1(\det(\mathbb{S}))$ is not torsion. In the case when this class is torsion, $\widehat{\mathcal{H}}_*$ is finitely generated in each degree; and the set of degrees where it is non-zero is bounded from above but unbounded from below.

The Seiberg-Witten Floer cohomology is defined by the dual differential on the $\mathbb{Z}$-module $\mathcal{C}^{SW} = \text{Hom}(\mathcal{C}_{SW}, \mathbb{Z})$. This differential, $\delta^*$, acts on any given cocycle C by $(\delta^* C)(\cdot) = C(\delta(\cdot))$. Note that the basis described above for $\mathcal{C}_{SW}$ supplies a canonical dual basis for $\mathcal{C}^{SW}$ and a $\mathbb{Z}/p\mathbb{Z}$ grading. This differential sends any given basis element $\mathfrak{c}$ of $\mathcal{C}^{SW}$ to

$$\delta^* \mathfrak{c} = \sum_{\mathfrak{c}' \in \mathcal{C}^{SW}} \sigma(\mathfrak{c}', \mathfrak{c}) \mathfrak{c}'.$$

(3.4)

Note that it increases the $\mathbb{Z}/p\mathbb{Z}$ degree by 1. The resulting cohomology groups are denoted in what follows by $\widehat{\mathcal{H}}^*$.

Keep in mind that the definition of these groups requires the choice of a function from a certain residual subset in $\mathcal{P}$. However, two such functions give isomorphic versions of Seiberg-Witten homology and cohomology. More is said below about the criteria for admission in this residual set.

**c) Contact forms and Seiberg-Witten equations**

Suppose now that a is a given contact 1-form on M. Fix a metric on M for which $*da = 2a$ and $|a| = 1$. Note that such a metric on TM is neither more nor less than an almost complex structure, J, on kernel(a) such $da(\cdot, J(\cdot))$ is a metric on the kernel of a. In particular, a pair (a, J) of contact form and almost complex structure in $\mathcal{J}_a$ supplies M with a canonical metric.

With the metric fixed, let $F \to M$ denote a $Spin_{\mathbb{C}}$ structure. The endomorphism cl(a) on $\mathbb{S}$ has square -1 and so its eigenspaces in each fiber define a splitting of $\mathbb{S}$ as the orthogonal, direct sum of two complex, Hermitian line bundles. This direct sum is written in what follows as $E \oplus EK^{-1}$ where $E \to \mathbb{S}$ and $K \to S$ are complex line bundles. The convention has cl(a) act as i on the first summand and -i on the second. The bundle $K^{-1} \to S$ is isomorphic as an SO(2) bundle to the kernel of a in TM with the orientation defined by da. Note that any given equivalence class of complex line bundle can arise in this manner from some $Spin_{\mathbb{C}}$ structure on M. Moreover, two $Spin_{\mathbb{C}}$ structures have isomorphic versions of E if and only if they are equivalent.

The contact form a determines a *canonical $Spin_{\mathbb{C}}$ structure*, this the $Spin_{\mathbb{C}}$ structure for which the spinor bundle decomposes as $\mathbb{S} = \mathbb{S}_I = I_{\mathbb{C}} \oplus K^{-1}$ where $I_{\mathbb{C}} \to M$ denotes the trivial complex line bundle. Fix a unit norm section $1_{\mathbb{C}}$ of $I_{\mathbb{C}}$. Such a section defines a *canonical connection* on $K^{-1} = \det(\mathbb{S}_I)$. This is the unique connection for which the section $\psi_I = (1_{\mathbb{C}}, 0)$ of $\mathbb{S}_I$ is annihilated by the corresponding Dirac operator. This canonical connection is written as $A_K$.

Let $\mathbb{S} = E \oplus EK^{-1}$ now denote the spinor bundle for some other $Spin_{\mathbb{C}}$ structure. Any given connection on $\det(\mathbb{S}) = E^2K^{-1}$ can be written as $A_K + 2A$ where A is a connection on E. With A a connection on E, the symbol $D_A$ is used to denote the operator $D^{A_K + 2A}$ that appears in (3.1)-(3.3). Conn(E) is used in what follows to denote the Frechêt space of smooth, Hermitian connections on E.

With the splitting $\mathbb{S} = E \oplus EK^{-1}$ given, the corresponding components of any given section $\psi$ of $\mathbb{S}$ are written as $(\alpha, \beta)$. Thus, $\alpha$ is a section of E and $\beta$ one of $EK^{-1}$.

The contact form is also used to define a certain family of perturbations for use in (3.2) and (3.3). This family is parametrized by $[1, \infty)$. To set the stage, view $\mathcal{P}$ now as a Banach space of functions on $Conn(E) \times C^{\infty}(M; \mathbb{S})$. Any given $r \in [1, \infty)$ version of these equations requires the choice of a function $\mathfrak{g}$ from $\mathcal{P}$. These equations, viewed now as equations for a pair $(A, \psi)$ from $Conn(E) \times C^{\infty}(M; \mathbb{S})$, read

- $B_A - r(\psi^{\dagger}\tau\psi - ia) - \mathfrak{T}|_{(A,\psi)} + \frac{1}{2} B_{A_K} = 0.$
- $D_A\psi - \mathfrak{S}|_{(A,\psi)} = 0.$

(3.5)

Here, $\mathfrak{T}$ and $\mathfrak{S}$ are defined from $\mathfrak{g}$ as before. Meanwhile, $B_{A_K}$ is the curvature 2-form for the connection $A_K$. The associated version of (3.3) for a map $s \to (A, \psi)$ from $\mathbb{R}$ to $\mathrm{Conn}(E) \times C^{\infty}(M; \mathbb{S})$ is

- $\frac{\partial}{\partial s} A + B_A - r(\psi^{\dagger}\tau\psi - ia) - \mathfrak{T}|_{(A,\psi)} + \frac{1}{2} B_{A_K} = 0$.
- $\frac{\partial}{\partial s} \psi + D_A\psi - \mathfrak{S}|_{(A,\psi)} = 0$.

(3.6)

Equations (3.5) and (3.6) can be made to look like (3.2) and (3.3) by replacing $\psi$ in the latter by $r\psi$. It is left to the reader to derive the relation between the respective versions of what is denoted by $\mathfrak{g}$.

**d) The Banach space $\mathcal{P}$**

As noted in Section 3b, the Seiberg-Witten Floer homology and cohomology can be defined only after choosing a function $\mathfrak{g}$ from a certain residual subset of $\mathcal{P}$. There are two criteria for membership in this set when the first Chern class of E is not torsion and three criteria when it is. The first concerns the linearized version of (3.5). To say more, fix any pair $\mathfrak{c} = (A, \psi) \in \mathrm{Conn}(E) \times C^{\infty}(M; \mathbb{S})$. Define an operator $\mathfrak{L}_{\mathfrak{c}}$ with domain and range $C^{\infty}(M; iT^*M \oplus \mathbb{S} \oplus i\mathbb{R})$ as follows: It sends any given triple $(b, \eta, \phi)$ in its domain to the section of $iT^*M \oplus \mathbb{S} \oplus i\mathbb{R}$ whose three components are:

- $*db - d\phi - 2^{-1/2} r^{1/2} (\psi^{\dagger}\tau\eta + \eta^{\dagger}\tau\psi) - \mathfrak{t}_{(A,\psi)}(b, \eta)$,
- $D_A\eta + 2^{1/2} r^{1/2}(\mathrm{cl}(b)\psi + \phi\psi) - \mathfrak{s}_{(A,\psi)}(b, \eta)$,
- $*d*b - 2^{-1/2} r^{1/2}(\eta^{\dagger}\psi - \psi^{\dagger}\eta)$ ,

(3.7)

where the pair $(\mathfrak{t}_{(A,\psi)}, \mathfrak{s}_{(A,\psi)})$ denotes the operator on $C^{\infty}(M; iT^*M \oplus \mathbb{S})$ that sends a given section $(b, \eta)$ to $(\frac{d}{dt} \mathfrak{T}(A+tb, \psi+t\eta), \frac{d}{dt} \mathfrak{S}(A+tb, \psi+t\eta))|_{t=0}$. The operator $\mathfrak{L}_{\mathfrak{c}}$ is symmetric and extends to $L^2(M; iT^*M \oplus \mathbb{S} \oplus i\mathbb{R})$ as an unbounded, self-adjoint operator with dense domain $L^2_1(M; iT^* \oplus \mathbb{S} \oplus i\mathbb{R})$. As such, it has pure point spectrum and each eigenvalue has finite multiplicity. Moreover, the spectrum is unbounded from above and from below. Finally, every eigenvector is smooth.

A function $\mathfrak{g}$ from $\mathcal{P}$ can be used to define the Seiberg-Witten Floer homology only when the following criteria is met:

> Criteria 1: *If $\mathfrak{c}$ is an irreducible solution to (3.5), then the operator $\mathfrak{L}_{\mathfrak{c}}$ has trivial kernel. If $\mathfrak{c}$ is a reducible solution to (3.5), then* $\mathrm{kernel}(\mathfrak{L}_{\mathfrak{c}})$ *consists of the constant sections of the $i\mathbb{R}$ summand of $iT^*M \oplus \mathbb{S} \oplus i\mathbb{R}$.*

(3.8)

The second required property for $\mathfrak{g}$ involves the operator on $\mathbb{R} \times M$ that arises from the linearized version of (3.6). To elaborate, suppose that $s \to \mathfrak{d}(s) = (A, \psi)$ is a smooth map from $\mathbb{R}$ into $\mathrm{Conn}(E) \times C^{\infty}(M; \mathbb{S})$ that has $s \to \pm\infty$ limits. Let $\mathfrak{c}_{\pm}$ denote the latter. Now define the operator $\mathfrak{D}_{\mathfrak{d}}$ from $C^{\infty}(\mathbb{R} \times M; iT^*M \oplus \mathbb{S} \oplus i\mathbb{R})$ to itself as follows: It sends a triple $(b, \eta, \phi)$ to the section whose respective three components are

- $\frac{\partial}{\partial s} b + *db - d\phi - 2^{-1/2} r^{1/2} (\psi^{\dagger}\tau\eta + \eta^{\dagger}\tau\psi) - \mathfrak{t}_{(A,\psi)}(b, \eta)$,
- $\frac{\partial}{\partial s} \eta + D_A \eta + 2^{1/2} r^{1/2}(\mathrm{cl}(b)\psi + \phi\psi) - \mathfrak{s}_{(A,\psi)}(b, \eta)$,
- $\frac{\partial}{\partial s} \phi + *d*b - 2^{-1/2} r^{1/2}(\eta^{\dagger}\psi - \psi^{\dagger}\eta)$.

(3.9)

This operator extends to define a bounded operator from the $L^2_1(\mathbb{R} \times M; iT^*M \oplus \mathbb{S} \oplus i\mathbb{R})$ to $L^2(\mathbb{R} \times M; iT^*M \oplus \mathbb{S} \oplus i\mathbb{R})$. If both $\mathfrak{c}_+$ and $\mathfrak{c}_-$ are irreducible and if both the $\mathfrak{c} = \mathfrak{c}_+$ and $\mathfrak{c} = \mathfrak{c}_-$ versions of $\mathfrak{L}_{\mathfrak{c}}$ have trivial kernel, then this extended version of $\mathfrak{D}_{\mathfrak{d}}$ is a Fredholm operator. This understood, what follows is the second requirement on $\mathfrak{g}$ for its use to define the Seiberg-Witten Floer homology and cohomology.

> Criteria 2: *Let $s \to \mathfrak{d}(s)$ denote an instanton solution to (3.6) such that both $|s| \to \infty$ limits are irreducible and such that their corresponding versions of $\mathfrak{L}_{\mathfrak{c}}$ have trivial kernel. Then $\mathfrak{D}_{\mathfrak{d}}$ has trivial kernel.*

(3.10)

The third requirement on $\mathfrak{g}$ for its use in defining the Seiberg-Witten Floer homology and cohomology concerns the operator $\mathfrak{D}_{\mathfrak{d}}$ when $s \to \mathfrak{d}(s)$ is an instanton solution to (3.6) with at least one $|s| \to \infty$ limit reducible. As these solutions play no essential role in what follows, this third criteria will not be stated explicitly. The reader can refer instead to [KM].

The space $\mathcal{P}$ has the following property: Let $\mathfrak{g} \in \mathcal{P}$. Just as the first derivatives of $\mathfrak{g}$ at any given $(A, \psi) \in \mathrm{Conn}(E) \times C^{\infty}(M; \mathbb{S})$ define the smooth section $(\mathfrak{T}|_{(A,\psi)}, \mathfrak{S}|_{(A,\psi)})$ of the bundle $iT^*M \oplus \mathbb{S}$, so the derivatives of $\mathfrak{g}$ to order $k \geq 1$ at $(A, \psi)$ define a smooth section of $\otimes_k (iT^*M \oplus \mathbb{S})$. Let $\mathfrak{g}_k|_{(A,\psi)}$ denote the latter. The derivatives of this section to any given order are bounded by an appropriate $(A, \psi)$-dependent multiple of $\|\mathfrak{g}\|_{\mathcal{P}}$. Here, $\|\cdot\|_{\mathcal{P}}$ denotes the Banach space norm on $\mathcal{P}$. For example, bounds of this sort appear in Theorem 11.6.1 of [KM]. See also Proposition 2.5 in [T3].

What follows describes some of the simplest non-constant functions in $\mathcal{P}$. To start, let $\mu$ denote a smooth, coclosed 1-form on M. Thus, $d*\mu = 0$. Use $\mathfrak{e}_{\mu}$ to denote the function on $\mathrm{Conn}(E)$ whose value on any given connection A is

$$\mathfrak{e}_{\mu}(A) = i\int_M \mu \wedge *B_A .$$

(3.11)

View $\mathfrak{e}_\mu$ as a function on $\mathrm{Conn}(E) \times C^\infty(M; \mathbb{S})$ that is independent of the second factor. Viewed in this light, $\mathfrak{e}_\mu$ is a candidate for a function from $\mathcal{P}$; and this is the case if $\mu$ comes from a certain dense set of coclosed 1-forms. To say more, let $\Omega_0$ denote the vector space of finite linear combinations of coclosed eigenfunctions of the operator $*d$, here sending $C^\infty(M; T^*M)$ to itself. Then $\mathcal{P}$ contains the linear space $\{\mathfrak{e}_\mu: \mu \in \Omega_0\}$. Moreover, the following is true: For each $k \geq 0$, there is a constant $c_k$ such that the $C^k$ norm of any given $\mu \in \Omega_0$ is bounded by $c_k \|\mathfrak{e}_\mu\|_\mathcal{P}$.

Note that these norm bounds imply that the function $\mu \to \|\mathfrak{e}_\mu\|_\mathcal{P}$ defines a norm on $\Omega_0$; and they imply that the completion of $\Omega_0$ with respect to this norm is a subspace of $C^\infty(M; T^*M)$. Use $\Omega$ in what follows to denote this completion. This norm on $\Omega$ is called the '$\mathcal{P}$ norm' in what follows.

The perturbations of particular interest in what follows have the from $\mathfrak{e}_\mu + \mathfrak{p}$ where $\mu \in \Omega$ and $\mathfrak{p} \in \mathcal{P}$ with $\|\mathfrak{p}\|_\mathcal{P}$ very much smaller than $\|\mathfrak{e}_\mu\|_\mathcal{P}$ and with $\|\mathfrak{e}_\mu\|_\mathcal{P} \ll 1$. Note that with $\mathfrak{g} = \mathfrak{e}_\mu$, the pair $(\mathfrak{T}, \mathfrak{S})$ in (3.6) is $(\mathfrak{T} = i*d\mu, \mathfrak{S} = 0)$. In this case, the terms $\mathfrak{t}$ and $\mathfrak{s}$ are absent in (3.7).

**e) Degrees and signs**

This subsection has two parts. The first elaborates on the $\mathbb{Z}/p\mathbb{Z}$ degree assignments to the generators of $\mathcal{C}_{SW}$, and the second elaborates on the signs that are used to define the differential on $\mathcal{C}_{SW}$.

*Part 1*: As mentioned above, the relative $\mathbb{Z}/p\mathbb{Z}$ degree between two generators of the Seiberg-Witten Floer complex is defined using the spectral flow for a 1-parameter family of unbounded, self-adjoint operators on $L^2(M; iT^*M \oplus \mathbb{S} \oplus i\mathbb{R})$. To elaborate on the relevant case, suppose that $\mathfrak{c}_-$ and $\mathfrak{c}_+$ are irreducible solutions to some $r$ and $\mathfrak{g}$ version of (3.5) and are such that the respective $\mathfrak{c} = \mathfrak{c}_-$ and $\mathfrak{c} = \mathfrak{c}_+$ versions of $\mathfrak{L}_\mathfrak{c}$ have trivial cokernel. Fix a path $s \to \mathfrak{d}(s) \in \mathrm{Conn}(E) \times C^\infty(M; \mathbb{S})$ parametrized by $[0, 1]$ such that $\mathfrak{d}(0) = \mathfrak{c}_-$ and $\mathfrak{d}(1) = \mathfrak{c}_+$. If the chosen path $s \to \mathfrak{d}(s)$ is sufficiently generic, then there will be at most one eigenvalue very near 0 at any $s \in [0, 1]$. Such an eigenvalue will have multiplicity 1 and vary smoothly as the parameter $s$ is changed as long as the eigenvalue is sufficiently close to 0. Moreover, if it changes sign as $s$ varies, it zero crosses zero with non-zero derivative. This understood, the spectral flow for the path is equal to the number of points in $(0, 1)$ where an eigenvalue crosses zero with positive derivative, minus the number where it crosses zero with negative derivative. (See, for example [T9].) This spectral flow is denoted in what follows by $f(\mathfrak{c}_-, \mathfrak{c}_+)$. Although a particular path must be chosen to compute this number, the number itself does not depend on the path. However, only the $\mathbb{Z}/p\mathbb{Z}$ reduction of $f(\mathfrak{c}_-, \mathfrak{c}_+)$ is gauge invariant. Granted the preceding, now view $\mathfrak{c}_-$ and $\mathfrak{c}_+$ as generators of $\mathcal{C}_{SW}$. Then

$$\text{degree}(\mathfrak{c}_+) - \text{degree}(\mathfrak{c}_-) = -f(\mathfrak{c}_-, \mathfrak{c}_+) \; mod(p) \,.$$

(3.12)

When $c_1(\det(\mathbb{S}))$ is torsion, there will be reducible solutions to (3.5). As noted above, the countable set of cycles that correspond to any given reducible generator can be labled by a set of the form $\{k, k-2, \ldots\}$ where $k \in \mathbb{Z}$. The relative grading between cycles $k - 2j$ and $k - 2j'$ is $2(j - j')$. Let $\mathfrak{c} \in \text{Conn}(E) \times C^\infty(M; \mathbb{S})$ denote a pair where $\mathfrak{L}_\mathfrak{c}$ has trivial kernel. Then the relative degree difference, $k - \text{degree}(\mathfrak{c})$, is defined to be minus the spectral flow for the family $\{\mathfrak{L}_{\mathfrak{d}(s)}\}_{s \in [0,1]}$ where $\mathfrak{d}: [0, 1] \in \text{Conn}(E) \times C^\infty(M; \mathbb{S})$ is a path that starts at $\mathfrak{c}$ and ends at an irreducible configuration that is very near the reducible one.

*Part 2*: The signs that appear in (3.4) are also defined using families of operators; in this case the operators that appear in (3.9). A digression is needed first to say more about how this is done (See Chapter 20 of [KM]). To start the digression, fix $r$ and $\mathfrak{g}$. Let $\mathfrak{c}_-$ and $\mathfrak{c}_+$ denote elements in $\text{Conn}(E) \times C^\infty(M; \mathbb{S})$ where $\mathfrak{L}_{(\cdot)}$ has trivial kernel. Introduce $\mathfrak{P} = \mathfrak{P}(\mathfrak{c}_-, \mathfrak{c}_+)$ to denote the space of piecewise differentiable maps from $\mathbb{R}$ to $\text{Conn}(E) \times C^\infty(M; \mathbb{S})$ that have $s \to -\infty$ limit which is gauge equivalent to $\mathfrak{c}_-$ and $s \to \infty$ limit which is gauge equivalent to $\mathfrak{c}_+$. Each $\mathfrak{d} \in \mathfrak{P}$ has its corresponding version of $\mathfrak{D}_\mathfrak{d}$ as given in (3.9); here viewed as a Fredholm operator from $L^2_1(\mathbb{R} \times M; iT^*M \oplus \mathbb{S} \oplus i\mathbb{R})$ to $L^2(\mathbb{R} \times M; iT^*M \oplus \mathbb{S} \oplus i\mathbb{R})$. Quillen [Q] explained how such a family of operators can be used to construct a real line bundle, $\det(\mathfrak{D}) \to \mathfrak{P}$. The fiber of $\det(\mathfrak{D})$ at a given $\mathfrak{d} \in \mathfrak{P}$ is canonical identified with $\wedge^{\max}(\text{kernel}(\mathfrak{D}_\mathfrak{d})) \times_\mathbb{R} (\wedge^{\max}\text{cokernel}(\mathfrak{D}_\mathfrak{d}))^*$ if either the kernel or cokernel of $\mathfrak{D}_\mathfrak{d}$ is non-trivial. As explained in [KM], this real line bundle is suitably gauge invariant and has a gauge invariant orientations. Use $\Lambda(\mathfrak{c}_-, \mathfrak{c}_+)$ to denote the 2-element set of orientations for $\det(\mathfrak{D}) \to \mathfrak{P}$; viewed here as a non-trivial $\mathbb{Z}/2\mathbb{Z}$ module.

Collectively, the modules $\Lambda(\cdot, \cdot)$ have the following three important features: To state the first, let $\mathfrak{c}_-$, $\mathfrak{c}_0$ and $\mathfrak{c}_+$ denote elements in $\text{Conn}(E) \times C^\infty(M; \mathbb{S})$ where $\mathfrak{L}_{(\cdot)}$ has trivial kernel. There is in this case the composition law $\Lambda(\mathfrak{c}_-, \mathfrak{c}_0) \otimes_{\mathbb{Z}/2\mathbb{Z}} \Lambda(\mathfrak{c}_0, \mathfrak{c}_+) = \Lambda(\mathfrak{c}_-, \mathfrak{c}_+)$. Second, $\Lambda(\mathfrak{c}_-, \mathfrak{c}_+)^* = \Lambda(\mathfrak{c}_+, \mathfrak{c}_-)$. Note that these last two properties imply that $\Lambda(\mathfrak{c}_-, \mathfrak{c}_+)$ can be written as $\Lambda(\mathfrak{c}_-) \otimes_{\mathbb{Z}/2\mathbb{Z}} \Lambda(\mathfrak{c}_+)^*$ where $\Lambda(\cdot)$ is a $\mathbb{Z}/2\mathbb{Z}$ module that is assigned to each gauge equivalence class of pair $\mathfrak{c} \in \text{Conn}(E) \times C^\infty(M; \mathbb{S})$ where $\mathfrak{L}_\mathfrak{c}$ has trivial kernel.

To state the final salient property, assume now that both $\mathfrak{c}_-$ and $\mathfrak{c}_+$ are irreducible solutions to (3.5) and that both the $\mathfrak{c} = \mathfrak{c}_-$ and $\mathfrak{c} = \mathfrak{c}_+$ versions of $\mathfrak{L}_\mathfrak{c}$ have trivial kernel. Assume that each $\mathfrak{d} \in \mathcal{M}(\mathfrak{c}_-, \mathfrak{c}_+)$ version of $\mathfrak{D}_\mathfrak{d}$ has trivial cokernel. Then the restriction of $\Lambda(\mathfrak{c}_-, \mathfrak{c}_+)$ to $\mathcal{M}(\mathfrak{c}_-, \mathfrak{c}_+)$ is canonically isomorphic to the latter's orientation sheaf.

With the digression now over, assume now that $r$ and $\mathfrak{g}$ obey the two criteria given in (3.8) and (3.10). What follows explains how the signs for the differential on $\mathcal{C}_{SW}$ are determined. Assign to each gauge equivalence class of irreducible solution to (3.5) an element, o(·), in the corresponding version of Λ(·). Suppose next that $c_-$ and $c_+$ are irreducible solutions to (3.5). Then $o(c_-)o(c_+) \in \Lambda(c_-, c_+)$ and so defines an orientation for each component of $\mathcal{M}(c_-, c_+)$. Meanwhile, each 1-dimensional component of this space is oriented by the generator of the $\mathbb{R}$ action that is induced by translation along the $\mathbb{R}$ factor of $\mathbb{R} \times M$. This understood, a given 1-dimensional component contributes +1 to $\sigma(c_-, c_+)$ when these two orientations agree; and it contributes -1 to $\sigma(c_-, c_+)$ when these two orientations disagree.

### f) Other functions on Conn(E) × $C^\infty$(M; $\mathbb{S}$)

Various functions on Conn(E) × $C^\infty$(M; $\mathbb{S}$) play a central role in subsequent parts of the story. The first of these is the gauge invariant function, E, from Conn(E) with value on $A \in$ Conn(E) given by

$$\mathrm{E}(A) = i \int_M a \wedge *B_A$$

(3.13)

The second function, the *Chern-Simons function*, also comes from Conn(E). Its definition requires first a choice, $A_E$, of a fiducial connection on E. It proves useful to choose the latter to be a connection whose curvature 2-form is harmonic. With $A_E$ chosen once and for all, the value of the Chern-Simons function on $A \in$ Conn(E) is given as follows: Write $A = A_E + \hat{a}_A$. Then

$$\mathrm{cs}(A) = -\int_M \hat{a}_A \wedge *d\hat{a}_A - 2 \int_M \hat{a}_A \wedge *(B_E + \tfrac{1}{2} B_{A_K}),$$

(3.14)

where $*B_E$ is the curvature of $A_E$ and $*B_{A_K}$ is that of $A_K$. Note that cs is fully gauge invariant only in the case where $c_1(\det(\mathbb{S}))$ is a torsion class.

The third of the four functions is denoted by $\mathfrak{a}$. It's critical points are the solutions to (3.5) and the maps that solve (3.6) parametrize the integral curves of its gradient vector field. This function is given by

$$\mathfrak{a} = \tfrac{1}{2}(\mathrm{cs} - r\mathrm{E}) + \mathfrak{g} + r\int_M \psi^\dagger D_A \psi.$$

(3.15)

The definition of the fourth of these functions requires first the choice of a section $\psi_E \in C^\infty(M; \mathbb{S})$ such that the $r = 1$, $\mathfrak{g} = 0$ and $c = (A_E, \psi_E)$ version of $\mathfrak{L}_c$ has trivial kernel. With this done, the fourth function, $f$, is a locally constant function defined off of a codimension 1 subvariety in Conn(E) × $C^\infty(M; \mathbb{S})$. Its value on any given $(A, \psi)$ is the spectral flow for the path of operators $\{\mathcal{L}^s\}_{s \in [0,1]}$ with $\mathcal{L}^s$ denoting the version of (3.7) that has $r^s$ in lieu of r, $s\mathfrak{g}$ in lieu of $\mathfrak{g}$, and $(A_E + s\hat{a}_A, \psi_E + s(\psi - \psi_E))$ in lieu of $(A, \psi)$.

The respective values of $\mathfrak{cs}$, $\mathfrak{a}\,f$ on any given pair $(A, \psi)$ are identical to those on $(A - u^{-1}du, u\psi)$ if $u \in C^\infty(M; S^1)$ is homotopically trivial or if $c_1(\det(\mathbb{S}))$ is a torsion class. However, the functions

$$\mathfrak{cs}^f = \mathfrak{cs} - 4\pi^2 f \quad and \quad \mathfrak{a}^f = \mathfrak{a} - 2\pi^2 f$$

(3.16)

are fully gauge invariant. This is to say that their values on any given $(A, \psi)$ are identical to those on $(A - u^{-1}du, u\psi)$ for all $u \in C^\infty(M; S^1)$. Note, however, that $\mathfrak{cs}^f$ and $\mathfrak{a}^f$ are only defined on the complement of the codimension 1 subvariety that consists of the elements $\mathfrak{c} \in \operatorname{Conn}(E) \times C^\infty(M; S^1)$ where $\mathcal{L}_\mathfrak{c}$ has a non-trivial kernel.

### g) Special cases

Of principle interest in what follows are the cases of (3.5) and (3.6) where the metric is such that $|a| = 1$ and $*da = 2a$. Thus, the metric comes via an almost complex structure J on the kernel of a such that $da(\cdot, J(\cdot))$ is a metric on kernel(a). Assume as well that the function $\mathfrak{g}$ has the form $\mathfrak{e}_\mu$ with $\mathfrak{e}_\mu$ as in (3.11) as defined using a given coclosed 1-form $\mu \in \Omega$. In this case, the equations in (3.5) read

- $B_A - r(\psi^\dagger \tau \psi - ia) - i*d\mu + \frac{1}{2} B_{A_K} = 0$.
- $D_A \psi = 0$.

(3.17)

With $\mu \in \Omega$ fixed, and $r \geq 1$ chosen, $\mathcal{M}^r$ henceforth denotes the space of gauge equivalence classes of solutions to (3.17).

In this case, the equations in (3.6) for instantons are

- $\frac{\partial}{\partial s} A + B_A - r(\psi^\dagger \tau \psi - ia) - i*d\mu + \frac{1}{2} B_{A_K} = 0$.
- $\frac{\partial}{\partial s} \psi + D_A \psi = 0$.

(3.18)

As noted above, the equations in (3.18) assert that $(A, \psi)$ is a critical point of the function

$$\mathfrak{a} = \tfrac{1}{2}(\mathfrak{cs} - rE) + \mathfrak{e}_\mu + r \int_M \psi^\dagger D_A \psi \, ;$$

(3.19)

and the equations in (3.19) assert that the path $s \to (A, \psi)|_s$ is an integral curve of the minus the gradient of $\mathfrak{a}$.

## 4. Proof of Theorem 1

The purpose of this section is to explain how Theorem 1 follows from a collection of theorems about the large r versions of (3.17) and (3.18). To set the stage, note that the definition of embedded contact homology requires the choice of a suitable pair (a, J) of contact 1-form from Lemma 2.1's residual set and almost complex structure from $\mathcal{J}_a$. Two different choices can, in principle, define different complexes and/or different differentials. (As noted above, Mike Hutchings conjectured that the resulting homology groups are isomorphic.) Likewise, the definition of the Seiberg-Witten Floer cohomology or homology requires choices. Such choices in this case consists a 4-tuple (â, $\hat{J}$, r, $\mathfrak{g}$) where â is a suitably chosen contact form and $\hat{J}$ is an almost complex structure from $\mathcal{J}_{\hat{a}}$. These are used to define the metric on M; and then â is used as the contact 1-form in the corresponding versions of (3.17) and (3.18). Meanwhile r ≥ 1 is a real number and $\mathfrak{g} \in \mathcal{P}$ is a suitable perturbation term. The choice of this data determines the chain complex and the differential. These can and will differ with differing choices of (â, $\hat{J}$, r, $\mathfrak{g}$). (In this case, it is known that the respective cohomology/homology groups as defined by any two choices of data are isomorphic.)

With the preceding understood, remark that Theorem 1 refers to isomorphisms that are constructed using maps between respective (co)chain complexes as defined by a given choice of the data (a, J) on the contact homology side, and a corresponding choice (â, $\hat{J}$, r, $\mathfrak{g}$) on the Seiberg-Witten side. In particular, it is necessary to first choose suitable pairs (a, J) and (â, $\hat{J}$, r, $\mathfrak{g}$) so as to embark on a proof of these theorems. (As it turns out, the obvious choice a = â and J = $\hat{J}$ is not, in general, the most useful.)

Granted all of the above, fix for the moment a pair (a, J) where a comes from Lemma 2.1's residual set and where J ∈ $\mathcal{J}_a$. This pair is used to define the embedded contact homology complex $\mathcal{C}_*$ and the associated differential. Meanwhile, the definition of the $\mathbb{Z}$-module for the Seiberg-Witten Floer cochain complex can be made for suitable choices of (â, $\hat{J}$), and then r ≥ 1 and a coclosed 1-form $\mu$ from the space $\Omega$ so as to define (3.17). The following proposition can be used to choose such pairs (r, $\mu$).

**Proposition 4.1**: *Fix a pair* (â, $\hat{J}$) *of contact 1-form and almost complex structure from* $\mathcal{J}_{\hat{a}}$. *Use them to define the metric on* M *and use* â *as the contact form in the corresponding version of (3.17). If* $c_1(\det(\mathbb{S}))$ *is a torsion class also fix* k ∈ $\mathbb{Z}$. *There exists* $r_* ≥ 1$ *and a residual set of coclosed 1-forms in* $\Omega$ *with the following properties: Suppose that* $\mu$ *is from this set and suppose that* $\mu$ *has* $\mathcal{P}$ *norm less than 1. Fix* r ∈ $[r_*, \infty)$ *and use* (â, $\hat{J}$) *with* r *and* $\mu$ *to define (3.17). There is a countable, non-accumulating set* $\mathcal{U} \subset [r_*, \infty)$ *such that if* r ∉ $\mathcal{U}$, *then*

- *If* $c_1(\det(\mathbb{S}))$ *is a torsion class, then all solutions to (3.17) with degree* k *or greater are irreducible.*

- *If $c_1(\det(\mathbb{S}))$ is not torsion, and if $\mathfrak{c}$ is a solution to (3.17) then the corresponding operator $\mathfrak{L}_\mathfrak{c}$ has trivial cokernel. If $c_1(\det(\mathbb{S}))$ is torsion, then such is the case if $\mathfrak{c}$ has degree* k *or greater.*
- *If $c_1(\det(\mathbb{S}))$ is not torsion, and if $\mathfrak{c}$ and $\mathfrak{c}'$ are gauge inequivalent solutions to (3.17) then $\mathfrak{a}(\mathfrak{c}) \neq \mathfrak{a}(\mathfrak{c}')$. If $c_1(\det(\mathbb{S}))$ is torsion, then such is the case if $\mathfrak{c}$ and $\mathfrak{c}'$ have degree* k *or greater.*

***Proof of Proposition 4.1***: This is Proposition 3.11 in [T3] and Proposition 2.1 in [T4].

Given the pair $(\hat{a}, \hat{J})$, a 1-form $\mu \in \Omega$, and then $r \geq 1$, use $\mathcal{M}^r$ in what follows to denote the set of gauge equivalence classes of solutions of solutions to (3.17). If $c_1(\det(\mathbb{S}))$ is torsion, then fix $k \in \mathbb{Z}$.

With $(\hat{a}, \hat{J})$ given, fix a 1-form $\mu$ as described in Proposition 4.1 and then fix $r$ in $[r_*, \infty) - \mathcal{U}$ so as to define the set $\mathcal{M}^r$. As explained in Section 3c of [T3] and in Section 2a of [T4], the elements in $\mathcal{M}^r$ can be used to label the generators of the Seiberg-Witten Floer complex when $c_1(\det(\mathbb{S}))$ is not torsion; and those with degree k or greater can be used to label the generators in degrees k or greater in the case when $c_1(\det(\mathbb{S}))$ is torsion. This being the case, they also label the generators of the dual cochain complex.

With (a, J) chosen, there is a concrete $\mathbb{Z}$ module for the embedded contact homology chain complex. Likewise, with $(\hat{a}, \hat{J})$, r and $\mathfrak{g} = \mathfrak{e}_\mu$ chosen as above, there is a concrete $\mathbb{Z}$ module for the Seiberg-Witten Floer cohomology complex, or at least the degree k or greater portion when $c_1(\det(\mathbb{S}))$ is torsion.

**a) From Reeb orbits to Seiberg-Witten solutions on M**

Suppose that $(\hat{a}, \hat{J})$ is a pair of contact structure from Lemma 2.1's residual set and almost complex structure in $\mathcal{J}_{\hat{a}}$. Fix L. The chain complex $\mathcal{C}_{\text{ech}}^L$ has generators that are equivalence classes of pairs of the form $(\Theta, \mathfrak{o})$ where $\Theta$ is a set whose typical element is a pair of Reeb orbit and positive integer subject to various constraints. It proves useful to introduce now $\mathcal{Z}_{\text{ech}}^L$ to denote the set of such $\Theta$. Assume now the following:

- *There is no element $\Theta \in \mathcal{Z}_{\text{ech}}^L$ with $\sum_{(\gamma,m) \in \Theta} m \ell_\gamma = L$*
- *Suppose that $\gamma$ is a Reeb orbit with $\ell_\gamma < L$. Then $\gamma$ has a tubular neighborhood map $\varphi \colon S^1 \times D \to M$ as described in Section 2a such that if $\gamma$ is hyperbolic with rotation number* k, *then $(\nu, \mu) = (\frac{1}{4}k, i\varepsilon e^{ikt})$ with $\varepsilon > 0$ but very small. Meanwhile, if $\gamma$ is hyperbolic, then its rotation number* R *is irrational. Furthermore,*
  i) *The pair $(\nu, \mu) = (\frac{1}{2}R, 0)$.*

ii) *The φ\*-pull back of* $T^{1,0}(\mathbb{R} \times M)$ *is spanned by* $ds + ia$ *and* $\frac{\ell_\gamma}{2\pi}(dz - iRz\,dt)$.
*Moreover, these two forms are orthogonal and have norm* $\sqrt{2}$.

(4.1)

The theorem that follows asserts the existence of a map from $\mathcal{Z}_{ech}{}^L$ into the version of $\mathcal{M}^r$ that is defined using (â, Ĵ), a sufficiently large r, and a reasonably small 1-form $\mu$. This map is used to define the isomorphism for Theorem 1.

**Theorem 4.2**: *Fix* $L \geq 1$ *and a pair* (â, Ĵ) *as above that obeys (4.1). There exists* $\kappa \geq 1$ *with the following significance: Define* $\mathcal{M}^r$ *using* $r \geq \kappa$ *and a 1-form* $\mu \in \Omega$ *with* $\mathcal{P}$ *norm bounded by* 1. *There exists a map* $\Phi^r: \mathcal{Z}_{ech}{}^L \to \mathcal{M}^r$ *with the three properties listed below.*

- $\Phi^r$ *is a bijection onto the subset in* $\mathcal{M}^r$ *of elements* $\mathfrak{c} = (A, \psi)$ *with* $\mathcal{E}(A) < 2\pi L$.
- *If* $\mathfrak{c} \in \mathcal{M}^r$ *is in the image of* $\Phi^r$, *then the operator* $\mathfrak{L}_\mathfrak{c}$ *has trivial kernel.*
- *Let* $\Theta$ *and* $\Theta'$ *denote any two elements in* $\Lambda^L$ *and let* $z$ *and* $z'$ *denote their respective* $\mathbb{Z}/p\mathbb{Z}$ *degrees. Meanwhile, let* $x$ *and* $x'$ *denote the respective* $\mathbb{Z}/p\mathbb{Z}$ *degrees of* $\Phi^r(\Theta)$ *and* $\Phi^r(\Theta')$. *Then* $x - x' = -(z - z')$ *modulo* $p\mathbb{Z}$.

The upcoming Section 5 gives some idea of what $\Phi^r$ looks like. The actual construction of $\Phi^r$ is in Paper II of this series, [T6]. The assertion that it defines a bijection as described by the Theorem's first bullet is proved in Paper IV of the series, [T8]. The assertions of the second and third bullets are proved in Paper III of the series, [T9].

A rather more complicated version of Theorem 4.2 holds when the second item in (4.1) is not invoked. In the latter case, each $\Theta \in \mathcal{Z}_{ech}{}^L$ parametrizes a subset of $\mathcal{M}^r$ such that the collection of these subsets partition the of elements $\mathfrak{c} = (A, \psi)$ with $\mathcal{E}(A) < 2\pi L$. More is said in Section 5 about this more general version of Theorem 4.2

**b) From pseudoholomorphic curves to Seiberg-Witten solutions on** $\mathbb{R} \times M$

Fix $L \geq 1$ and then fix $\mu \in \Omega$ with $\mathcal{P}$ norm less than 1. Any given large r version of Theorem 4.2's map $\Phi^r$ identifies $\mathcal{Z}_{ech}{}^L$ with a subset in $\mathcal{M}^r$. This being the case, $\Phi^r$ can be used to define a monomorphism from $\mathcal{C}_{ech}{}^L$ into the Seiberg-Witten cochain complex. This is done as follows: Let $\Theta \in \mathcal{Z}_{ech}{}^L$. Order the subset of pairs $(\gamma, 1) \in \Theta$ for which $\gamma$ is hyperbolic with even rotation number. Doing so for all such $\Theta$ identifies $\mathcal{Z}_{ech}{}^L$ with a set of generators of $\mathcal{C}_{ech}{}^L$. The image of $\mathcal{Z}_{ech}{}^L$ via $\Phi^r$ defines a set of generators of $\mathcal{C}^{SW}$. Extend this map of generators in a $\mathbb{Z}$-linear fashion. The monomorphisms so constructed is canonical up to precomposing by an isomorphism from $\mathcal{C}_{ech}{}^L$ to itself that changes the sign of some of the generators. A monomorphism that is obtained from $\Phi^r$ in this way is denoted by $T_\Phi$.

The next theorem asserts in part that there is a choice for $T_\Phi$ that intertwines the embedded contact homology differential with the Seiberg-Witten Floer cohomology differential. This theorem reintroduces the space $\mathcal{M}_1(\Theta_-, \Theta_+)$ from Section 2c and the space $\mathcal{M}_1(\mathfrak{c}_-, \mathfrak{c}_+)$ from Section 3b. Note that both spaces admit a canonical $\mathbb{R}$-action, this induced by the action of $\mathbb{R}$ on $\mathbb{R} \times M$ as the group of constant translations of the $\mathbb{R}$ factor in $\mathbb{R} \times M$.

**Theorem 4.3**: *Fix $L \geq 1$ and a pair $(\hat{a}, \hat{J})$ as above that obeys (4.1). There exists $\kappa \geq 1$ with the following significance: Define $\mathcal{M}^r$ using $r \geq \kappa$ and a 1-form $\mu \in \Omega$ with $\mathcal{P}$ norm bounded by 1.*

- *Let $\Theta_-$ and $\Theta_+$ denote any two elements in $\mathcal{Z}_{ech}^L$. Use $\mathfrak{c}_-$ and $\mathfrak{c}_+$ to denote solutions to (3.17) whose gauge equivalences classes are the respective images in $\mathcal{M}^r$ of $\Theta_-$ and $\Theta_+$ via Theorem 4.2's map $\Phi^r$.*
  i) *The space $\mathcal{M}_1(\mathfrak{c}_-, \mathfrak{c}_+)$ has a finite set of components, and each component is an orbit of the canonical $\mathbb{R}$ action. In addition, if $\mathfrak{d} \in \mathcal{M}_1(\mathfrak{c}_-, \mathfrak{c}_+)$, then the corresponding $\mathfrak{g} = \mathfrak{e}_\mu$ version of (3.9) has trivial cokernel.*
  ii) *There is an $\mathbb{R}$-equivariant diffeomorphism, $\Psi^r$, from $\mathcal{M}_1(\Theta_-, \Theta_+)$ to $\mathcal{M}_1(\mathfrak{c}_-, \mathfrak{c}_+)$.*
- *There is a choice for $T_\Phi$ such that if $\Theta_-$ and $\Theta_+$ now denote any given pair of generators of $\mathcal{C}_{ech}^L$, then the contribution, +1 or -1, of any given component in $\mathcal{M}_1(\Theta_-, \Theta_+)$ to the embedded contact homology differential is the same as the contribution of its $\Psi^r$ image to the Seiberg-Witten Floer cohomology differential.*

A rough picture of $\Psi^r$ is given in Section 5. The full construction is in Paper II of this series, [T6]. The proof that $\Psi^r$ is a bijection is given in Paper IV, [T8]. Meanwhile, the proof of that $\Psi^r$ is an $\mathbb{R}$-equivariant diffeomorphism is given in Paper III, [T3]; as is the proof of the assertion in the theorem's second bullet. Theorem 4.3 has a replacement of sorts when the second item in (4.1) is not present. The latter is vastly more complicated to state, let alone prove. More is said on this score in Section 5.

Theorem 4.3 implies that there is a version of the monomorphism $T_\Phi$ that intertwines the embedded contact homology differential with the differential that defines the Seiberg-Witten Floer cohomology. As such, this $T_\Phi$ identifies $\mathcal{C}_{ech}^L$ as a subcomplex in the Seiberg-Witten Floer cochain complex. The assertion made by the next theorem implies that any sufficiently large r version of $T_\Phi(\mathcal{C}_{ech}^L)$ is mapped to itself by the Seiberg-Witten Floer cohomology differential.

**Theorem 4.4**: *Fix $L \geq 1$ and a pair $(\hat{a}, \hat{J})$ as above that obeys (4.1). There exists $\kappa \geq 1$ with the following significance: Define $\mathcal{M}^r$ using $r \geq \kappa$ and a 1-form $\mu \in \Omega$ with $\mathcal{P}$ norm* *bounded by* 1. *Suppose that* $c_+ \in \Phi^r(\mathcal{Z}_{ech}{}^L)$ *and that* $c_- \in \mathcal{M}^r$ *is such that* $\mathcal{M}_1(c_-, c_+) \neq \emptyset$. *Then* $c_-$ *is also in* $\Phi^r(\mathcal{Z}_{ech}{}^L)$.

Theorem 4.4 is proved in Section 4h, below

**c) Proof of Theorem 1**

Fix a pair (a, J) of contact 1-form from Lemma 2.1's residual set and almost complex structure from $\mathcal{J}_a$. Use this data to define the embedded contact homology chain $\mathbb{Z}$-module $\mathcal{C}_{ech}$ and its differential. The latter has a filtration $\{\mathcal{C}_{ech}{}^L\}_{L\geq 1}$ with the corresponding homology groups. Fix $L \geq 1$ such that the top item in (4.1) holds for the pair (a, J). Fix small $\delta > 0$ and Proposition 2.5 supplies a $(\delta, L)$ approximation, $(\hat{a}, \hat{J})$, to (a, J). The latter defines the analogous set of $\mathbb{Z}$-modules $\{\hat{\mathcal{C}}_{ech}{}^{L'}\}_{L'\geq 1}$ and corresponding homology groups. Let $\hat{\mathcal{Z}}_{ech}{}^L$ denote $\hat{\mathcal{C}}_{ech}{}^L$ analog of $\mathcal{Z}_{ech}{}^L$. In this regard, the elements in $\hat{\mathcal{Z}}_{ech}{}^L$ are geometrically identical to those $\mathcal{Z}_{ech}{}^L$. As noted in Proposition 2.4, the corresponding versions of $\mathcal{ECH}^L$ are canonically isomorphic.

Theorem 1 requires the following additional observations:

**Theorem 4.5**: *Fix a pair* (a, J) *where* a *is a contact 1-form from Lemma 2.1's residual set, and where* $J \in \mathcal{J}_a$. *Fix* $L_* \geq 1$ *such that the* (a, J) *and* $L_*$ *version of the first item in (4.1) holds. Also, fix* $k \in \mathbb{Z}$ *when* $c_1(\det(\mathbb{S}))$ *is a torsion class. There exists* $\kappa \geq 1$ *with the following significance: Fix* $L \geq \kappa$ *such that the* (a, J) *and* L *version of the first item in (4.1) holds. Fix* $\delta \in (0, \kappa^{-1})$; *and then choose a* $(\delta, L)$ *approximation*, $(\hat{a}, \hat{J})$, *to the pair* (a, J). *Choose a 1-form* $\mu \in \Omega$ *with* $\mathcal{P}$ *norm bounded by 1 as described in Proposition 4.1. Then take* r *very large so as to define* $\mathcal{M}^r$ *and* $\Phi^r$ *using the data* $(\hat{a}, \hat{J})$, r *and* $\mu$. *Likewise, use* $T_\Phi$ *to define the* $(\hat{a}, \hat{J})$ *version of the monomorphism given by the second bullet in Theorem 4.3. Let* $\theta$ *denote any given Seiberg-Witten Floer cohomology class, but of degree* k *or greater if* $c_1(\det(\mathbb{S}))$ *is a torsion class.*

- $\theta$ *has a representative cocycle that lies in* $T_\Phi(\hat{\mathcal{C}}_*{}^L)$.
- *Suppose that* $\upsilon \in \hat{\mathcal{C}}_*{}^{L_*}$ *and that* $T_\Phi(\upsilon)$ *is a coboundary. If* $c_1(\det(\mathbb{S}))$ *is a torsion class, also assume that* $\upsilon$ *has degree* -k *or less. Then* $\upsilon$ *is a boundary in* $\hat{\mathcal{C}}_*{}^L$.

The proof of Theorem 4.5 uses, for the most part, analysis that was introduced in [T3] and [T4]. This theorem is also proved in Section 4h.

*Proof of Theorem 1*: If $c_1(\det(\mathbb{S}))$ is torsion, fix $k \in \mathbb{Z}$. Fix $(\delta, L)$ and then $(\hat{a}, \hat{J})$ as in Theorem 4.5. Use the latter with $\mu \in \mathcal{P}$ and a sufficiently large choice of r to define $\mathcal{M}^r$ and $\Phi^r$. Let $T_\Phi$ denote the monomorphism given in the second bullet of the $(\hat{a}, \hat{J})$ version of Theorem 4.3. According to what is said in Theorem 4.4, if r is large enough, then $T_\Phi(\hat{\mathcal{C}}_*^L)$, or at least the degree k or greater part if $c_1(\det(\mathbb{S}))$ is torsion, is a submodule of the Seiberg-Witten Floer cochain complex that is mapped to itself by the latter's differential. According to what is said in the first bullet of Theorem 4.5, the induced map from the cohomology of this submodule into the full Seiberg-Witten Floer cohomology is surjective. Given what is said by Proposition 2.4, this implies that the Seiberg-Witten Floer cohomology can be viewed as a summand inside the version of $\mathcal{ECH}^L$ that is defined by (a, J). The second bullet with Proposition 2.4 implies that this summand is all that remains when the direct limit is taken to define the embedded contact homology

**d) Filtrations of the Seiberg-Witten Floer cohomology**

This section with Sections 4e-g supply the background material that is used subsequently to prove Theorems 4.4 and 4.5. To begin, suppose that $\mathcal{A} \in \mathbb{R}$ has been given. Let (a, J) denote a pair consisting of a contact 1-form and compatible almost complex structure. Use the latter to define the metric on M and the corresponding version of (3.17). Fix a form $\mu \in \Omega$ with $\mathcal{P}$ norm less than 1 as in Proposition 4.1. If $c_1(\det(\mathbb{S}))$ is a torsion class, also fix $k \in \mathbb{Z}$. Now suppose that there exists $r_\mathcal{A} \in [1, \infty)$ with the following property:

*Suppose that $r \geq r_\mathcal{A}$. If $c_1(\det(\mathbb{S}))$ is not torsion, then there are no solutions to (3.17) with $\mathfrak{a}^f = r\mathcal{A}$. If $c_1(\det(\mathbb{S}))$ is torsion, then there are no solutions to (3.17) with degree k or greater with $\mathfrak{a} = r\mathcal{A}$.*

(4.2)

Fix $r \in [r_* + r_\mathcal{A}, \infty) - \mathcal{U}$ so as to define the Seiberg-Witten Floer cochain complex. Given that $\delta^*$ increases $\mathfrak{a}^f$, this cochain complex has as a subcomplex the set $\mathcal{C}^{SW,\mathcal{A}}$ that is generated by the elements in $\mathcal{M}^r$ with $\mathfrak{a}^f > r\mathcal{A}$ or $\mathfrak{a} \geq r\mathcal{A}$ as the case may be. The short exact sequence

$$0 \to \mathcal{C}^{SW,\mathcal{A}} \to \mathcal{C}^{SW} \to \mathcal{C}^{SW}/\mathcal{C}^{SW,\mathcal{A}} \to 0$$

(4.3)

induces a corresponding long exact sequence in cohomology, thus

$$\cdots \to H^j(\mathcal{C}^{SW,\mathcal{A}}) \to H^j(\mathcal{C}^{SW}) \to H^j(\mathcal{C}^{SW}/\mathcal{C}^{SW,\mathcal{A}}) \to H^{j+1}(\mathcal{C}^{SW,\mathcal{A}}) \to \cdots$$

(4.4)

where the restriction $j \geq k$ is assumed in the case when $c_1(\det(\mathbb{S}))$ is torsion. To elaborate, recall from Section 3 of [T3] and Section 2 of [T4] that with r fixed, the differential on $\mathcal{C}^{SW}$ is defined from the instanton solutions to the version of (3.6) that is defined by taking the function $\mathfrak{g}$ in (3.15) to equal $\mathfrak{e}_\mu + \mathfrak{p}$ where $\mathfrak{p}$ is a suitably generic element in $\mathcal{P}$ with very small norm. In particular, there will be no solutions with $\mathfrak{a}^f = r\mathcal{A}$ or $\mathfrak{a} = r\mathcal{A}$ to the version of (3.5) that is defined by $\mathfrak{p}$ if the latter is chosen from a sufficiently small radius ball about the origin in $\mathcal{P}$. A sufficiently generic choice from this ball can then be used to define the differential on $\mathcal{C}^{SW}$. This differential will preserve $\mathcal{C}^{SW,\mathcal{A}}$ and then give the resulting exact sequence in (4.4). The same argument used in [KM] to prove the $\mathfrak{p}$-independence of $H^k(\mathcal{C}^{SW})$ proves that each of the terms in (4.4) is also $\mathfrak{p}$-independent if $\mathfrak{p}$ is chosen from a sufficiently small radius ball about the origin in $\mathcal{P}$.

To say more about this, note that distinct versions are compared by choosing a suitable path of perturbations that interpolate from one to the other. If the path is chosen so that each element has small norm, then none of the corresponding versions of (3.17) will have solutions with $\mathfrak{a}^f = r\mathcal{A}$ or $\mathfrak{a} = r\mathcal{A}$. Moreover, if the path is chosen in a sufficiently generic fashion, then there will be only a finite set of points along the path where either the complex in (4.3) can not be defined, or where the differential can not be defined (this follows from what is done in Section 7 of [T3].) Fix any two points in any given interval in the complement of this bad set. Then the corresponding complexes are canonically isomorphic via an isomorphism that preserves the differentials (see Sections 3 in [T3] and Section 2 in [T4].) A comparison between the versions of (4.4) that are defined on consecutive intervals uses cobordisms to construct a chain equivalence. The resulting chain maps need not increase $\mathfrak{a}^f$ or $\mathfrak{a}$, but in any event, they can be constructed so that if $\mathfrak{a}^f$ or $\mathfrak{a}$ is not increased, then it is decreased by a very small amount. In fact, the amount of decrease can be assumed to be as desired (see, the arguments in Section 7 of [T3]). This understood, these chain maps induce the desired equivalence between the corresponding versions of (4.4).

The fact that the terms in (4.4) and the arrows do not depend on the perturbation term $\mathfrak{p}$ if the latter has small norm implies that the terms and homomorphisms in (4.4) are also independent of r for $r > r_{\mathcal{A}}$. This perturbation independence also holds for perturbations of (a, J) and $\mu$ that differ from the latter by terms with sufficiently small norm.

To make a precise statement of this last point for use below, suppose now that $\{(a_\tau, J_\tau)\}_{\tau \in [0,1]}$ is a smoothly parametrized family with each term consisting of a contact 1-form and compatible almost contact structure. Let $\{\mu_\tau\}_{\tau \in [0,1]}$ denote a corresponding family in $\Omega$ such that each member has $\mathcal{P}$ norm less than 1. Assume that $(a_0, J_0) = (a, J)$ and that $\mu_0 = \mu$. Assume in addition that $\mu_1$ is chosen from the $(A_1, J_1)$ version of Proposition 4.1's residual subset. Finally, assume that the assertion made in (4.2) holds when (a, J) and $\mu$ are replaced by any given $(a_\tau, J_\tau)$ and $\mu_\tau$.

As (4.2) holds for $(a_1, J_1)$ and $\mu_1$, there is a corresponding version of (4.3) and also (4.4). The next lemma asserts that the original version and the new version of (4.4) are equivalent. This lemma uses $r_1$ and $r_2$ to denote the respective $\tau = 0$ and $\tau = 1$ versions of the constant $r_*$ that appears in Proposition 4.1; and it uses $\mathcal{U}_0$ and $\mathcal{U}_1$ to denote the corresponding versions of the set $\mathcal{U}$ that is described in Proposition 4.1.

**Lemma 4.6**: *Suppose that $\{(a_\tau, J_\tau), \mu_\tau\}_{\tau \in [0,1]}$ has all of the properties just described. Fix $r \in [r_1 + r_2 + r_{\mathcal{A}}, \infty) - (\mathcal{U}_0 \cup \mathcal{U}_1)$. There exist isomorphisms between corresponding terms in the respective $\tau = 0$ and $\tau = 1$ versions of (4.2) that intertwine the homomorphisms in these exact sequences.*

*Proof of Lemma 4.6*: To start, the assumptions imply the existence of a ball $\mathcal{B} \subset \mathcal{P}$ about the origion of small radius such that if $\mathfrak{p} \in \mathcal{B}$, then the following is true: Fix $\tau \in [0, 1]$ and use $(a_\tau, J_\tau)$, the given value for r and the perturbation term $\mathfrak{g} = \mathfrak{e}_{\mu_\tau} + \mathfrak{p}$ to define (3.5). Then there will be no solutions to the latter with $\mathfrak{a}^f = r\mathcal{A}$ or $\mathfrak{a} = r\mathcal{A}$ as the case may be with $\mathfrak{a}$ defined as in (3.15) using either $\mathfrak{g}$ or just $\mathfrak{e}_{\mu_\tau}$. If $\mathfrak{p}$ is chosen in a sufficiently generic manner, then the corresponding solutions to (3.5) can be used to define the Seiberg-Witten cochain complex and the solutions to (3.6) will define the the differential on this complex. Since the differential increases $\mathfrak{a}^f$ and $\mathfrak{a}$, it preserves the corresponding version of (4.3) and so there is a corresponding version of (4.4). The respective versions of (4.4) that are defined by two suitably generic choices from $\mathcal{B}$ will be isomorphic.

Now consider comparing the distinct $\tau$ and $\tau'$ versions of (4.4). The story here is complicated slightly by the fact that changes in $\tau$ change $(a_\tau, J_\tau)$ and thus the metric on M. Such changes do not arise as perturbations from $\mathcal{P}$. To circumvent this point, suppose that $\mathfrak{p}$ is chosen so that all solutions to the $(a_\tau, J_\tau)$, r and $\mathfrak{g} = \mathfrak{e}_{\mu_\tau} + \mathfrak{p}$ version of (3.5) obey (3.7) and such that (3.8) is also obeyed. If $\varepsilon > 0$, is sufficiently small, then such will be the case for all $(a_{\tau'}, J_{\tau'})$, r and $\mathfrak{g} = \mathfrak{e}_{\mu_{\tau'}} + \mathfrak{p}$ solutions as well if $|\tau - \tau'| < \varepsilon$. As a consequence, there is a canonical identification between the corresponding $\tau$ and $\tau'$ versions of the terms in (4.3) that intertwines the associated differentials. Hence, the corresponding $\tau$ and $\tau'$ versions (4.4) are isomorphic.

**e) Min-max**

This section introduces various additional notions that were used in [T3], [T4]. Minor modification are made on these in preparation for their use in the upcoming proof of Theorem 4.5.

To set the stage, fix a pair (a, J) of contact 1-form and compatible almost complex structure. If $c_1(\det(\mathbb{S}))$ is torsion, also fix an integer k. Let $\mu \in \Omega$ denote a 1-form with $\mathcal{P}$ norm less than 1 of the sort described by Proposition 4.1. Use the (a, J), $\mu$ and $r \gg 1$ to define (3.17) and the space $\mathcal{M}^r$.

Suppose that $\theta$ is a non-zero Seiberg-Witten Floer cohomology class in a given degree, this k or greater if $c_1(\det(\mathbb{S}))$ is torsion. Let $\mathfrak{n} = \sum_\mathfrak{c} z_\mathfrak{c} \mathfrak{c}$ denote a representative of this class in the Seiberg-Witten Floer cochain complex. Here, $z_\mathfrak{c} \in \mathbb{Z}$ and $\mathfrak{c} \in \mathcal{M}^r$. Define $\mathfrak{a}^f[\mathfrak{n}, r]$ to be the minimum of the values of $\mathfrak{a}^f$ on the set $\{\mathfrak{c}: z_\mathfrak{c} \neq 0\}$. Set $\mathfrak{a}^f_\theta[r]$ to denote the maximum of $\{\mathfrak{a}^f[\mathfrak{n}, r]: \mathfrak{n} \text{ represents } \theta\}$. Note the $\mathfrak{a}^f_\theta[\cdot]$ is defined by first taking a minimum and then a maximum; and that this is opposite to the order used in [T3] and [T4]. The order is switched because $\theta$ is a *cohomology* class rather than a *homology* class.

**Proposition 4.7:** *Fix a pair* (a, J), *and fix* $k \in \mathbb{Z}$ *if* $c_1(\det(\mathbb{S}))$ *is a torsion class. Fix* $\mu \in \Omega$ *with* $\mathcal{P}$ *norm less than 1 as described by Proposition 4.1. For each* $r \in [r_*, \infty) - \mathcal{U}$ *use the data* (a, J), $\mu$ *and* r *to define the corresponding version of the Seiberg-Witten Floer cochain complex and the associated cohomology. These various versions of the Seiberg-Witten Floer cohomology groups (in degrees k or greater when* $c_1(\det(\mathbb{S}))$ *is torsion) can be identified so that the following is true: If* $\theta$ *is any given non-zero cohomology class, then the assignment* $r \to \mathfrak{a}^f_\theta[r]$ *is the restriction of continuous, piecewise differentiable function on the half line* $[r_*, \infty)$.

*Proof of Proposition 4.7*: Except for the ordering change with regards to 'min' and 'max', the argument is the same as that used to prove Proposition 4.2 in [T3] when $c_1(\det(\mathbb{S}))$ is torsion. But for this same ordering change, the argument is essentially that used for Proposition 2.5 in [T4] when $c_1(\det(\mathbb{S}))$ is not torsion. With regards to the ordering of min and max, note that if $\mathfrak{c}$ is a given generator from $\mathcal{M}^r$, then a generator $\mathfrak{c}'$ appears on the right hand side of (3.3)'s definition of $\delta^*\mathfrak{c}$ only if $\mathcal{M}_1(\mathfrak{c}', \mathfrak{c}) \neq \emptyset$. As such, $\delta^*\mathfrak{c}$ is a sum of generators on which $\mathfrak{a}^f(\cdot) > \mathfrak{a}^f(\mathfrak{c})$.

Now fix $\mathcal{A} \in \mathbb{R}$ and assume that (4.2) holds. A similar min-max construction can be done for classes in $H^*(\mathcal{C}^{SW, \mathcal{A}})$ and also for those in $H^*(\mathcal{C}^{SW}/\mathcal{C}^{SW, \mathcal{A}})$. To say more, suppose that $\theta$ is a non-zero class in either of these groups. Assume that $\theta$ has degree k or greater in the case when $c_1(\det(\mathbb{S}))$ is torsion. Given a cocycle $\mathfrak{n}$ that represents $\theta$, define $\mathfrak{a}^f[\mathfrak{n}, r]$ to be the mimimum of $\mathfrak{a}^f(\cdot)$ on the generators that represent $\theta$. Then define $\mathfrak{a}^f_\theta[r]$ to be the maximum of the elements in $\{\mathfrak{a}^f[\mathfrak{n}, r]: \mathfrak{n} \text{ represents } \theta\}$. The analog of Proposition 4.7 in this new context is given by the next proposition. It's proof is essentially the same as that for Proposition 4.7 so the details are left to the reader.

**Proposition 4.8:** *Fix a pair* (a, J), *and fix* $k \in \mathbb{Z}$ *if* $c_1(\det(\mathbb{S}))$ *is a torsion class. Choose a form* $\mu \in \Omega$ *with* $\mathcal{P}$ *norm less than 1 as described by Proposition 4.1. Fix* $\mathcal{A} \in \mathbb{R}$ *and suppose that (4.2) holds. For each* $r \in [r_* + r_\mathcal{A}, \infty) - \mathcal{U}$, *use the data* (a, J), $\mu$ *and* r *to* *define the corresponding versions of (4.3) and (4.4). As* r *varies, the corresponding versions of (4.4) can be identifed (in degrees* k *or greater when* $c_1(\det(\mathbb{S}))$ *is torsion) so that if* θ *denotes any given cohomology class in any of the three cohomology groups, then the assignment* r → $\mathfrak{a}^f_\theta[r]$ *is the restriction of continuous, piecewise differentiable function on the half line* $[r_*, \infty)$.

The continuity and piecewise differentiability of $\mathfrak{a}^f_\theta[\cdot]$ is exploited in the next subsection.

### f) Bounds on E from $\mathfrak{a}^f$ and vice versa.

The next proposition plays one of the key roles in the proof of Theorem 4.5. To set the terminology, suppose that (a, J) is given to define the metric on M. Fix $\mu \in \Omega$ with $\mathcal{P}$ norm less than 1 as described in Proposition 4.1. Suppose in addition that $\mathcal{A} \in \mathbb{R}$ has been specified and that (4.2) holds. The various $r \in [r_* + r_\mathcal{A}, \infty) - \mathcal{U}$ versions of (4.4) are implicitly identified using one of the identifications provided by Proposition 4.8. Such an identification should be understood when reference is made to a class in one of the groups in (4.4) with no reference to the precise value of r.

**Proposition 4.9**: *Fix a pair* (a, J), *and if* $c_1(\det(\mathbb{S}))$ *is torsion, fix an integer* k. *Choose* $\mu \in \Omega$ *with* $\mathcal{P}$ *norm less than 1 as described in Proposition 4.1. Fix* $\mathcal{A} \in \mathbb{R}$ *and suppose that (4.2) holds. There exists* $\mathcal{K} \geq 1$ *with the following significance: Suppose that* θ *is a non-zero cohomology class of fixed degree (k or greater if* $c_1(\det(\mathbb{S}))$ *is torsion) in any of the cohomology groups that appear in (4.4). There exists an increasing, unbounded sequence* $\{r_i\}_{i=1,2,...}$ *in* $[r_* + r_\mathcal{A}, \infty) - \mathcal{U}$ *such that the* $r = r_i$ *version of the class* θ *has a representative cocycle with* $E(\cdot) \leq 2\pi\mathcal{K}$ *on each generator that appears with a non-zero coefficient.*

The proof of this last proposition makes use of the following restatement of the following results from [T3] and [T4].

**Proposition 4.10**: *Let* (a, J) *denote a pair consisting of a contact 1-form and compatible almost complex structure. There exists* $\kappa \geq 1$ *such that the following is true: Fix* $\mu \in \Omega$ *with* $\mathcal{P}$ *norm bounded by 1. Suppose that* $\mathfrak{c} = (A, \psi)$ *is a solution to the* (a, J), $\mu$ *and* r *version of (3.17). Then*
- $|\mathfrak{cs}^f| < \kappa r^{31/16}$.
- *If* $c_1(\det(\mathbb{S}))$ *is torsion, then* $|\mathfrak{cs}| \leq \kappa r^{2/3}(1 + E^{4/3})$.
- *If* $c_1(\det(\mathbb{S}))$ *is not torsion, then* $|\mathfrak{cs}^f| \leq \kappa\, r^{2/3}(\ln r)^\kappa (1 + E^{4/3})$

*Proof of Proposition 4.9*: But for notation, the arguments in the case $c_1(\det(\mathbb{S}))$ is torsion are those used in Propositions 4.6 and Corollary 4.7 in [T3] with Propositions 4.8 and 4.10 added. In the case that $c_1(\det(\mathbb{S}))$ is not torsion, the arguments are the same but for notation as those used in Section 2c of [T4] but with Propositions 4.8 and 4.10 added.

*Proof of Proposition 4.10*: The bound in the first bullet follows from Proposition 5.1 in [T3] in the case when $c_1(\det(\mathbb{S}))$ is torsion, and Proposition 1.10 of [T4] when $c_1(\det(\mathbb{S}))$ is not torsion. The bound in the second bullet follows from (4.2) in [T3 and Lemma 2.4 in [T3]. The bound in the third bullet restates Proposition 1.9 of [T4].

Proposition 4.10 has an additional very important corollary, this stated by

**Proposition 4.11**: *Let* (a, J) *denote a pair consisting of a contact 1-form and compatible almost complex structure. There exists $\kappa \geq 1$ such that the following is true: Fix $\mu \in \Omega$ with $\mathcal{P}$ norm bounded by* 1.

- *Suppose that $c_1(\det(\mathbb{S}))$ is torsion and that $\mathfrak{c} = (A, \psi)$ is a solution to the* (a, J), $\mu$ *and r version of (3.17) with $\mathfrak{a}(\mathfrak{c}) > -\kappa^{-1} r^{31/16}$. Then $|2r^{-1} \mathfrak{a}(\mathfrak{c}) + \mathrm{E}(\mathfrak{c})| \leq \kappa r^{-1/50} \mathrm{E}(\mathfrak{c})$.*
- *Suppose that $c_1(\det(\mathbb{S}))$ is not torsion and that $\mathfrak{c} = (A, \psi)$ is a solution to the* (a, J), *$\mu$ and r version of (3.17) with $\mathfrak{a}^f(\mathfrak{c}) > -\kappa^{-1} r^{31/16}$. Then $|2r^{-1} \mathfrak{a}^f(\mathfrak{c}) + \mathrm{E}(\mathfrak{c})| \leq \kappa r^{-1/50} \mathrm{E}(\mathfrak{c})$*

The proof of this proposition introduces a convention that is used throughout this paper and its sequels: In all appearances, '$c_0$' denotes a constant greater than 1 whose value is independent of r, $\mu$, and any given (A, $\psi$). Subsequent appearances of $c_0$ are allowed to have different values, but these can be assumed to increase from one appearance to the next.

*Proof of Proposition 4.11*: Since $|\mathfrak{cs}^f| \leq c_0 r^{31/16}$, and since it is assumed that $\mathfrak{a}$ or $\mathfrak{a}^f$ is greater than $-c_0 r^{31/16}$, it follows that $\mathrm{E} \leq c_0 r^{15/16}$. Hold this last bound for the moment. In the case when $c_1(\det(\mathbb{S}))$ is torsion, use the second bullet of Proposition 4.10 to see that $|\mathfrak{cs}| \leq c_0 r^{2/3} \mathrm{E}^{4/3}$. In the case when $c_1(\det(\mathbb{S}))$ is not torsion, use the third bulled in Proposition 4.10 to see that $|\mathfrak{cs}^f| \leq c_0 r^{2/3} \mathrm{E}^{4/3} (\ln r)^{1/c_0}$. These bounds and the bound just derived for $\mathrm{E}$ implies that $|\mathfrak{cs}|$ or $|\mathfrak{cs}^f|$ is no greater than $c_0 r^{-1/50} (r\mathrm{E})$. These imply the asserted bounds.

### g) Bounds on E and $\mathfrak{a}^f$ for families

The next result asserts a parameterized version of what is proved in Section 6d of [T3]. To set things up, suppose that $\{(a_\tau, J_\tau)\}_{\tau \in [0,1]}$ is a smoothly parametrized family of pairs consisting of a contact 1-form and compatible almost complex structure. This family is assumed now to have one additional attribute. Suppose that $L \geq 1$ has been specified such that the following condition holds for each $\tau \in [0, 1]$:

> *Let $\Theta$ denote a set of pairs of the form $(\gamma, m)$ with $\gamma$ Reeb orbit as defined by $a_\tau$ and m a positive integer. Assume that distinct pairs have distinct Reeb orbit component. Then $\sum_{(\gamma,m) \in \Theta} m \ell_\gamma \neq L$.*

(4.5)

Let $\{\mu_\tau\}_{\tau \in [0,1]}$ denote a corresponding family of 1-forms, each in $\Omega$ and each with $\mathcal{P}$ norm bounded by 1.

**Proposition 4.12**: *Given the data $\{(a_\tau, J_\tau), \mu_\tau\}_{\tau \in [0,1]}$ and L, there exists $\kappa$ with the following significance: Fix $\tau \in [0, 1]$. Suppose that $r \geq \kappa$ and that $\mathfrak{c} = (A, \psi)$ is a solution to the version of (3.17) that is defined by the data $(a_\tau, J_\tau)$, $\mu_\tau$ and r. Assume in addition that $E(A) \leq 2\pi L + \kappa^{-1}$. Then $E(A) < 2\pi L - \kappa^{-1}$.*

This last proposition leads to the following important observation.

**Proposition 4.13**: *Fix the data $\{(a_\tau, J_\tau), \mu_\tau\}_{\tau \in [0,1]}$ and L, and fix $k \in \mathbb{Z}$ if $c_1(\det(\mathbb{S}))$ is torsion. There exists $\kappa$ with the following significance: Fix $\tau \in [0, 1]$ and $r \geq \kappa$. Suppose that $\mathfrak{c} = (A, \psi)$ is a solution to the version of (3.17) that is defined by $(a_\tau, J_\tau)$, $\mu_\tau$ and r. Then $\mathfrak{a}^f(\mathfrak{c}) \neq -\pi L r$.*

*Proof of Proposition 4.12*: In what follows, the spinor bundle $\mathbb{S}$ is written as $E \oplus EK^{-1}$, and corresponding components of a given section $\psi$ are denoted by $(\alpha, \beta)$.

To start the proof, suppose that the proposition is not true. Then there exists an unbounded sequence $\{r_n\}_{n=1,2,\ldots} \subset [1, \infty)$, a convergent sequence $\{\tau_n\}_{n=1,2,\ldots} \in [0, 1]$, and a corresponding sequence $\{(A_n, \psi_n)\}_{n=1,2,\ldots}$ where any given $(A_n, \psi_n)$ obeys (3.17) as defined using $r = r_n$ and the $\tau = \tau_n$ version of $(a_\tau, J_\tau)$ and $\mu_\tau$. Let $\tau \in [0, 1]$ denote the limit point of the sequence $\{\tau_n\}_{n=1,2,\ldots}$. The four steps that follow explain why the existence of such a sequence leads to nonsense and so proves the lemma.

<u>Step 1</u>: The arguments used in Section 6d of [T3] can be repeated with only cosmetic changes to find a possibly empty set, $\Theta$, of the following sort: First, the typical element in $\Theta$ is a pair $(\gamma, m)$ where $\gamma$ is an $a_\tau$ Reeb orbit and m is a positive integer.

Moreover, distinct pairs have distinct Reeb orbit components. Finally, there exists some subsequence of $\{(A_n, \psi_n = (\alpha_n, \beta_n))\}_{n=1,2,...}$, hence renumbered consecutively from 1, and there is a sequence $\{\varepsilon_n\}_{n=1,2,...} \subset (0, \frac{1}{100}]$ with limit zero such that

- $|\alpha_n|^2 \geq 1 - \varepsilon_n$ *at all points with distance* $\varepsilon_n$ *or greater from* $\cup_{(\gamma,m)\in\Theta} \gamma$.
- *Fix* $(\gamma, m) \in \Theta$ *and let* $\varphi: S^1 \times D \to M$ *denote a coordinate chart for a tubular neighborhood of* $\gamma$ *of the sort that is depicted in* (2.1). *Assume that the closure of the image of* $\varphi$ *is disjoint from all other Reeb orbits from* $\Theta$. *Then* $\alpha_n$ *vanishes on any given constant* $t \in S^1$ *slice of* $\varphi(S^1 \times D)$ *with degree* m.

(4.6)

<u>Step 2</u>: Fix a smooth function $\chi: \mathbb{R} \to [0, 1]$ such that $\chi = 1$ on $(-\infty, \frac{5}{16}]$ and $\chi = 0$ on $[\frac{7}{16}, \infty)$. Given a pair $(A, \alpha)$ of connection on E and section of E, introduce the connection

$$\hat{A} = A - \tfrac{1}{2}(1 - \chi(|\alpha|^2))|\alpha|^{-2}(\bar\alpha \nabla_A \alpha - \alpha \nabla_A \bar\alpha).$$

(4.7)

Note that $\hat{A}$ is flat where $|\alpha|^2 > \tfrac{1}{2}$; and here $\alpha/|\alpha|$ is covariantly constant. Use $\hat{A}_n$ in what follows to denote the $(A_n, \alpha_n)$ version of $\hat{A}$. It consequence of the first item in (4.6) that $\hat{A}_n$ is flat at all points in M with distance $\varepsilon_n$ or greater from $\cup_{(\gamma,m)\in\Theta} \gamma$. Let $(\gamma, m) \in \Theta$. Let $\varphi: S^1 \times D \to M$ denote the coordinate chart map from the second item of (4.6). It is a consequence of this second item that

$$i \int_{S^1 \times D} dt \wedge \varphi^*(*B_{\hat{A}_n}) = m.$$

(4.8)

What with (2.1) and Lemma 2.2 in [T3], the equality in (4.8) implies that

$$i \int_M a_\tau \wedge *B_{\hat{A}_n} = 2\pi \sum_{(\gamma,m)\in\Theta} m\ell_\gamma + \mathfrak{e}$$

(4.9)

where $|\mathfrak{e}| \leq c_0 \varepsilon_n (1 + |E(A_n)|)$. Given that the sequence $\{\tau_n\}_{n=1,2,...}$ converges to $\tau$, a second appeal to Lemma 2.2 in [T3] now applied to (4.9) finds

$$i \int_M a_{\tau_n} \wedge *B_{\hat{A}_n} = 2\pi \sum_{(\gamma,m)\in\Theta} m\ell_\gamma + \mathfrak{e}_n$$

(4.10)

where $|\mathfrak{e}_n| \leq c_0 \delta_n (1 + K)$ with $\{\delta_n\} \in (0, 1)$ a sequence with limit zero as $n \to \infty$.

<u>Step 3</u>: Integrate by parts to see that

$$\left| \int_M a_{\tau_n} \wedge (*B_{\hat{A}_n} - *B_{A_n}) \right| = \left| \int_M da_{\tau_n} \wedge (1 - \chi(|\alpha_n|^2)) |\alpha_n|^{-2} (\bar{\alpha}_n \nabla_{A_n} \alpha_n - \alpha_n \nabla_{A_n} \bar{\alpha}_n) \right|$$
(4.11)

To bound the right hand side of (4.11), first define the function g on the domain $[0, \infty)$ by setting the rule

$$g(t) = -\int_t^2 (1 - \chi(s))s^{-1}ds$$
(4.12)

Since $(1 - \chi(|\alpha|^2))|\alpha|^{-2}(\bar{\alpha}\nabla_A \alpha + \alpha\nabla_A \bar{\alpha}) = d(g(|\alpha|^2))$, an integration by parts on the right hand side of (4.11) identifies the latter with

$$2\left| \int_M da_{\tau_n} \wedge (1 - \chi(|\alpha_n|^2)) |\alpha_n|^{-2} \alpha_n \nabla_{A_n} \bar{\alpha}_n \right| \ .$$
(4.13)

Next, use the Dirac equation to identify the covariant derivative of $\alpha_n$ along the Reeb vector field for $a_{\tau_n}$ with derivatives of $\beta_n$. Make this identification, and then use Holder's inequality to bound (4.13) by $\|\nabla_{A_n}\beta_n\|_2$, where the subscript indicates the $L^2$ norm.

To complete the bound on the right hand side of (4.11), integrate both sides of what is written in the second line of Equation (6.4) in [T3] over M and use the latter's Lemma 2.2 again to bound the $L^2$ norm of $\nabla_{A_n}\beta_n$ by $c_0 r_n^{-1/2}$.

<u>Step 4</u>: Given that the right hand side of (4.11) is bounded by $c_0 r_n^{-1/2}$, it thus follows from (4.10) that $E(A_n) \leq 2\pi \sum_{(\gamma,m) \in \Theta} m\ell_\gamma + c_0 \delta_n'$ where $\{\delta_n'\}_{n=1,2,...}$ is a sequence with limit zero as $n \to \infty$. However, this then implies that $E(A_n) < 2\pi L$ for all n sufficiently large. Given the assumptions, this is nonsense.

*Proof of Proposition 4.13*: Suppose for the moment that the constant $\kappa$ that appears in Proposition 4.11 can be chosen so that Proposition 4.11's conclusions hold with the same constant $\kappa$ for all $\tau \in [0, 1]$ versions of $(a_\tau, J_\tau)$. Suppose that there exists $\tau \in [0, 1]$ and a solution to (3.17) as defined with the data $(a_\tau, J_\tau)$, $\mu_\tau$ and r with $\mathfrak{a}^f = -\pi L r$. In the case when $c_1(\det(\mathbb{S}))$ is torsion, this implies that $\mathfrak{a} \geq -\pi L r - 2\pi^2 k$. In any case, an appeal to Proposition 4.11 finds that $E \leq 2\pi L + c_0 r^{-1/50}$. A further appeal to Proposition 4.12 finds that $E < 2\pi L - c_0^{-1}$ if $r \geq c_0$. Another look at Proposition 4.11 finds that $\mathfrak{a}^f \geq -\pi(L - c_0^{-1})r$.

The constant $\kappa$ that appears in Proposition 4.11 comes from the version that appears in Proposition 4.10. The value of $\kappa$ that makes the first item of Proposition 4.10 true comes from Proposition 5.1 in [T3] and Proposition 1.10 in [T4]. A look at the proofs of Proposition 5.1 in [T3] and Proposition 1.10 in [T4] finds that the latter version of $\kappa$ can be chosen so as to hold for every $\tau \in [0, 1]$ version of $(a_\tau, J_\tau)$ and $\mu_\tau$.

The value of $\kappa$ that makes the second item in Proposition 4.10 true comes from Lemma 2.4 in [T3]; and the value of $\kappa$ that makes the third item true comes via Proposition 1.9 in [T4]. A close look at the proofs of both of these propositions shows that their contributions to $\kappa$ can be assumed to be $\tau$ independent.

### h) Proof of Theorems 4.4 and 4.5

The preceding three subsections supply all of the heavy machinery for the proof. It remains only to put the machinery to use.

***Proof of Theorem 4.4***: Keep in mind that $E < 2\pi L$ on $\Phi^r(\mathcal{Z}_{ech}^L)$. If $c_1(\det(\mathbb{S}))$ is torsion, then it follows from the second item of Proposition 4.10 that $\mathfrak{a} \geq -\pi L r$ on $\Phi^r(\mathcal{Z}_{ech}^L)$ if $r \geq c_0$. If $c_1(\det(\mathbb{S}))$ is not torsion, then the third item of Proposition 4.10 finds $\mathfrak{a}^f > -\pi L r$ on $\Phi^r(\Lambda^L)$ if $r \geq c_0$. Now let $\mathfrak{d} \in \mathcal{M}_1(\mathfrak{c}_-, \mathfrak{c}_+)$. Since $s \to \mathfrak{a}(\mathfrak{d}(s))$ is a decreasing function on $\mathbb{R}$, it follows that $\mathfrak{a}(\mathfrak{c}_-) > -\pi L r$ if $c_1(\det(\mathbb{S}))$ is torsion, and $\mathfrak{a}^f(\mathfrak{c}_1) > -\pi L r$ otherwise. This understood, Proposition 4.11 implies that $E(\mathfrak{c}_-) < 2\pi L$ if $r \geq c_0$. According to Theorem 4.2, this implies that $\mathfrak{c}_- \in \Phi^r(\mathcal{Z}_{ech}^L)$.

***Proof of Theorem 4.5***: The proof has three parts.

*Part 1*: The lemma that follows summarizes the contribution from this part.

**Lemma 4.14**: *Let $(a, J)$ be a pair of contact 1-form from Lemma 2.1's residual set and almost complex structure from $\mathcal{J}_a$. Fix $k \in \mathbb{Z}$ when $c_1(\det(\mathbb{S}))$ is torsion. Then there exists $\kappa > 1$ and $\mathcal{K} \geq 1$ with the following significance: Fix $\mu_0 \in \Omega$ with $\mathcal{P}$ norm less than 1 as described in Proposition 4.1. Use $(a, J)$, $\mu_0$ and $r \in [r_*, \infty) - \mathcal{U}$ to define $\mathcal{M}^r$ and the Seiberg-Witten Floer cochain complex.*
- *Let $\theta$ denote a non-zero class in $H^*(\mathcal{C}^{SW})$, but with degree $k$ or greater if $c_1(\det(\mathbb{S}))$ is torsion. Then $\theta$ has a representative cocycle with $\mathfrak{a}^f$ greater than $-\pi \mathcal{K} r$ and with $E < 2\pi \mathcal{K}$ on each generator that appears with non-zero coefficient.*
- *Let $\mathcal{A} = -\pi \mathcal{K}$. There exists $r_{\mathcal{A}} \geq 1$ such that (4.2) holds.*
- *Let $\theta$ denote a non-zero class in $H^*(\mathcal{C}^{SW,\mathcal{A}})$, but with degree $k$ or greater if $c_1(\det(\mathbb{S}))$ is torsion. Then $\theta$ has a representative cocycle with $\mathfrak{a}^f > -\pi \kappa \mathcal{K} r$ and with $E < 2\pi \kappa \mathcal{K}$ on each generator that appears with non-zero coefficient.*
- *There exists $r_{\kappa \mathcal{A}}$ such that (4.2) using $r_{\kappa \mathcal{A}}$ and $\kappa \mathcal{A}$ in lieu of $r_{\mathcal{A}}$ and $\mathcal{A}$.*

***Proof of Lemma 4.14***: Let θ denote a non-zero class in H*($\mathcal{C}^{SW}$), but of degree k or greater when $c_1(\det(\mathbb{S}))$ is torsion. It follows from Proposition 4.10 that there exists $\mathcal{K} \geq 1$ and an unbounded set $\{r_i\}_{i=1,2,...} \in [r_*, \infty) - \mathcal{U}$ such that when $r = r_i$, then θ has a cocycle representative, $\mathfrak{n}$, with $E \leq 2\pi\mathcal{K}$ on each generator. Given that there the cohomology in question is finitely generated, the constant $\mathcal{K}$ can be taken so as to be independent of the choice for θ. Note that $\mathcal{K}$ can be chosen so that (4.5) holds with $L = \mathcal{K}$.

Fix $r_i$ and let $\mathfrak{n}$ denote the cocycle representative described above. It follows from Proposition 4.12 that $E < 2\pi\mathcal{K} - c_0^{-1}$ on each generator that appears $\mathfrak{n}$. The second and third items of Proposition 4.10 then find that $|\mathfrak{cs}^f| \leq c_0 r_i^{-1/50} (r_i E)$ on each generator that appears in $\mathfrak{n}$. This implies that $\mathfrak{a}^f_\theta > -2\pi(\mathcal{K}-c_0)r_i$. Given that $\mathfrak{a}^f_\theta[\cdot]$ is a continuous function, it follows from Proposition 4.13 that $\mathfrak{a}^f_\theta[r] > -\pi\mathcal{K}$ for $r > c_0^{-1}$. It then follows from Propositions 4.11 and 4.12 that each $r > c_0^{-1}$ version of θ has a cocycle representative with the property that $\mathfrak{a}^f > -\pi\mathcal{K}$ and $E < 2\pi\mathcal{K} - c_0^{-1}$ on each generator.

Proposition 4.13 implies what is asserted by the second bullet of the lemma. To obtain the third bullet's assertion, let θ denote a non-zero class in H*($\mathcal{C}^{SW,\mathcal{A}}$), but of degree k or greater when $c_1(\det(\mathbb{S}))$ is torsion. Repeat the arguments for the first bullet using this class θ to find $c_0 > 1$ such that $\mathfrak{a}^f > -\pi c_0 \mathcal{K}$ and $E < 2\pi c_0 \mathcal{K} - c_0^{-1}$. The final bullet again follows from Proposition 4.13.

*Part 2*: Let K denote the constant given in Lemma 4.14 and set $\mathcal{A} = -\pi\mathcal{K}$. It follows from Proposition 4.13 that (4.2) holds, and so (4.3) and (4.4) are well defined. Lemma 4.14 implies that the homomorphism $H^j(\mathcal{C}^{SW,\mathcal{A}}) \to H^j(\mathcal{C}^{SW})$ is surjective, at least if $j < k$ in the case when $c_1(\det(\mathbb{S}))$ is torsion.

Set $L = \mathcal{K}$ and fix some very small, but positive δ. Let $(\hat{a}, \hat{J})$ denote a $(\delta, L)$ approximation to $(a, J)$. Fix $\mu \in \Omega$ with $\mathcal{P}$ norm less than 1 as describe in the $(\hat{a}, \hat{J})$ version of Proposition 4.1. It follows from Proposition 4.13 that (4.2) is also obeyed by $(\hat{a}, \hat{J})$ and $\mu$. Let $\hat{\mathcal{C}}^{SW,\mathcal{A}}$ and $\hat{\mathcal{C}}^{SW}$ denote the $\mathbb{Z}$ modules that appear in the latter's version of (4.3). It is a consequence of Lemma 4.6 that $H^j(\hat{\mathcal{C}}^{SW,\mathcal{A}}) \approx H^j(\mathcal{C}^{SW,\mathcal{A}})$, at least for $j > k$ when $c_1(\det(\mathbb{S}))$ is torsion. Meanwhile, the $(\hat{a}, \hat{J})$ and $\mu$ version of Propositions 4.11 and 4.12 imply that $E < 2\pi L - c_0^{-1}$ on all generators of $\hat{\mathcal{C}}^{SW,\mathcal{A}}$ when $r > c_0$. Conversely, any generator of $\hat{\mathcal{C}}^{SW}$ with $E < 2\pi L + c_0^{-1}$ is a generator of $\hat{\mathcal{C}}^{SW,\mathcal{A}}$. This follows from the second and third bullets in Proposition 4.10.

Granted the preceding, the $(\hat{a}, \hat{J})$ and $\mu$ versions of Theorems 4.2-4.4 identify $H^j(\hat{\mathcal{C}}^{SW,\mathcal{A}}) \approx H_{-j}(\hat{\mathcal{C}}_{ech}^L)$, at least for $j > k$ when $c_1(\det(\mathbb{S}))$ is torsion. This is what is claimed in the first bullet of Theorem 4.5.

*Part 3*: To address the second bullet of Theorem 4.5, let $T_\Phi$ denote the monomorphism given by the $(\hat{a}, \hat{J})$ version of Theorem 4.3. Suppose that $\upsilon \in \hat{\mathcal{C}}_{ech}^{L_*}$ is such that $T_\Phi(\upsilon)$ is a coboundary. Let $\mathcal{A}_* = -\pi L_*$. Given that $T_\Phi$ identifies $H^j(\hat{\mathcal{C}}^{SW,\mathcal{A}_*})$ with $H_{-j}(\hat{\mathcal{C}}_{ech}^{L_*})$, it is sufficient to assume that the class $\hat{\lambda} \in H^*(\hat{\mathcal{C}}^{SW,\mathcal{A}_*})$ represented by $T_\Phi(\upsilon)$ is the image via (4.4)'s connecting homomorphism of a non-zero class in $H^{j-1}(\hat{\mathcal{C}}^{SW}/\hat{\mathcal{C}}^{SW,\mathcal{A}_*})$. It follows from Lemma 4.6 that this class corresponds to one in $H^{j-1}(\mathcal{C}^{SW}/\mathcal{C}^{SW,\mathcal{A}_*})$. Let $\theta$ denote the latter. If $r \geq c_0$, Propositions 4.12 and 4.13 find $\mathcal{K}_* \geq 1$ and a representative cocycle with $E < 2\pi\mathcal{K}_*$ and $\mathfrak{a}^f > -\pi\mathcal{K}_*$ on each generator. If $L$ is such that $L > \mathcal{K}_*$, then it follows that $\theta$ is represented by a cocycle in $\mathcal{C}^{SW,\mathcal{A}}$. Meanwhile, $\hat{\lambda}$ corresponds via the isomorphism from Lemma 4.6 to a class $\lambda \in H^j(\mathcal{C}^{SW,\mathcal{A}_*})$. The fact that $\theta$ is represented by a cocycle in $\mathcal{C}^{SW,\mathcal{A}}$ implies that $\lambda$ is zero in $H^j(\mathcal{C}^{SW,\mathcal{A}})$. This implies (again via the isomorphism from Lemma 4.6) that $\hat{\lambda}$ is zero in $H^j(\hat{\mathcal{C}}^{SW,\mathcal{A}})$. Given the identification of the latter with $H_{-j}(\hat{\mathcal{C}}_{ech}^L)$, this means that $\upsilon$ is a boundary in $\hat{\mathcal{C}}_{ech}^L$.

## 5. Theorems 4.2 and 4.3

The constructions that lead to $\Phi^r$ and $\Psi^r$ and the arguments for Theorems 4.2 and 4.3 are modifications of those used in [T1] and [T2] to prove the equivalence between the Gromov and Seiberg-Witten invariants of compact, symplectic 4-manifolds. What follows in this section is a brief description of what is involved. The details are contained in Papers II, III and IV of this series, [T6]-[T8].

### a) Vortices on $\mathbb{C}$

Both $\Phi^r$ and $\Psi^r$ use the solutions to the vortex equations on $\mathbb{C}$ to construct pairs of connection on $E$ and section of $\mathbb{S}$ from Reeb orbit or pseudoholomorphic curves as the case may be. Solutions to the vortex equations played a similar role in the contructions done in [T2]; see its article *Gr => SW*. What follows in this subsection provides a brief summary of the vortex part of the story.

The vortex moduli spaces are labeled by a non-negative integer, with the integer n version of the vortex moduli space denoted by $\mathfrak{C}_n$. The latter consists of certain equivalence classes of pairs $(A, \alpha)$, where $A$ is a hermitian connection on the trivial complex line bundle over $\mathbb{C}$, and where $\alpha$ is a section of this bundle. A pair $\mathfrak{c} = (A, \alpha)$ is in $\mathfrak{C}_n$ if and only if the curvature of $A$ and the $A$-covariant derivative of $\alpha$ satisfy

- $*F_A = -i(1 - |\alpha|^2)$.
- $\bar{\partial}_A \alpha = 0$.
- $|\alpha| \leq 1$.
- *The function* $(1 - |\alpha|^2)$ *is integrable on* $\mathbb{C}$ *and* $\int_\mathbb{C} (1 - |\alpha|^2) = 2\pi n$

(5.1)

Here, $\bar{\partial}_A$ denotes the d-bar operator that is defined by the connection A. The equivalence relation that defines a point in $\mathfrak{C}_n$ identifies pairs (A, α) and (A´, α´) when $A´ = A - u^{-1}du$ and α´ = uα where u is a smooth map from $\mathbb{C}$ to $S^1$.

The space $\mathfrak{C}_0$ consists of a single element, this the gauge equivalence class of the pair (A = 0, α = 1). When n ≥ 1, the space $\mathfrak{C}_n$ has the structure of a smooth, complex manifold that is biholomorphic to $\mathbb{C}^n$. This holomorphic identification is realized as follows: As is explained by the author in Section 2 of the article *Gr => SW* from [T2], [T10] and [JT], if (A, α) solves (5.1), then α has precisely n zero counting multiplicities. Let $\mathfrak{z}_{\mathfrak{c}} = \{z_1, \ldots, z_n\}$ denote the resulting set in the n'th symmetric product of $\mathbb{C}^n$. A holomorphic diffeomorphism from $\mathfrak{C}_n$ to $\mathbb{C}^n$ sends $\mathfrak{c}$ to the point in $\mathbb{C}^n$ whose q'th coordinate is $\sum_{1 \le j \le n} z_j^q$. As it turns out

$$\sum_{1 \le j \le n} z_j^q = \tfrac{1}{2\pi} \int_{\mathbb{C}} z^q (1 - |\alpha|^2).$$

(5.2)

With regards to the integral in (5.2), note that

$$0 < 1 - |\alpha|^2 < c_0 \sum_{1 \le j \le m} e^{-\sqrt{2}|z - z_j|} \quad and \quad |\nabla_A \alpha|^2 \le c_0 \sum_{1 \le j \le m} e^{-\sqrt{2}|z - z_j|}$$

(5.3)

where $c_0$ is a constant that is independent of n and (A, α) ∈ $\mathfrak{C}_m$.

This holomorphic identification to $\mathbb{C}^n$ does <u>not</u> provide the natural Riemannian metric on $\mathfrak{C}_n$ when n > 1. The relevant metric is described momentarily. To set the stage for the story on the metric, remark that the (1, 0) tangent space at a given $\mathfrak{c}$ = (A, α) to $\mathfrak{C}_n$ is naturally isomorphic to a certain vector space of pairs (x, ι) where x is a $\mathbb{C}$-valued function on $\mathbb{C}$ and ι is a section of the trivial bundle. To lie in $T_{1,0}\mathfrak{C}_n|_{\mathfrak{c}}$, both x and ι must be square integrable on $\mathbb{C}$ and obey the coupled system

$$\partial x + 2^{-1/2} \bar{\alpha} \iota = 0 \quad and \quad \bar{\partial}_A \iota + 2^{-1/2} \alpha x = 0.$$

(5.4)

In this equation, $\partial$ is shorthand for $\tfrac{\partial}{\partial z}$. The pair (x, ι) provides the first order change in (A, α) that adds $2^{-1/2}x$ to the (0, 1) part of A and adds ι to α. The relevant metric on $\mathfrak{C}_n$ is defined so that the metric norm of (x, ι) is $\pi^{-1/2}$ times its $L^2$ norm as defined by integration on $\mathbb{C}$. This metric is a Kahler metric with respect to the complex structure.

**b) Theorem 4.2**

This subsection is meant to give an indication of what is involved in constructing the map $\Phi^r$. To set the stage, write the spinor bundle $\mathbb{S}$ as $E \oplus EK^{-1}$. With L > 1 fixed, introduce as notation $\mathcal{Z}^L$ to denote the set whose typical element, Θ, consists of pairs of the form (γ, m) where γ is a Reeb orbit and m a positive integer. Require that distinct pairs from Θ have distinct Reeb orbit components. In addition, require $\sum_{(\gamma,m)\in\Theta} m\ell_\gamma < L$; and require that the formal sum $\sum_{(\gamma,m)\in\Theta} m\gamma$ define a cycle whose class in $H_1(M; \mathbb{Z})$ is Poincaré dual to $c_1(E)$.

The subsequent description of $\Phi^r$ has four parts.

*Part 1*: Fix $\Theta \in \mathcal{Z}^L$. Assign to each $(\gamma, m) \in \Theta \in \mathcal{Z}^L$ a smooth map, $\mathfrak{c}_\gamma: S^1 \to \mathfrak{C}_m$. The first step to defining $\Phi^r$ is to construct a pair in $\mathrm{Conn}(E) \times C^\infty(M; \mathbb{S})$ from the data $\{\mathfrak{c}_\gamma\}_{(\gamma,m)\in\Theta}$. The map $\mathfrak{c}_\gamma: S^1 \to \mathfrak{C}_m$ is lifted so as to give a connection on, and section of the product complex line bundle over $S^1 \times \mathbb{C}$. The pull-back of this pair to any given constant $t \in S^1$ slice of $S^1 \times \mathbb{C}$ is a solution to (5.1). Define $\hat{r}_\gamma: S^1 \times \mathbb{C} \to S^1 \times \mathbb{C}$ so that $\hat{r}_\gamma^*(t, z) = (t, (\frac{\ell_\gamma}{2\pi} r)^{1/2} z)$. Use $(A^{(\gamma)}, \alpha^{(\gamma)})$ to denote the pull-back via $\hat{r}_\gamma$ of the chosen lift of $\mathfrak{c}_\gamma$. As might be expected from (5.3), the connection $A^{(\gamma)}$ is nearly flat where $|z| \gg r^{-1/2}$ on $S^1 \times \mathbb{C}$. Meanwhile, $\alpha^{(\gamma)}$ is nearly $A^{(\gamma)}$-covariantly constant with norm 1 on this same part of $S^1 \times \mathbb{C}$.

Suppose now that a tubular neighborhood map has been chosen for each Reeb orbit from Θ of the sort that is described in Section 2a. Use such a map to identify a tubular neighborhood of each Reeb orbit with $S^1 \times D \subset S^1 \times \mathbb{C}$. Then $(A^{(\gamma)}, \alpha^{(\gamma)})$ can be viewed as a pair of connection on, and section of the product complex line bundle over a tubular neighborhood of γ in M with $A^{(\gamma)}$ near flat near the boundary of this tubular neighborhood, and with $\alpha^{(\gamma)}$ having norm 1 and nearly covariantly constant near this same boundary.

This understood, the collection $\{(A^{(\gamma)}, \alpha^{(\gamma)})\}_{(\gamma,m)\in\Theta}$ are pasted together using 'bump' functions so as to obtain a pair $(A^*, \alpha^*)$ of connection on and section of a complex line bundle over M. The latter is isomorphic to E. Note that $A^*$ is flat except very near the Reeb orbits from Θ; and likewise $\alpha^*$ is covariantly constant with norm 1 except very near these same Reeb orbits.

Write $\mathbb{S} = E \oplus EK^{-1}$ and define $(A^*, \psi^* = (\alpha^*, 0)) \in \mathrm{Conn}(E) \times C^\infty(M; \mathbb{S})$. A calculation shows that $(A^*, \psi^*)$ comes reasonably close to solving (3.17) if r is large. Note that such is the case by virtue of the fact that the pairs in Θ involve <u>Reeb</u> <u>orbits</u>. Indeed, the construction just described can be applied to a set such as Θ whose typical element is a pair (γ, m) where γ is an embedded loop in M. If γ is not a Reeb orbit, then the resulting $(A^*, \psi^*)$ will not come close to solving (3.17) when r is large.

The plan is to look for a solution to (3.17) near $(A^*, \psi^*)$. Such a solution can be found when r is large if the collection $\{\mathfrak{c}_\gamma: S^1 \to \mathfrak{C}_m\}_{(\gamma, m)\in\Theta}$ are suitably constrained.

*Part 2*: What follows describes the constraints on $\{\mathfrak{c}_\gamma\}_{(\gamma,m) \in \Theta}$. To this end, return to the vortex moduli space $\mathfrak{C}_m$. Let $(\nu, \mu)$ denote a pair consisting of a real number and a complex number. Any such pair defines a function, $\mathfrak{h}$, on $\mathfrak{C}_m$ given by

$$\mathfrak{h} = \tfrac{1}{4\pi} \int_\mathbb{C} (2\nu |z|^2 + (\mu \bar{z}^2 + \bar{\mu} z^2))(1 - |\alpha|^2) \ .$$

(5.5)

As with any function on $\mathfrak{C}_m$, this one defines a Hamiltonian vector field. Now suppose that $\nu$ and $\mu$ are respectively, a real valued function on $S^1$ and a $\mathbb{C}$-valued function on $S^1$. Then (5.5) defines a 1-parameter family of Hamiltonian vector fields on $\mathfrak{C}_m$. Of interest are the closed, integral curves of the latter. These are maps $\mathfrak{c} \colon S^1 \to \mathfrak{C}_m$ that obey at each t $\in S^1$ the equation

$$\tfrac{i}{2} \mathfrak{c}' + \nabla^{(1,0)} \mathfrak{h}|_\mathfrak{c} = 0 \ ,$$

(5.6)

where $\mathfrak{c}'$ is shorthand for the (1, 0) part of $\mathfrak{c}_* \tfrac{d}{dt}$, and where $\nabla^{(1,0)} \mathfrak{h}$ denotes the (1, 0) part of the gradient of $\mathfrak{h}$.

Now suppose that $\gamma$ is a Reeb orbit. Fix a tubular neighborhood map for $\gamma$ of the sort described in Section 2a. Then $\gamma$ has an associated pair $(\nu, \mu)$ for use in (5.5), this the pair that appears in (2.1) and (2.3). With the preceding understood, what follows is the key observation. Suppose that the following is true:

*Each $(\gamma, m) \in \Theta$ version of $\mathfrak{c}_\gamma$ is a solution to the corresponding version of (5.6). In addition, the linearized version of the left hand side of (5.6) at this $\mathfrak{c}_\gamma$ defines an operator with trivial kernel.*

(5.7)

Note that the linearization of (5.6) at a given map $\mathfrak{c} \colon S^1 \to \mathfrak{C}_m$ defines a first order, elliptic and symmetric operator on $C^\infty(S^1, \mathfrak{c}^* T^{1,0} \mathfrak{C}_m)$.

Under the assumption in (5.7), perturbation theory can be employed to modify $(A^*, \psi^*)$ when r is large so that the result solves the corresponding version of (3.17). To put this in a more formal way, introduce $\mathfrak{C}\Theta$ for $\Theta \in \mathcal{Z}^L$ to denote the set whose typical element assigns to each $(\gamma, m) \in \Theta$ a corresponding solution to (5.6). Say that a solution to (5.6) is *nondegenerate* when the linearization of the left hand side of (5.6) at the solution has trivial kernel. Use $\mathfrak{C}\Theta^*$ to denote the subset where each assigned solution is nondegenerate. Given $L \geq 1$, let $\mathfrak{C}\mathcal{Z}^L$ denote $\{\mathfrak{C}\Theta \colon \sum_{(\gamma,m) \in \Theta} m \ell_\gamma < L\}$. Use $\mathfrak{C}\mathcal{Z}^{L*}$ to denote the subset $\{\mathfrak{C}\Theta^* \colon \Theta \in \mathcal{Z}^L\} \subset \mathfrak{C}\mathcal{Z}^L$.

*If the contact form comes from Lemma 2.1's residual subset, if $\mathfrak{C}\mathcal{Z}^{L*} = \mathfrak{C}\mathcal{Z}^L$, and if r is sufficiently large, then perturbation theory defines an injective map*

$$\Phi^r: \mathfrak{CZ}^L \to \mathcal{M}^r \text{ whose image consists of the set of elements with } \mathrm{E} < 2\pi\mathrm{L}.$$
(5.8)

This map is constructed in [T6]; and [T8] proves that it is injective and surjective onto the $\mathrm{E} < 2\pi\mathrm{L}$ subset in $\mathcal{M}^r$. It is fair to say that these parts of [T6] and [T8] do little more than reinterpret parts of the respective $Gr \Rightarrow SW$ and $Gr = SW$ articles in [T2].

If $\mathfrak{CZ}^L \neq \mathfrak{CZ}^{L*}$, then perturbation theory constructs, for each sufficiently large r, an injective map from a certain subset of $\mathfrak{CZ}^L$ into $\mathcal{M}^r$ whose image consists of the set of elements with $\mathrm{E} < 2\pi\mathrm{L}$. The latter is also denoted by $\Phi^r$ in what follows.

*Part 3*: What follows says some things about the space of solutions to (5.6). To start, note that the solution space to (5.6) is compact if $(\nu, \mu)$ is either hyperbolic or m-elliptic. This is proved in [T8]. In either case, there exists a unique solution for m = 1; this is the vortex with $\alpha^{-1}(0) = \{0\}$.. The unique m = 1 solution to (5.6) is nondegenerate if $(\nu, \mu)$ is nondegenerate. In this m = 1 case, the corresponding linear operator is the operator L that is depicted in (2.3).

It is not known whether the solution space to the m > 1 versions of (5.6) consists of solely nondegenerate solutions, even if $(\nu, \mu)$ are chosen in a generic fashion. However, if each Reeb orbit with $\ell_\gamma \leq \mathrm{L}$ is either hyperbolic or m-elliptic, and if each of the corresponding versions of (5.6) has solely nondegenerate solutions when $\ell_\gamma \leq \mathrm{L}$, then $\mathfrak{CZ}^L$ is a finite set.

All solutions to (5.6) are known for some specific choices of $\nu$ and $\mu$:

- *If $(\nu, \mu) = (\frac{1}{4}k, i\varepsilon e^{ikt})$ is hyperbolic with rotation number k. Here, $\varepsilon > 0$ but very small. Then there are no m > 1 solutions to (5.6).*
- *If $(\nu, \mu) = (\frac{1}{2}\mathrm{R}, 0)$ with R irrational, then there is a unique solution to (5.6) for each m, this the vortex with $\alpha^{-1}(0) = \{0\}$. Moreover, the latter is nondegenerate.*

(5.9)

These last facts are proved in [T6]. If (5.9) holds for each Reeb orbit $\gamma$ with $\ell_\gamma < \mathrm{L}$, then the set $\mathfrak{CZ}^L$ is precisely the set $\mathcal{Z}_{ech}{}^L$ that gives the generators of the embedded contact homology subcomplex $\mathcal{C}_{ech}{}^L$. In this case, $\Phi^r$ is a map from $\mathcal{Z}_{ech}{}^L$ into $\mathcal{M}^r$; this is the map used in Theorem 4.2.

*Part 4*: As it turns out, the proof that embedded contact homology is isomorphic to Seiberg-Witten Floer cohomology does not require knowledge of all solutions to (5.6). Knowledge of the corresponding Hamiltonian Floer cohomology groups is sufficient.

To elaborate, Floer [F1] [F2] introduced his celebrated 'Floer (co)homology' to resolve a famous conjecture of Arnold that concerned closed orbits of time-dependent Hamiltonian vector fields on symplectic manifolds. What is written in (5.6) is an example of just such a Hamiltonian dynamical system. In particular, if $(\nu, \mu)$ is a nondegenerate pair, then there are well defined, $\mathbb{Z}$-graded Floer homology and cohomology groups whose generators are solutions to a suitably generic compactly supported (on $\mathfrak{C}_m$) perturbation of (5.6).

There is one subtle point here, this involving the instantons that define the differentials. In this context, an instanton is a smooth map from $\mathfrak{c}: \mathbb{R} \times S^1 \to \mathfrak{C}_m$ that obeys the equation

$$\bar{\partial}\mathfrak{c} + \nabla^{1,0}h|_\mathfrak{c} = 0 ,$$

(5.10)

where $\bar{\partial}$ is a suitably defined version of the d-bar operator on $\mathfrak{c}^*T^{1,0}\mathfrak{C}_m \to S^1 \times \mathbb{R}$. The map $\mathfrak{c}$ must also limit as $s \to \pm\infty$ to a solution of (5.6). In order to obtain a well defined differential for the Hamiltonian Floer (co)homology, it is necessary to prove that the moduli space of instanton solutions to (5.10) can be compactified by adding 'broken trajectories'. This can be done when $(\nu, \mu)$ are either hyperbolic or m-elliptic.

In any event, it can be shown that there are well defined Hamiltonian Floer (co)homology groups for (5.6) when $(\nu, \mu)$ are either hyperbolic or m-elliptic. Furthermore, it can be shown (using (5.9)) that the Hamiltonian Floer cohomology groups are as follows:

- $\mathbb{Z}$ *if* $m = 1$.
- $0$ *if* $m > 1$ *and* $(\nu, \mu)$ *is hyperbolic*.
- $\mathbb{Z}$ *if* $m > 1$ *and* $(\nu, \mu)$ *is m-elliptic*.

(5.11)

Suppose that each $\Theta \in \mathcal{Z}^L$ contains only pairs $(\gamma, m)$ such that $\gamma$ is either hyperbolic or m-elliptic. Granted only this assumption, it is none-the-less the case that the points in $\mathfrak{C}\Theta$ that are mapped by $\Phi^r$ into $\mathcal{M}^r$ carry, in a suitable sense, the product of the Hamiltonian Floer cohomology groups as defined by the various $(\gamma, m) \in \Theta$ versions of (5.6). Somewhat more is said about this in the next subsection.

**c) Theorem 4.3**

This section is meant to give a rough indication of how $\Psi^r$ is constructed. There are three parts to what follows. The final two parts say something about what is involved when (5.9) does not hold.

*Part 1*: The idea is to mimick as much as possible what is done in the *Gr => SW* article in [T2]. As done there, the first step constructs an approximate solution to (3.18) such that the connection is flat except very near to the given curve in $\mathcal{M}_1(\Theta_-, \Theta_+)$, the section of $\mathbb{S}$ lies only in the E summand, and this section of E is covariantly constant with norm 1 except very near the given curve. Step 2 uses perturbation theoretic techniques to find an honest solution to (3.18) that differs little from the approximate one. There are, however, serious new issues that do not arise in [T2], these relating to the behavior of the elements of $\mathcal{M}_1(\Theta_-, \Theta_+)$ where |s| is large on $\mathbb{R} \times M$.

To elaborate on this last point, suppose that $\Sigma \subset \mathbb{R} \times M$ is an embedded, pseudoholomorphic curve. As such, $\Sigma$ has a well defined normal bundle, $N \to \Sigma$, and a fixed radius disk bundle $N_1 \subset N$ with an exponential map $e_\Sigma \colon N_1 \to \mathbb{R} \times M$ that immerses $N_1$ and embeds a neighborhood of the zero section as a neighborhood of $\Sigma$. Even so, there need not exist a fixed radius disk subbundle of N that is everywhere embedded by $e_\Sigma$. The point being that the constant s slices of distinct ends of $\Sigma$ can limit as $s \to \pm$ to the same Reeb orbit. In addition, the constant s slices of any given end need not define a degree 1 braid in the tubular neighborhood of the nearby Reeb orbit.

These remarks about the fixed radius disk bundle are relevant because the constructions in the article *Gr =>SW* from [T2] require an embedding of just such a bundle. When a fixed radius disk bundle is embedded by $e_C$, then the constructions in [T2] can be copied with only minor changes to produce an approximate solution to (3.18) and then a deformation of the latter to an honest solution. In this regard, the approximate solution has the following appearance: Use the exponential map to identify the fixed radius subbundle of $N_1$ with a tubular neighborhood of $\Sigma$ in $\mathbb{R} \times M$. The pull-back of the connection and the section of E to any given fiber of $N_1$ differs little from the pullback of a solution to (5.1) via the map from $\mathbb{C}$ to $\mathbb{C}$ that sends any given $z \in \mathbb{C}$ to $r^{1/2}z$.

In general, only the following can be guaranteed: Given $R > 1$, there exists $\rho > 0$ such that $e_\Sigma$ embeds the radius $\rho$ disk bundle in $N_1$ where |s| < R. This understood, the constructions in the article *Gr =>SW* from [T2] need modifications at large |s| on $\mathbb{R} \times M$. The full details are given in [T6] and [T8]; they account for the length of these papers. Said briefly, the approximate solution at moderate values of |s| is constructed as in the article *Gr =>SW* from [T2]. At points where s << -1, the curve $\Sigma$ is very near the s << -1 part of a union of $\mathbb{R}$-invariant cylinders, each of the form $\mathbb{R} \times \gamma$ with $\gamma \subset M$ a Reeb orbit. The constructions in *Gr =>SW* from [T2] are applied using these cylinders in lieu of $\Sigma$ to obtain an approximate solution where s << -1. Likewise, the constructions in *Gr =>SW* from [T2] are applied to the s >> 1 part of another union of $\mathbb{R}$ invariant cylinders to obtain the approximate solutions on this same part of $\mathbb{R} \times M$. The approximate solutions on these three regions in $\mathbb{R} \times M$ are then glued together where the regions overlap so as to obtain an approximate solution on the whole of $\mathbb{R} \times M$. It is a consequence of (5.3) that a gluing of this sort will result in a pair $(A^*, \psi^*)$ that nearly solves (3.18) when r is large.

With the approximate solution in hand, a perturbative construction finds a nearby (A, ψ) that obeys (3.18) on the nose. The latter construction is somewhat more complicated than that in *Gr =>SW* from [T2].

*Part 2*: The assumption in (5.9) greatly simplifies matters. The analog of Theorem 4.3 when (5.9) is not assumed is very much more complicated. The complications are two fold: First, $\Phi^r$ now associates to each $\Theta \in \mathcal{Z}^L$ a set, $\Phi_\Theta$, of elements in $\mathcal{M}^r$. These elements do not all have the same degree and there will, in general, be instanton solutions to (3.19) with both $s \to -\infty$ limit and $s \to +\infty$ limit in $\Phi_\Theta$. With degrees and signs taken into account, these sorts of instantons compute the product of the Hamiltonian Floer cohomology given in (5.11) for the various pairs $(\gamma, m) \in \Theta$.

Meanwhile, if $\Theta_-$ and $\Theta_+$ are distinct elements in $\mathcal{Z}^L$, there may be instanton solutions with $s \to -\infty$ and $s \to +\infty$ limits in respectively $\Phi_{\Theta_-}$ and $\Phi_{\Theta_+}$. However, the set of such solutions is not necessarily in 1-1 correspondence with $\mathcal{M}_1(\Theta_-, \Theta_+)$. If (5.9) does not hold, then each instanton in $\mathcal{M}_1(\Theta_-, \Theta_+)$ can determine a number of instantons solutions to (3.19), even when both $\Phi_{\Theta_-}$ and $\Phi_{\Theta_+}$ consist of a single element. To say more about this last point, recall that an approximate solution to (3.19) for, say $s \ll -1$, is constructed using as template what is done in *Gr =>SW* from [T2] with the pseudoholomorphic curve taken to be a product of cylinders. The template from the article *Gr => SW* in [T2] requires a solution to (5.10) for each such cylinder. In this regard, the $s \to -\infty$ limit of the solution for a cylinder $\mathbb{R} \times \gamma$ with $(\gamma, m) \in \Theta_-$ must be a solution $c_\gamma$ to (5.6). However, the $s \to \infty$ limit must be quite different since it has to match up with what is done at moderate values of $|s|$ using the template from *Gr =>SW* from [T2] as applied to the given $\Sigma \subset \mathcal{M}_1(\Theta_-, \Theta_+)$. The precise behavior of $\alpha^{-1}(0)$ at large s is determined by the various versions of (2.6) that come from the ends of $\Sigma$ whose constant s slices converge as $s \to -\infty$ to $\gamma$. In particular, the solution to (5.10) must be such that $\sup_{t \in S^1} \{\text{dist}(\alpha^{-1}(0), 0)|_{(s,t)}\}$ diverges as $s \to \infty$. The story is even more complicated if there are two or more ends involved and they have distinct versions of what is denoted as $q_E$ in Section 2b.

In the case when (5.9) holds, there is but a single relevant solutions to (5.10); and the story, though still long, is more or less straightforward. If (5.9) does not hold, then there may be many relevant solutions to (5.9), and then each will determine a distinct instanton solution to (3.19).

*Part 3*: The upshot of all of this is that when (5.9) does not hold, the proof that the embedded contact homology is isomorphic to Seiberg-Witten Floer cohomology requires much more work, both on the analytic side, and on the algebraic side. What follows is meant to give a rough indication of what is involved on the algebraic side:

Each element in $\mathcal{Z}^L$ determines some number of generators in $\mathcal{C}^{SW}$. This is the case even for elements that are not in $\mathcal{Z}_{ech}^L$ and so are not considered generators of $\mathcal{C}_{ech}^L$. Note that an element in $\mathcal{Z}^L - \mathcal{Z}_{ech}^L$ pairs one or more hyperbolic Reeb orbits with an integer greater than 1. In any event, with the differentials taken into account, each $\Theta \in \mathcal{Z}^L$ determines a submodule, $\mathbb{Z}\Phi_\Theta \subset \mathcal{C}^{SW}$. Let $L_1 < L_2, < \cdots$ denote the ordered set of numbers that can be obtained as $\sum_{(\gamma,m) \in \Theta} m \ell_\gamma$ with $\Theta \in \mathcal{Z}^L$. If r is sufficiently large, then the $\mathfrak{a}^f \geq -\pi L r$ subcomplex of the Seiberg-Witten Floer cochain complex can be filtered as

$$\cdots \subset \oplus_{\Theta \in \mathcal{Z}^{L_k}} \mathbb{Z}\Phi_\Theta \subset \oplus_{\Theta \in \mathcal{Z}^{L_{k+1}}} \mathbb{Z}\Phi_\Theta \subset \cdots .$$

(5.12)

The $E_2$ term of the corresponding spectral sequence is isomorphic to the free $\mathbb{Z}$ module generated by the generators of $\mathcal{C}_{ech}^L$. This follows from the aforementioned fact that the cohomology of $\mathbb{Z}\Phi_\Theta$ is isomorphic to the product of the various $(\gamma, m)$ versions of (5.11). In particular, this cohomology is isomorphic to either $\mathbb{Z}$ or 0 with $\mathbb{Z}$ arising if and only if $\Theta$ is in $\mathcal{Z}_{ech}^L$ and so gives a generator of $\mathcal{C}_{ech}^L$. The induced differential on the $E_2$ term of the spectral sequence corresponding to (5.12) should be identical to the differential on $\mathcal{C}_{ech}^L$.

**Appendix: Proof of Proposition 2.5**

Let $\mathcal{R}_L$ denote the set of Reeb orbits for a with symplectic action less than L. To set the stage for the constructions that follow, agree to associate a tubular neighborhood map from $S^1 \times D$ as described in Sections 2a and 2b for each Reeb orbit in $\mathcal{R}_L$. Since there are but a finite number of such Reeb orbits, no generality is lost by assuming that these tubular neighborhood maps have pairwise disjoint image. When $\gamma \in \mathcal{R}_L$, the associated tubular neighborhood map is denoted by $\varphi_\gamma$. This map is used to identify its image and domain so as to view the 1-form a near $\gamma$ as the 1-form on $S^1 \times D$ that is given by $\ell_\gamma$ times what is written on the right hand side of (2.1).

Use $\varphi_\gamma$ to define the functions $(\nu, \mu)$ that appear in (2.1). This done, fix a homotopy $\{\tau \to (\nu_\tau, \mu_\tau)\}_{\tau \in [0,1]}$ as described in Lemma 2.3 with $(\nu_0, \mu_0) = (\nu, \mu)$. Choose this homotopy to be independent of $\tau$ near 0 and near 1.

What follows describes how the proof proceeds: Let Q denote a very large integer. The plan is to define a sequence $\{(a_k, J_k)\}_{k=0,1,\ldots,Q}$ such that $(a_0, J_0) = (a, J)$ and such that for each $k \in \{0, \ldots, Q\}$, all but the fourth item in (2.11) are satisfied if $\hat{a} = a_k$ and $\hat{J} = J_k$. For $k \geq 1$, the latter is replaced by:

- $\frac{2\pi}{\ell_\gamma} \varphi^* a_k = (1 - 2\nu_{\tau=k/N}|z|^2 - \mu_{\tau=k/N}\bar{z}^2 - \bar{\mu}_{\tau=k/Q}z^2) dt + \frac{i}{2}(z d\bar{z} - \bar{z} dz)$

- $\varphi^*T^{1,0}(\mathbb{R}\times M)$ *is spanned by* $ds + i\varphi^*a_k$ *and* $\frac{\ell_\gamma}{2\pi}(dz - 2i(\nu_{\tau=k/Q}z + \mu_{\tau=k/Q}\bar{z})dt)$

(A.1)

An inductive argument is used to make these constructions. A lower bound for the integer Q is described below.

Two facts play a prime roll in the construction of the sequence $\{(a_k, J_k)\}_{k=1,2,...}$. Here is the first:

*The assignment of $(\tau, t) \in [0, 1] \times S^1$ to any given Reeb orbit's version of the pair $(\nu_\tau(t), \mu_\tau(t))$ defines a smooth map to $\mathbb{R} \times \mathbb{C}$. The derivatives of this map to any given fixed order enjoy a uniform bound that is independent of $\gamma \in \mathcal{R}_L$.*

(A.2)

To state the second key fact, introduce $L_0$ to denote the smallest of the lengths of all closed Reeb orbits.

*For any given positive integer $q \leq L_0^{-1}L + 1$, there is a positive lower bound, independent of $\gamma \in \mathcal{R}_L$ and $\tau \in [0, 1]$, to the absolute value of any eigenvalue of the corresponding $(\nu_\tau, \mu_\tau)$ version of (2.3) on the space of $2\pi q$ periodic functions.*

(A.3)

The statement of the third fact requires a digression to set the notation. To start it, fix generators $\Theta_-$ and $\Theta_+$ from $\mathcal{C}_*^L$. An element from J's version of $\mathcal{M}_1(\Theta_-, \Theta_+)$ consists of some number of $\mathbb{R}$-invariant cylinders with integer weights, and one non $\mathbb{R}$-invariant, irreducible submanifold. Let C denote either one of these cylinders, or the non $\mathbb{R}$-invariant submanifold. Deformations of C that preserve to first order the J invariance of its tangent space can be viewed with the help of a suitable exponential map as sections of the normal bundle $N \to C$ that obey a certain first order, $\mathbb{R}$-linear elliptic equation. The linear operator that defines this equation is denoted by $\mathcal{D}_C$. This is the operator in (2.3) that was briefly described in the paragraph that follows Lemma 2.2. As noted there, it defines a bounded, Fredholm map from $L^2_1(C; N)$ to $L^2(C; N \otimes T^{0,1}C)$. In this guise, its cokernel is trivial, and so it has an inverse that gives a bounded, linear map to the $L^2$ orthogonal complement of its kernel. This kernel is trivial if $C = \mathbb{R} \times \gamma$. The inverse of $\mathcal{D}_C$ is denoted by $\mathcal{D}_C^{-1}$.

With this notation set, here is the third point:

*There is a bound, independent of $(\Theta_-, \Theta_+)$ and C from $\mathcal{M}_1(\Theta_-, \Theta_+)$, to the norm of $\mathcal{D}_C^{-1}$.*

(A.4)

To explain, remember that there are but a finite number of pairs $(\Theta_-, \Theta_+)$ to choose from, and for each, there are but a finite number of components in $\mathcal{M}_1(\Theta_-, \Theta_+)$. Thus, up to the action of $\mathbb{R}$, there are but a finite number of possible choices for C. The lower bound for the integer N is determined by the bounds in (A.2), (A.3) and (A.4). To this end, let $\lambda_0$ denote the bound that is alluded to in (A.3) and let $\sigma_0$ denote the bound that is alluded to in (A.4).

To initiate the induction, note that $(a_0, J_0)$ are such that $(\hat{a} = a_0, \hat{J} = J_0)$ satisfy (A.1) and all but the fourth item of (2.11). Suppose that $k \in \{0, 2, \ldots, Q-1\}$ and that a pair $(a_k, J_k)$ have been defined so that $(\hat{a} = a_k, \hat{J} = J_k)$ satisfies all but the fourth item of (2.11) and (A.1) if $k \geq 1$. The assertion in (A.4) also holds when $\mathcal{M}_1(\cdot, \cdot)$ is defined using $J_k$. This is because the set of subvarieties in question is finite.

The induction from $k$ to $k+1$ requires one additional and crucial input, this a constant, $\sigma_* \geq 1$, whose definition follows. To start, let $\gamma \in \mathcal{R}_L$ and let $C = \mathbb{R} \times \gamma$. For each $\tau \in [0, 1]$, use $L_\tau$ to denote the version of (2.2) that has $\gamma$'s version of $(\nu, \mu)$ replaced by $(\nu_\tau, \mu_\tau)$. For each integer $q \in \{1, \ldots, L_0/L + 1\}$, let $\mathcal{D}_{C,\tau,q}$ denote the operator $\frac{\partial}{\partial s} + L_\tau$ with domain the space of complex valued, $L^2_1$ functions on $\mathbb{R} \times (\mathbb{R}/2\pi q\mathbb{Z})$ and range the space of complex valued, $L^2$ functions on $\mathbb{R} \times (\mathbb{R}/2\pi q\mathbb{Z})$. By virtue of (A.3), this operator is invertible. Take $\sigma_*$ to be a $q$, $\tau$ and $\gamma \in \mathcal{R}_L$ independent upper bound for the norm of this inverse. This constant $\sigma_*$ is determined by $\lambda_0$ and a sup norm bound for all $\gamma \in \mathcal{R}_L$ versions of the pair $(\nu, \mu)$.

The completion of the induction from $k$ to $k+1$ is presented below in eight parts. Before starting, take note of the convention used here that '$c_0$' denotes a constant that is independent of the relevant variables. It's value is greater than 1 and it can be assumed to increase between subsequent appearances.

*Part 1*: To construct a candidate for $a_{k+1}$, remark first that there exists, by assumption, some $\rho_k \ll \delta$ with the following significance: Let $\gamma \in \mathcal{R}_L$; and use $\varphi_\gamma$ again to identify a tubular neighborhood of $\gamma$ with $S^1 \times D$. The 1-form $a_k$ on the $|z| < \rho_k$ part of $S^1 \times D$ is given by the version of (2.1) that has $(\nu, \mu)$ replaced by $(\nu_{\tau=k/Q}, \mu_{\tau=k/Q})$.

The construction of $a_{k+1}$ also requires a smooth, increasing function, $\chi: [0, \infty) \to [0, 1]$ with value 1 on $[0, \frac{5}{16}]$ and value 0 on $[\frac{7}{16}, \infty)$. This function should be fixed once and for all. Given $\rho > 0$ with $\rho \ll \rho_k$, define a function $\tau_\rho$ on $D$ to equal

$$\tau_\rho = \tfrac{k}{Q} + \tfrac{1}{Q}\chi(\tfrac{1}{\rho}|z|).$$

(A.5)

Thus, $\tau_\rho = \tfrac{k+1}{Q}$ where $|z| \leq \tfrac{1}{4}\rho$ and $\tau_\rho = \tfrac{k}{Q}$ where $|z| > \rho$. Note as well that

$$|d\tau_\rho| \leq c_0 \tfrac{1}{Q}\rho^{-1} \quad \text{and} \quad |\nabla d\tau_\rho| \leq c_0 \tfrac{1}{Q}\rho^{-2}.$$

(A.6)

With $\rho$ chosen as above, define the 1-form on $S^1 \times D^2$ by the formula

$$\tfrac{2\pi}{\ell_\gamma} a^\rho = (1 - 2\nu_{\tau_\rho}|z|^2 - \mu_{\tau_\rho}\bar{z}^2 - \bar{\mu}_{\tau_\rho} z^2)\, dt + \tfrac{i}{2}(z\, d\bar{z} - \bar{z}\, dz) + (1-\chi(\tfrac{1}{\rho}|z|))(\cdots)$$

(A.7)

where the terms that are indicated by the three dots on the right are identical to those that appear in (2.1). It is a consequence of (A.6) that what is written in (A.7) defines a contact 1-form when ρ is sufficiently small. To see this, remark that $a^\rho = a_k$ where $|z| > \rho$, and where $|z| \le \rho$,

- $\frac{2\pi}{\ell_\gamma} a^\rho = dt + \frac{i}{2}(z\, d\bar{z} - \bar{z}\, dz) + \mathfrak{r}_0$
- $\frac{2\pi}{\ell_\gamma} da^\rho = i\, dz \wedge d\bar{z} - 2(\nu_{\tau_\rho} z + \mu_{\tau_\rho} \bar{z})\, d\bar{z} \wedge dt - 2(\nu_{\tau_\rho} \bar{z} + \bar{\mu}_{\tau_\rho} z)\, dz \wedge dt + \mathfrak{r}_1$

(A.8)

where $|\mathfrak{r}_0| \le c_0 |z|^2$ and $|\mathfrak{r}_1| \le c_0 \frac{1}{Q} |z|$.

Define $a_{k+1,\rho}$ as follows: If $\gamma \in \mathcal{R}_L$, set $a_{k+1,\rho}$ to equal $\gamma$'s version of $a^\rho$ on the image of $\varphi_\gamma$. Meanwhile, set $a_{k+1,\rho} = a_k$ on the complement of the union of these tubular neighborhoods.

**Lemma A.1**: *There exists $\kappa > 1$ with the following significance: If $Q > \kappa$ and if $\rho$ is sufficiently small, then $a_{k,\rho}$ satisfies the first and second items in (2.11) plus (A.1).*

*Proof of Lemma A.1*: The proof has four steps.

Step 1: Let $v_k$ denote the Reeb vector field for the contact form $a_k$ and let $v_{k+1,\rho}$ denote the Reeb vector field for $a_{k+1,\rho}$. The latter agrees with $v_k$ except near a Reeb orbit from $\mathcal{R}_L$. Let $\gamma$ denote such an orbit. As before, use $\varphi_\gamma$ to view a neighborhood of $\gamma$ as a neighborhood of $S^1 \times \{0\}$ in $S^1 \times D$. This done, then $v_{k+1,\rho}$ agrees with $v_k$ except on the part of $S^1 \times D$ where $|z| < \rho$. Meanwhile, it follows from the second item in (2.19) that

$$\frac{\ell_\gamma}{2\pi} v_{k+1,\rho} = \frac{\partial}{\partial t} + 2i(\nu_{\tau_\rho} z + \mu_{\tau_\rho} \bar{z})\frac{\partial}{\partial z} - 2i(\nu_{\tau_\rho} \bar{z} + \bar{\mu}_{\tau_\rho} z)\frac{\partial}{\partial \bar{z}} + \mathfrak{v}$$

(A.9)

where

$$|\mathfrak{v}| \le c_0 \frac{1}{Q} |z| \quad and \quad |\nabla \mathfrak{v}| \le c_0.$$

(A.10)

The formula for $v_{k+1,\rho}$ indicates first that $\gamma$ is also an $v_{k+1,\rho}$ Reeb orbit, and that it is elliptic or hyperbolic as an $v_{k+1,\rho}$ Reeb orbit for all small $\rho$. It also indicates that $\gamma$'s rotation number as an $v_{k+1,\rho}$ Reeb orbit is independent of $\rho$ and identical to its rotation number as $v_k$ and $v$ orbit.

Step 2: This step proves that every $v_{k,\rho}$ Reeb orbit with symplectic action less than L lies in a tubular neighborhood of some Reeb orbit from $\mathcal{R}_L$. To this end, suppose that there exists a sequence $\{\rho_j\}_{j=0,1,...}$ with limit zero and a corresponding sequence $\{\gamma_j\}$ of Reeb orbits for the $\rho = \rho_j$ version of $v_{k+1,\rho}$, all with symplectic action as defined by the $\rho = \rho_j$ version of $a_{k+1,\rho}$ bounded by L. View these loops as images of maps from $S^1$ into M. The bounds in (A.10) on $v_{k+1,\rho}$ and its first derivative guarantee the existence of a convergent subsequence in $C^1(S^1; M)$ whose limit map has the following property: It's image is a closed integral curve of $v_k$ with symplectic action less than L. Thus, each large j version of $\gamma_j$ must lie in the image of $\varphi_\gamma$ that is associated to some Reeb orbit $\gamma \in \mathcal{R}_L$.

The following is a direct consequence: There exists $\rho_0$ such that if $\rho$ is less than $\rho_0$ and if $\gamma'$ is an $v_{k+1,\rho}$ Reeb orbit with symplectic action less than L, then $\gamma'$ lies in the image of the tubular neighborhood map $\varphi_\gamma$ that is associated to some Reeb orbit $\gamma \in \mathcal{R}_L$.

<u>Step 3</u>: A virtual repeat of what is said in Step 2 strengthens Step 2's conclusions as follows: Given $\sigma > 0$, there exists $\rho_\sigma$ such that if $\rho < \rho_\sigma$ and if $\gamma'$ is an $v_{k+1,\rho}$ Reeb orbit with symplectic action less than L, then $\gamma'$ lies in the image of the tubular neighborhood map $\varphi_\gamma$ that is associated to a Reeb orbit $\gamma \in \mathcal{R}_L$. Moreover, if $\gamma'$ is such an $v_{k+1,\rho}$ Reeb orbit, and if it lies in $\varphi_\gamma(S^1 \times D)$, then the coordinate z for D obeys $|z| \le \sigma$ on $\gamma'$.

<u>Step 4</u>: Let $\gamma \in \mathcal{R}_L$ and let $\gamma'$ denote a $v_{k+1,\rho}$ Reeb orbit with symplectic action less than L that lies in $\varphi_\gamma(S^1 \times D)$. Suppose for the sake of argument that $\gamma' \ne \gamma$. Let q denote the winding number of $\gamma'$ in $S^1 \times D$. It follows from (A.9) that $\gamma'$ can be viewed as a $2\pi q$ periodic map from $\mathbb{R}$ to $S^1 \times D$ by parametrizing it so that the pull-back of dt is the Euclidean 1-form on $\mathbb{R}$. This done, use $z' : \mathbb{R}/(2\pi q\mathbb{Z}) \to \mathbb{C}$ to denote the function that sends $t \to z(\gamma'(t)) \in \mathbb{C}$.

Let $\hat{v} : \mathbb{R}/(2\pi q\mathbb{Z}) \to \mathbb{R}$ denote the map whose value at any given point $t \in \mathbb{R}/(2\pi q\mathbb{Z})$ is that of $v_{\tau_\rho}$ at $(t, z'(t))$. Define $\hat{\mu}$ in an similar fashion. It then follows that

$$|\hat{v} - v_{\tau=k/Q}| + |\hat{\mu} - \mu_{\tau=k/Q}| \le c_0 \tfrac{1}{Q}.$$

(A.11)

Here, $c_0$ depends solely on the first derivative bounds that are alluded to in (A.2).

The preceding inequality implies that the function $z'$ is a $2\pi q$-periodic solution to an equation that has the schematic form

$$\tfrac{i}{2} \tfrac{d}{dt} z' + v_{\tau=k/Q} z' + \mu_{\tau=k/Q} \overline{z}' = \mathfrak{r} \quad where \quad |\mathfrak{r}| \le c_0(\tfrac{1}{Q} |z'| + |z'|^2).$$

(A.12)

<u>Step 5</u>: Let $\lambda_0 > 0$ again denote the eigenvalue bound that is alluded to in (A.3). Let $c_0$ denote the specific constant that appears in (A.12). Then (A.12) requires $z' = 0$ if Q is chosen so that $Q^{-1} \le (100 c_0)^{-1} \lambda_0$. With this choice of Q, each Reeb orbit for any $\rho << \rho_k$ version of $v_{k+1,\rho}$ with symplectic action less than L must be a loop from $\mathcal{R}_L$.

Given Lemma A.1, each $\rho \ll \rho_k$ version of $a_{k+1,\rho}$ is a candidate for $a_{k+1}$.

*Part 2*: This part of the proof defines an almost complex structure that is compatible with each such small $\rho$ version of $a_{k+1,\rho}$. This almost complex structure is denoted in what follows by $J_{k+1,\rho}$.

To start, set $J_{k+1,\rho}$ equal to $J_k$ on the complement of images of the tubular neighborhood maps for the Reeb orbits in $\mathcal{R}_L$. Now let $\gamma$ denote such a Reeb orbit, and let $\varphi_\gamma$ denote its tubular neighborhood map. As before, use $\varphi_\gamma$ to identify $S^1 \times D$ with a neighborhood of $\gamma$. Since $J_{k+1,\rho} \frac{\partial}{\partial s} = v_{k+1,\rho}$, the only ambiguity concerns the action of $J_{k+1,\rho}$ on the kernel of $a_{k+1,\rho}$. A look at (A.1) and (A.7) indicates that $J_{k+1,\rho}$ can be chosen so as to have the following properties: First, $J_{k+1,\rho} = J_k$ except at points where $|z| \leq \rho$. Second,

$$|J_{k+1,\rho} - J_k| \leq c_0 \tfrac{1}{Q} |z| \quad \text{with} \quad |\nabla(J_{k+1,\rho} - J_k)| \leq c_0 \tfrac{1}{Q} \quad \text{and} \quad |\nabla^2(J_{k+1,\rho} - J_k)| \leq c_0 \tfrac{1}{Q} \rho^{-1}.$$
(A.13)

Third, the $\varphi$ pull-back of the $J_{k+1,\rho}$ version of $T^{1,0}(\mathbb{R} \times M)$ is spanned by $ds + i\varphi^* a_{k+1,\rho}$ and by $\tfrac{\ell_\gamma}{2\pi}(dz - 2i(v_{\tau_\rho} z + \mu_{\tau_\rho} \bar{z})dt)$ at points where $|z| \leq \tfrac{1}{4}\rho$. Fix $J_{k+1,\rho}$ with these properties.

*Part 3*: The next task is to construct a 1-1 map from the set of components of the $J_k$ version of $\mathcal{M}_1(\Theta_-, \Theta_+)$ to those of the $J_{k+1,\rho}$ version. To start this task, from $\Theta_-$ and $\Theta_+$ $\mathcal{C}_*^L$ and $\Sigma$ from the $J_k$ version of $\mathcal{M}_1(\Theta_-, \Theta_+)$. Then each $\mathbb{R}$-invariant cylinder from $\Sigma$ is $J_{k+1,\rho}$ pseudoholomorphic, this because each has the form $\mathbb{R} \times \gamma$ with $\gamma \in \mathcal{R}_L$. Let $C \subset \Sigma$ denote the component that is not $\mathbb{R}$-invariant. Then $C$ is $J_{k+1,\rho}$ pseudoholomorphic except where it intersects the product of $\mathbb{R}$ with a tubular neighborhood of a Reeb orbit in $\mathcal{R}_L$. To say more on this, let $\gamma$ denote such a Reeb orbit. Again use $\varphi_\gamma$ to identify $S^1 \times D$ with its $\varphi_\gamma$ image. Given $R > 1$, there exists $\rho_R \ll \rho_k$ such that if $\rho < \rho_R$, then the intersection of $C$ with the $|z| < \rho$ part of $\mathbb{R} \times (S^1 \times D)$ can occur only in the following two ways:

- *Intersection occurs in a disk of radius $R^{-1}$ in C centered around each of point from the finite set where C intersects $\mathbb{R} \times \gamma$.*
- *Intersection can occur where $|s| > R$ on those ends of C that are labeled by a pair from $\Theta_- \cup \Theta_+$ whose Reeb orbit component is $\gamma$.*

(A.14)

Let $\chi_{C,\rho}$ denote the characteristic function of the support of $|J_{k+1,\rho} - J_k|$ on $C$.

In any event, (A.13) implies that $C$ is nearly $J_{k+1,\rho}$ pseudoholomorphic in that each of its tangent planes is nearly $J_{k+1,\rho}$ invariant. Moreover, $C$ is nearly $J_{k+1,\rho}$ pseudoholomorphic in an $L^2$ sense. To elaborate, reintroduce the normal bundle $N \to C$ and let π: T(ℝ × M)|_C → N denote the orthogonal projection. It then follows from (2.6), (A.13) and (A.14) that

$$\int_C |\pi \circ J_{k+1,\rho}|^2 \leq c_0 \frac{1}{Q^2} \rho^2.$$

(A.15)

The relatively small $L^\infty$ and $L^2$ norms of $\pi \circ J_{k+1,\rho}$ suggest a perturbative construction of a 1-1 map from the set of components of the $J_k$ version of $\mathcal{M}_1(\Theta_-, \Theta_+)$ to the set of components of the $J_{k+1,\rho}$ version that pairs components so as to satisfy the fourth item in (2.11). Such a construction is given in the three steps that follow.

  Step 1: This step sets up this perturbative construction. To start, fix a component of the $J_k$ version of $\mathcal{M}_1(\Theta_-, \Theta_+)$ and fix a point Σ in this component. A partner for Σ in the $J_{k+1,\rho}$ version of $\mathcal{M}_1(\Theta_-, \Theta_+)$ is described next. This partner has the same set of ℝ-invariant cylinders with the same integer weights as does Σ. Let C ⊂ Σ denote the component that is not ℝ-invariant. The analogous component of the partner to Σ is constructed as a deformation of C that comes via a section of C's normal bundle N by composing the section with a suitably chosen exponential map from a disk subbundle in N to ℝ × M.

  To say more, suppose that $N_1 \subset N$ is a constant radius disk subbundle and suppose that $e_C : N_1 \to \mathbb{R} \times M$ is an exponential map that embeds each fiber disk as a pseudoholomorphic disk. Such maps are constructed in Section 5d of *SW => Gr* from [T2]. Note that $e_C$ can not embed the whole of $N_1$ unless each pair in $\Theta_- \cup \Theta_+$ has its second component equal to 1. In any event, given $e_C$, let η denote a section of $N_1$ over Σ which has $|s| \to \infty$ limit equal to zero. Then $e_C \circ \eta$ is $J_{k+1,\rho}$-pseudoholomorphic if and only if it obeys an equation that has the schematic form

$$\mathcal{D}_C \eta + \mathfrak{p}_1 \cdot \eta + (\mathfrak{R}_1(\eta) + \mathfrak{p}_2) \cdot \nabla_C \eta + \mathcal{R}_0(\eta) + \mathfrak{p}_0 = 0$$

(A.16)

To elaborate for a moment on the notation, $\mathcal{D}_C$ denotes the operator that appears in (2.8). Meanwhile, $\mathfrak{p}_1$ is a zero'th order, ℝ-linear operator that obeys $|\mathfrak{p}_1| \leq c_0 \frac{1}{Q}$ and with support in the two regions that are listed in (A.14). What is called $\mathfrak{p}_2$ in (A.16) is a homomorphism with support where $J_{k+1,\rho} \neq J_k$. It has norm $|\mathfrak{p}_2| \leq c_0 \frac{1}{Q} \rho$, this because it is bounded by $c_0 |J_{k+1,\rho} - J_k|$. What is denoted by $\mathfrak{p}_0$ is obtained from $\pi \circ J_{k+1,\rho}$ by restricting the latter to the (0, 1) tangent space of C. Finally, $\mathfrak{R}_1$ denotes a fiber preserving map from $N_1$ to Hom($T^{1,0}C, T^{0,1}C$) and $\mathfrak{R}_2$ denotes a fiber preserving map from $N_1$ to $N \otimes T^{0,1}C$. By virtue of (2.6) and (A.13), these two maps obey

- $|\mathfrak{R}_1(b)| \leq c_0 \frac{1}{Q} |b|$ *and* $|\mathfrak{R}_0(b)| \leq c_0 |b|^2$.

- $|\nabla \mathfrak{R}_1| \le c_0 \frac{1}{Q}$   *and*   $|\nabla \mathfrak{R}_0| \le c_0 \frac{1}{Q} |b|$ .

(A.17)

Step 2: Granted what was just said, a contraction mapping argument can be used to find small normed solutions to (A.10) when ρ is small given that the linear term in (A.16) is invertible as a map between suitable Banach spaces, and given that $\mathfrak{p}_0$ has suitably small norm as an element in the range Banach space.

**Lemma A.2**: *There exists* κ > 1 *with the following significance: Suppose that* Q ≥ κ, *that* k ∈ {1, …, Q} *and that* $\{(a_i, J_i)\}_{1 \le i \le k}$ *has been constructed. If* ρ *is sufficiently small, then* $\mathcal{D}_C + \mathfrak{p}_2 \cdot \nabla_C + \mathfrak{p}_1$ *has bounded inverse as a map from the* $L^2$-*orthogonal complement of the kernel of* $\mathcal{D}_C$ *in* $L^2_1(C; N)$ *to* $L^2(C; N \otimes T^{0,1}C)$. *Moreover, such is also the case for the operator* $\mathcal{D}_C + \sigma(\mathfrak{p}_2 \cdot \nabla_C + \mathfrak{p}_1)$ *for each* σ ∈ [0, 1].

Note that the σ ≠ 1 version of the lemma is needed only to compare respective ±1 weights that are used to define the embedded contact homology differential. This lemma is proved below in Part 7. Assume it for now.

The $L^2_1$ norm does not dominate the $L^\infty$ norm, and this makes $L^2_1$ unsuitable as the Banach space for the contraction mapping. However, a slightly stronger norm can be used to define a suitable Banach space. To say more, introduce a norm on the space of compactly supported sections of either N or $N \otimes T_\mathbb{C} C$ as follows: Its square assigns to a section, ζ, the number

$$\int_C |\zeta|^2 + \sup_{z \in C} \sup_{x \in (0,1)} x^{-1/100} \int_{\text{dist}(z,\cdot) < x} |\zeta|^2$$

(A.18)

The Banach space for the contraction mapping argument is the completion of the space of compactly supported sections of N using the norm whose square assigns to any given section η the sum of three terms. The first is the square of the $L^2_1$ norm, and the next two are the respective ζ = η and ζ = ∇η versions of (A.18). This space is denoted by $\mathcal{B}_1$ and its norm is denoted by $\|\cdot\|_*$. An appeal to Theorem 3.5.2 in [M] finds a constant, $c_C$, such that $|\cdot| \le c_C \|\cdot\|_*$. Note that this constant $c_C$ depends on the curve C. Let $\mathcal{B}_1^\perp \subset \mathcal{B}_1$ denote the subspace of elements that are $L^2$-orthogonal to the kernel of $\mathcal{D}_C$.

Use $\mathcal{B}_0$ to denote the completion of the space of compactly supported sections of $N \otimes T^{0,1}C$ using the norm whose square is is depicted in (A.18). If $\mathcal{D}_C + \mathfrak{p}_2 \cdot \nabla_C + \mathfrak{p}_1$ has a bounded inverse mapping $L^2(C; N \otimes T^{0,1}C)$ to $L^2_1(N)$, then an argument using Theorem 5.4.1 of [M] finds that the inverse of $\mathcal{D}_C + \mathfrak{p}_2 \cdot \nabla_C + \mathfrak{p}_1$ restricts to $\mathcal{B}_0$ so as to define a bounded, linear map from $\mathcal{B}_0$ to $\mathcal{B}_1^\perp$.

With the preceding understood, fix $\sigma > 0$ so that elements of $\mathcal{B}_1$ with $\|\cdot\|_* < 2\sigma$ define sections of the disk bundle $N_1$. Let $\mathcal{U}_\sigma \subset \mathcal{B}_1^\perp$ denote the ball of radius $\sigma$ centered on the origin and define the map $\mathcal{T}: \mathcal{U}_\sigma \to \mathcal{B}_1^\perp$ by setting

$$\mathcal{T}(\eta) = - (\mathcal{D}_C + \mathfrak{p}_2 \cdot \nabla_C + \mathfrak{p}_1)^{-1}(\mathfrak{R}_1(\eta) \cdot \nabla_C \eta + \mathcal{R}_0(\eta) + \mathfrak{p}_0) \,.$$

(A.19)

**Lemma A.3**: *There exists $\sigma' \in (0, \sigma)$ such that if $\rho$ is sufficiently small then the following is true: Suppose that $C$ is a component of a submanifold in the $J_k$ version of $\mathcal{M}_1(\Theta_-, \Theta_+)$. Then $\mathcal{T}$ defines a contraction mapping from $\mathcal{U}_{\sigma'}$ to itself. For such $\rho$, the map $\mathcal{T}$ has a unique fixed point in $\mathcal{U}_{\sigma'}$. Moreover, this fixed point has $\|\cdot\|_*$ norm bounded by $c_C \rho$, and it is a smooth section of $N$ that obeys (A.16). Here, $c_C$ is independent of $\rho$ but depends on $C$. In the case that $C = \mathbb{R} \times \gamma$, this fixed point is $\eta = 0$.*

*Proof of Lemma A.3*: It follows from (A.13) and the first line in (A.17) that

$$\|\mathcal{T}(\eta)\|_* \leq c_{C1}(\|\eta\|_*^2 + \rho) \,,$$

(A.20)

where $c_{C1}$ is a constant that is independent of $\rho$ but dependent on $C$. This last bound implies that $\mathcal{T}$ maps the ball in $\mathcal{B}_1^\perp$ of radius $\frac{1}{4} c_{C1}^{-1}$ to itself when $\rho < \frac{1}{8} c_{C1}^{-2}$. Meanwhile, the second line in (A.17) implies that

$$\|\mathcal{T}(\eta) - \mathcal{T}(\eta')\|_* \leq c_{C2}(\|\eta\|_* + \|\eta'\|_*)\|\eta - \eta'\|_*$$

(A.21)

where $c_{C2}$ is a second $C$ dependent but $\rho$ independent constant. This last bound implies that $\mathcal{T}$ maps the ball of radius $\sigma'$ to itself as a contraction mapping if $\sigma' < \frac{1}{4}(c_{C1} + c_{C2})^{-1}$ and if $\rho$ is sufficiently small. The remaining assertions of the lemma follow using standard elliptic regularity arguments as can be found in Chapter 5 of [M].

Step 3: Let $C$ now denote the non $\mathbb{R}$-invariant component of $\Sigma$. With $\rho$ very small, let $C' \subset \mathbb{R} \times M$ denote the immersed subvariety that is obtained from $C$ using the section $\eta$ given by Lemma A.3. This subvariety is $J_{k+1,\rho}$ pseudoholomorphic by construction. Introduce $\Sigma'$ to denote the union of $C'$ and the other $\mathbb{R}$-invariant elements in $\Sigma$ with their associated integer weights. As a parenthetical remark, note that Proposition 11.4 in [HS] implies the following: The union of the subvarieties that comprise $\Sigma'$ is an embedded subvariety in $\mathbb{R} \times M$. In any event, $\Sigma'$ defines an element in the $J_{k+1,\rho}$ version of $\mathcal{M}_1(\Theta_-, \Theta_+)$.

The association of Σ to Σ´ is an injective map from the set of components of the $J_k$ version of $\mathcal{M}_1(\Theta_-, \Theta_+)$ into the set of components of the $J_{k+1,\rho}$ version of $\mathcal{M}_1(\Theta_-, \Theta_+)$. This map is denoted in what follows by $\mathcal{F}$.

*Part 4*: This part verifies that if ρ is sufficiently small, then the components of the $J_{k+1,\rho}$ version of $\mathcal{M}_1(\Theta_-, \Theta_+)$ that lie in the image of $\mathcal{F}$ are smooth points in this version of $\mathcal{M}_1(\Theta_-, \Theta_+)$ in the sense of Lemma 2.2. It also verifies that the correspondence that is defined by $\mathcal{F}$ satisfies the the third item in (2.11) if $\hat{J}$ is replaced there by $J_{k+1,\rho}$.

**Lemma A.4**: *Suppose that ρ is very small. Let Σ´ denote a subvariety from a component of the $J_{k+1,\rho}$ version of $\mathcal{M}_1(\Theta_-, \Theta_+)$ that is in the image of the map $\mathcal{F}$. Let C´ denote a component of Σ´. The associated deformation operator $\mathcal{D}_{C´}$ has trivial cokernel. Thus, the component of Σ´ in the $J_{k+1,\rho}$ version of $\mathcal{M}_1(\Theta_-, \Theta_+)$ is an orbit of the $\mathbb{R}$-action on this space. Moreover, the sign that this component would contribute to the embedded contact homology differential is the same as that of its $\mathcal{F}$-inverse image in the $J_k$ version of $\mathcal{M}_1(\Theta_-, \Theta_+)$.*

*Proof of Lemma A.4*: The invertibility of $\mathcal{D}_{C´}$ when C´ is an $\mathbb{R}$-invariant cylinder from Σ´ is automatic since this operator doesn't change when $(a_k, J_k)$ is replaced by $(a_{k+1,\rho}, J_{k+1,\rho})$. Let C´ denote the non $\mathbb{R}$-invariant component of Σ´. Let Σ denote the subvariety from the $J_k$ version of $\mathcal{M}_1(\Theta_-, \Theta_+)$ that gives rise to Σ´, and let $C \subset \Sigma$ denote the non $\mathbb{R}$-invariant component that is used to construct C´ via Lemma A.3. Since C´ is the image via the exponential map of a section of C's normal bundle, it follows that C´ is immersed and so has a normal bundle, N´ → C´. Note for reference momentarily that Lemma A.2 asserts that $\mathcal{D}_C + \sigma(\mathfrak{p}_2 \cdot \nabla_C + \mathfrak{p}_1)$ is invertible for any constant $\sigma \in [0, 1]$.

Let π: N → C again denote the normal bundle to C. Use the exponential map $e_C$ to view C´ as the graph in N of the section η. This identifies the normal bundle of C´ with the restriction of the bundle π*N to this graph. The view of C´ as the graph of η also supplies an $\mathbb{R}$-linear isomorphism between $T^{0,1}C$ and $T^{0,1}C´$. These identifications allow $\mathcal{D}_{C´}$ to be viewed as a bounded operator from $L^2_1(C; N)$ to $L^2(C; N \otimes T^{0,1}C)$. As such, it has the form $\mathcal{D}_C + \mathfrak{p}_2 \cdot \nabla_C + \mathfrak{p}_1 + \mathfrak{r}$ where $\mathfrak{r}$ has operator norm as a map from $\mathcal{B}_1$ to $\mathcal{B}_0$ that is bounded by $c_C \varepsilon(\rho)$ where $\rho \to \varepsilon(\rho)$ is a decreasing function with limit 0 as $\rho \to 0$. Meanwhile, $c_C$ is independent of ρ but not C. Indeed, the latter fact follows because the derivative terms in $\mathfrak{r}$ have coefficients bounded by $c_C|\eta|$ and the zero'th order terms in $\mathfrak{r}$ have coefficients bounded by $c_C|\nabla\eta|$. The operator norms of these terms can be bounded by $c_C\|\eta\|_*$ using, respectively, Theorem 3.5.2 and Lemma 5.4.1 in [M].

Granted that $\mathcal{D}_C + \mathfrak{p}_2 \cdot \nabla_C + \mathfrak{p}_1$ is invertible, it follows from what was just said about $\mathfrak{r}$ that any sufficiently small ρ version of $\mathcal{D}_{C´}$ has trivial cokernel when viewed as a map from $L^2_1(C'; N')$ to $L^2(C'; N' \otimes T^{0,1}C')$. The fact that $\mathcal{D}_{C'}$ has trivial cokernel implies that the $\Sigma'$ is a smooth point of the $J_{k+1,\rho}$ version of $\mathcal{M}_1(\Theta_-, \Theta_+)$ in the sense of Lemma 2.2. Thus, its component is isomorphic to $\mathbb{R}$ with tangent vector field the generator of the $\mathbb{R}$ action that is induced by the constant translations on the $\mathbb{R}$ factor of $\mathbb{R} \times M$.

The small norm of $\mathcal{D}_{C'} - (\mathcal{D}_C + \mathfrak{p}_2 \cdot \nabla_C + \mathfrak{p}_1)$ and the fact that $(\mathcal{D}_C + \sigma(\mathfrak{p}_2 \cdot \nabla_C + \mathfrak{p}_1))$ is invertible for all $\sigma \in [0, 1]$ imply that the sign contribution of the component of $\Sigma'$ to the embedded contact homology differential is the same as that of the component of $\Sigma$. Here is why: These signs are defined using the determinant line bundles for the operators $\mathcal{D}_{C'}$ and $\mathcal{D}_C$. The linear interpolation between $\mathcal{D}_{C'}$ and $\mathcal{D}_C + \mathfrak{p}_2 \cdot \nabla_C + \mathfrak{p}_1$ provides a canonical isomorphism between the respective determinant lines, as does the linear interpolation between $\mathcal{D}_C + \mathfrak{p}_2 \cdot \nabla_C + \mathfrak{p}_1$ and $\mathcal{D}_C$. Meanwhile, the orientations that are defined by the $\mathbb{R}$ actions on the components of the respective $J_k$ and $J_{k+1,\rho}$ versions of $\mathcal{M}_1(\Theta_-, \Theta_+)$ are compatible as the construction of $C'$ from $C$ is $\mathbb{R}$-equivariant.

*Part 5*: Make a very small perturbation $a_{k+1,\rho}$ and $J_{k+1,\rho}$ so that the resulting contact structure is in Lemma 2.1's residual set and so that the complex structure comes from the set $\mathcal{J}_{a_{k+1,\rho}}$. The perturbation of $a_{k+1,\rho}$ can be as small as desired as measured with respect to any $q \gg 1$ version of the $C^q$ norm, but it should have support away from the Reeb orbits in $\mathcal{R}_L$, and also away from the pseudoholomorphic curves that are supplied by the map $\mathcal{F}$. The perturbation of $J_{k+1,\rho}$ can also be as small as desired as measured with respect to any given $q \gg 1$ version of the $C^q$ norm. In addition, it can and should be made so that the new version of the almost complex structure agrees with the original near the cylinders that comprise $\mathbb{R} \times \mathcal{R}_L$ and on all pseudoholomorphic curves that appear in elements from $\mathcal{F}$'s image. Moreover, the 2-jets of the new and original versions should also agree on the latter. Finally, the new and old versions should differ by an endomorphism with $C^3$ norm less than $\rho^3$. Note in this regard that the projection to M of the union of the curves that appear in elements from $\mathcal{F}$'s image defines a codimension 1 subvariety in M. The arguments to justify that such perturbations exist are very much like those used in Section 4 of [HT2] and will not be presented.

Agree to use $(a_{k+1,\rho}, J_{k+1,\rho})$ henceforth to denote this slightly perturbed version of the almost complex structure given in Part 2.

*Part 6:* This part proves that the small $\rho$ versions of the map $\mathcal{F}$ are onto. To start the proof, suppose, to the contrary that such is not the case so as to derive some nonsense. Under this assumption, there exists a decreasing sequence $\{\rho_\upsilon\}_{\upsilon=1,2,\ldots}$ with limit zero, and a for each index $\upsilon$, there would exist a point $\Sigma_\upsilon$ in the $J_{k+1,\rho_\upsilon}$ version of $\mathcal{M}_1(\Theta_-, \Theta_+)$ that does not arise as described from a a point in the $J_k$ version of $\mathcal{M}_1(\Theta_-, \Theta_+)$. What are now standard compactness arguments can be used to prove that the sequence $\{\Sigma_\upsilon\}$ has a subsequence that converges to what is often called a *broken trajectory*. Indeed, the limiting behavior is defined by a finite, ordered set $\Lambda = \{\Sigma_1, \ldots, \Sigma_p\}$ where any given $\Sigma_j$ consists of a finite set of pairs of the form (S, m) with $S \subset \mathbb{R} \times M$ an irreducible, $J_k$ pseudoholomorphic subvariety and m a positive integer. Moreover, there are constraints on the $|s| \to \infty$ limits of the pairs that comprise any given $\Sigma_j$. The digression that follows is needed to describe these constraints.

To start the digression, suppose that $(S, m) \in \Sigma_j$. The large $|s|$ slices of any end E $\subset$ S converge as $|s| \to \infty$ as a multiple cover of some Reeb orbit. If $\gamma$ is such a limit Reeb orbit, define $m_{\gamma,S_-}$ to denote the multiplicity of this covering. Set $m_{\gamma,S_-} = 0$ if $\gamma$ is not multiply covered by the $s \to -\infty$ limit of the constant s slices of any negative end of S. Now associate to $\Lambda_j$ the set $\Theta_{j-}$ whose elements are pairs of the form $(\gamma, q)$ where $\gamma$ is a Reeb orbit that is multiply covered by the $s \to -\infty$ limit of the constant s slices of some end in $\cup_{(S,m) \in \Lambda_j} S$, and where $q = \sum_{(S,m) \in \Lambda_j} m\, m_{\gamma,S_-}$. Likewise define $\Theta_{j+}$.

What follows are the constraints on the pairs that comprise the elements from $\Lambda$:

$$\Theta_{1-} = \Theta_-, \quad \Theta_{p+} = \Theta_+ \text{ and } \Theta_{j+} = \Theta_{j+1-} \text{ for each } j \in \{1, \ldots, p-1\}.$$
(A.22)

Note that these constraint imply that each $\Theta_{j\pm} \in \mathcal{C}_*^L$. Given Equation (102) in [HS] and Propositions 11.4 and Corollary 11.5 in [HS], it follows that $\Lambda$ has just one element and that this element is in the $J_k$ version $\mathcal{M}_1(\Theta_-, \Theta_+)$. The argument here is essentially the same as that used to Theorem 1.8 in [Hu]. See also the proof of Lemma 7.19 in [HT1]. Keep in mind that this limit element in the $J_k$ version of $\mathcal{M}_1(\Theta_-, \Theta_+)$ consists of a set of disjoint cylinders with weights, and one non-$\mathbb{R}$ invariant submanifold that is disjoint from the cylinders. Let $\Sigma$ denote the element in question.

The manner of convergence of the subsequence of $\{\Sigma_\upsilon\}_{\upsilon=1,2,\ldots}$ to $\Sigma$ is described next There is a sequence $\{s_\upsilon\}_{\upsilon=1,2,\ldots} \in \mathbb{R}$ such that the translation $s \to s_\upsilon$ of $\Sigma_\upsilon$ along the $\mathbb{R}$ factor of $\mathbb{R} \times M$ gives a new subsequence, now renamed $\{\Sigma_\upsilon\}_{\upsilon=1,2,\ldots}$, such that

- $\lim_{\upsilon \to \infty} (\sup_{z \in (\cup_{(C,m) \in \Sigma_\upsilon} C)} \text{dist}(z, \Sigma) + \sup_{z \in (\cup_{(C,m) \in \Sigma} C)} \text{dist}(z, \Sigma_\upsilon)) = 0$.
- $\lim_{\upsilon \to \infty} \sum_{(C,m) \in \Sigma_\upsilon} m \int_C w = \sum_{(C,m) \in \Sigma} m \int_C w$ *for any 2-form* w *on* $\mathbb{R} \times M$ *with compact support*.

(A.23)

The sort of convergence that is dictated in (A.23) requires the existence of a 1-1 correspondence between the set of components of any sufficiently large $\upsilon$ version of $\Sigma_\upsilon$ and the set of components of $\Sigma$. This correspondence is such that if (S, m) is a component of $\Sigma$ and $(S_\upsilon, m_\upsilon) \subset \Sigma_\upsilon$, then the convergence in (A.23) holds with $(S_\upsilon, m_\upsilon)$ replacing $\Sigma_\upsilon$ and (S, m) replacing $\Sigma$.

Granted this last point, the contraction mapping theorem argument from Lemma A.3 proves that the partner from any sufficiently large $\upsilon$ version of $\Sigma_\upsilon$ to an $\mathbb{R}$-invariant cylinder in $\Sigma$ must coincide with this cylinder and its respective $\Sigma_\upsilon$ and $\Sigma$ integer weights must agree. Here is why: Let C denote a pseudoholomorphic subvariety for any given almost complex structure compatible to any given contact 1-form. Suppose that C is not $\mathbb{R}$ invariant. View $\pi(C)$ as a cycle in M and write the boundary of this cycle as $\partial_+ C - \partial_- C$, where $\partial_\pm C$ are the respective positive integer weighted sums of Reeb orbits that arise by taking the $s \to \pm\infty$ limits of the constant s slices of C. Then these two weighted sums can not be equal because the integral over C of the exterior derivative of the contact form is strictly positive.

Meanwhile, any partner in a sufficiently large $\upsilon$ version $\Sigma_\upsilon$ to the non $\mathbb{R}$-invariant component of $\Sigma$ must be that given via Lemma A.3 and its contraction mapping. These conclusions contradict the assumptions made at the outset as they imply that the component of each large $\upsilon$ version of $\Sigma_\upsilon$ in the $J_{k+1,\rho_\upsilon}$ version of $\mathcal{M}_1(\Theta_-, \Theta_+)$ is in the image of the 1-1 map from the set of components of the $J_k$ version of $\mathcal{M}_1(\Theta_-, \Theta_+)$.

*Part 7*: This step supplies the

***Proof of Lemma A.2***: The assertion follows by construction when C´ is an $\mathbb{R}$-invariant cylinder. Consider the case when C´ is not $\mathbb{R}$-invariant. Let $\Sigma$ denote the subvariety in the $J_k$ version of $\mathcal{M}_1(\Theta_-, \Theta_+)$ that is paired by $\mathcal{F}$ with $\Sigma'$ and let $C \subset \Sigma$ denote the non $\mathbb{R}$-invariant component that gives rise to C´ via Lemma A.3's contraction mapping.

As noted previously, the operator norm of $\mathfrak{p}_2 \cdot \nabla_C$ is bounded by $c_{C1}\rho$ where $c_{C1}$ depends on C but not on $\rho$. This understood, the Lemma A.2 follows with a proof that $\mathcal{D}_C + \sigma\,\mathfrak{p}_1$ also obeys its conclusions if $\rho$ is small.

To prove the latter assertion, note that $\mathfrak{p}_1$ is non-zero only on the domains that are described by (A.14). The contribution to the operator norm of $\mathfrak{p}_1$ from the disks that are described in the first item of (A.14) is bounded by $c_0 R^{-1}$ since $\mathfrak{p}_1$ has a $\rho$-independent point wise bound and the disks have area $R^{-2}$. Granted that such is the case, the Sobolev theorems in dimension 2 imply the following: Given $\upsilon \in (0, 2)$, there exists a constant $c_0(\upsilon) \geq 1$ such that the contribution of this part of $\mathfrak{p}_1$ to the norm is less than $c_0(\upsilon) R^{-\upsilon}$.

This understood, let $\mathfrak{p}_{1+}$ denote the part of $\mathfrak{p}_1$ with support far out on the ends of C. Granted what was just said about $\mathfrak{p}_1 - \mathfrak{p}_{1*}$, it is sufficient to prove that all sufficiently small $\rho$ versions of $\mathcal{D}_C + \mathfrak{p}_{1+}$ obey the conclusions of Lemma A.4 when $\rho$ is small.

To prove the latter assertion, fix $R \gg 1$ so that the $|s| \geq R$ portion of C is far out on C's ends. Take $\rho$ small so as to guarantee that $\mathfrak{p}_{1+}$ has support only where $|s| > 100R$. Let $u_R$ denote the function on $\mathbb{R}$ that equals 0 where $|s| < R$, equals $(|s|/R - 1)$ where $|s| \in [R, 2R]$ and equals 1 where $|s| > 2R$. Note that this function is Lipschitz. Now write

$$\|(\mathcal{D}_C + \mathfrak{p}_{1+})\eta\|_2^2 = \|u_R^2(\mathcal{D}_C + \mathfrak{p}_{1+})\eta\|_2^2 + \|(1-u_R)^2\mathcal{D}_C \eta\|_2^2$$
(A.24)

Note that $\mathfrak{p}_1$ is absent from the far right term in (A.24) by virtue of the fact that $\mathfrak{p}_1$ is zero where $1 - u_R$ is not. Commute the functions $u_R$ and $(1 - u_R)$ past the derivatives to obtain

$$\|(\mathcal{D}_C + \mathfrak{p}_{1+})\eta\|_2^2 \geq (1 - c_0 R^{-1}) \|(\mathcal{D}_C + \mathfrak{p}_{1+})(u_R\eta)\|_2^2 + \|\mathcal{D}_C((1-u_R)\eta)\|_2^2 - c_0 R^{-1}\|\eta\|_2^2 .$$
(A.25)

The next point is that $u_R\eta$ has support far out on the ends of C. This is to say that each component of the support of $u_R$ in C sits where C is represented as a multi-valued graph over either the $s \geq R$ or $s \leq -R$ part of some $\mathbb{R}$-invariant cylinder; this as depicted in (2.6). To see what this implies, suppose that $\mathcal{E} \subset C$ is an end, and let $\gamma$ denote the associated Reeb orbit. Represent $\mathcal{E}$ as in (2.6) where the eigenfunctions and eigenvalues are those of the operator $L_k$ that is given by replacing $\gamma$'s version of $(\nu, \mu)$ in (2.3) with $(\nu_{\tau=k/Q}, \mu_{\tau=k/Q})$. The operator $\mathcal{D}_C$ on $\mathcal{E}$ differs from $\frac{\partial}{\partial s} + L_k$ by terms that are bounded by $c_0$ $e^{-2\lambda s}$ with $\lambda$ an eigenvalue of $L_k$ that is respectively positive or negative when $\mathcal{E}$ is positive or negative. Since $|\mathfrak{p}_1| \leq c_0 \frac{1}{Q}$, this implies that

$$\|(\mathcal{D}_C + \mathfrak{p}_{1+})(u_R\eta)\|_2^2 \geq (1 - c_0(\tfrac{1}{Q} + R^{-\lambda}))\sigma_*^{-2}\|u_R\eta\|_{2,1}^2 .$$
(A.26)

Here, $\|\cdot\|_{2,1}$ denotes the $L_1^2$ norm. Meanwhile, with $\Pi_C$ denoting the $L^2$ projection orthogonal to the kernel of $\mathcal{D}_C$,

$$\|\mathcal{D}_C((1-u_R)\eta)\|_2^2 \geq \sigma_k^{-2}\|\Pi_C((1-u_R)\eta)\|_{2,1}^2 .$$
(A.27)

Here, $\sigma_k$ is a bound on the inverse of $\mathcal{D}_C$. These last three equations imply that any sufficiently small $\rho$ version of $\mathcal{D}_C + \mathfrak{p}_1$ is invertible as a map from the $L^2$-orthogonal complement in $L_1^2(C; N)$ of the kernel of $\mathcal{D}_C$ to $L^2(C; N \otimes T^{0,1}C)$.

*Part 8*: Given what has been said in the preceding parts, all sufficiently small $\rho$ versions of the pair $(a_{k+1,\rho}, J_{k+1,\rho})$ are such that all but the fourth item of (2.11) are obeyed with $\hat{a} = a_{k+1,\rho}$ and $\hat{J} = J_{k+1,\rho}$. As (A.1) is obeyed, and as $\sigma_{k+1}$ is less than what appears on the right hand side of (A.5), the induction can proceed with $a_{k+1}$ and $J_{k+1}$ set equal to any very small $\rho$ version of $a_{k+1,\rho}, J_{k+1,\rho}$.